\documentclass[pdflatex,sn-mathphys-num]{sn-jnl}

\usepackage{graphicx}%
\usepackage{multirow}%
\usepackage{amsmath,amssymb,amsfonts}%
\usepackage{amsthm}%
\usepackage{mathrsfs}%
\usepackage[title]{appendix}%
\usepackage{xcolor}%
\usepackage{textcomp}%
\usepackage{manyfoot}%
\usepackage{booktabs}%
\usepackage{algorithm}%
\usepackage{algorithmicx}%
\usepackage{algpseudocode}%
\usepackage{listings}%

\usepackage{lineno}
\usepackage{hyperref}
\usepackage{bm}			
\usepackage{indentfirst}
\usepackage{pdflscape}
\usepackage{graphicx}
\usepackage{tabularx,tabulary,multirow,array,hhline,booktabs}
\usepackage{pict2e}
\usepackage{enumerate}
\usepackage{nicematrix}
\usepackage{adjustbox} 
\usepackage{siunitx}
\usepackage{afterpage}
\usepackage{algpseudocode}
\usepackage{enumitem}

\newcolumntype{d}[1]{D{.}{.}{#1}}
\newcolumntype{M}[1]{>{\centering\arraybackslash}m{#1}}
\newcolumntype{N}{@{}m{0pt}@{}}
\DeclareMathOperator*{\argmin}{\arg\min}

\graphicspath{{./figures_arXiv/}{./}}

\makeatletter
\newcommand{\namelabel}[1]{%
\phantomsection
\renewcommand{\@currentlabel}{#1}
\label{#1}
}
\makeatother

\theoremstyle{thmstyleone}%
\theoremstyle{thmstyletwo}%
\newtheorem{remark}{Remark}%

\theoremstyle{thmstylethree}%

\usepackage{ulem}
\newcommand{\pf}[1]{\textcolor{black}{#1}}
\newcommand{\pfb}[1]{\textcolor{black}{#1}}
\definecolor{jhnew}{RGB}{230,120,0}
\newcommand{\nn}[1]{\textcolor{black}{#1}}

\definecolor{jhvio}{RGB}{112,48,160}
\newcommand{\vv}[1]{\textcolor{black}{#1}}

\begin{document}

\title[Eikonal boundary condition for level set method]{Eikonal boundary condition for level set method}


\author*[1]{\fnm{Jooyoung} \sur{Hahn}}\email{jooyoung.hahn@stuba.sk}

\author[2]{\fnm{Karol} \sur{Mikula}}\email{karol.mikula@stuba.sk}

\author[2]{\fnm{Peter} \sur{Frolkovi\v{c}}}\email{peter.frolkovic@stuba.sk}

\affil[1]{\orgdiv{Department of Software Engineering}, \orgname{Faculty of Nuclear Sciences and Physical Engineering, Czech Technical University in Prague}, \orgaddress{\street{Trojanova 339/13}, \city{Prague}, \postcode{120 00}, \country{Czech Republic}}}

\affil[2]{\orgdiv{Department of Mathematics and Descriptive Geometry}, \orgname{Faculty of Civil Engineering, Slovak University of Technology in Bratislava}, \orgaddress{\street{Radlinsk{\'e}ho 11}, \city{Brtislava}, \postcode{810 05}, \country{Slovak Republic}}}


\abstract{In this paper, we propose to use the eikonal equation as a boundary condition when advective or normal flow equations in the level set formulation are solved numerically on polyhedral meshes in the three-dimensional domain. Since the level set method can use a signed distance function as an initial condition, the eikonal equation on the boundary is a suitable choice at the initial time. Enforcing the eikonal equation on the boundary for later times can eliminate the need for inflow boundary \nn{values}, which are typically required for transport equations \nn{and are in general not known}. In selected examples where exact solutions are available, we compare the proposed method with the method using the Dirichlet boundary condition. The numerical results confirm that the use of the eikonal boundary condition provides comparable accuracy and robustness in surface evolution compared to the use of the Dirichlet boundary condition, which is generally not available. We also present numerical results of evolving a general closed surface.}




\maketitle

\section{Introduction}\label{sec:intro}

Curve and surface evolution modeling is a critical area of research with applications in various fields such as \pfb{ multiphase fluid simulation\pfb{~\cite{ref:SSO94,ref:CHMO96,ref:SP00}}, image segmentation\pfb{~\cite{ref:MSV95,ref:CKS97,ref:CV01}}, forest fire propagation\pfb{~\cite{ref:MKF09,ref:MBK11}}, crystal growth\pfb{~\cite{ref:SS92,ref:CMOS97,ref:GFCO03}}, combustion and turbulent flow dynamics~\cite{ref:HASL17,ref:MAOSNL12,ref:SALSP22,ref:KAW88,ref:P99},} and more. These models are used to track the interface dynamics and morphological changes of curves and surfaces over time. Two primary approaches are employed to solve these problems: Lagrangian and Eulerian methods. In the Lagrangian framework, the interface is represented explicitly by geometrical objects that move with the interface. This approach has been extensively studied and applied, yielding significant results in curve and surface evolution. Notable works in this area include the development of finite element methods for \vv{the Willmore flow and related geometric evolution equations}, which provide accurate approximations of interface dynamics by parametrizing the surface~\cite{ref:BGR08}\vv{, and the parametric approximation of mean curvature and anisotropic curvature flows~\cite{ref:DDE05}}. The mathematical foundation for finite element approximation of curve evolution equations~\cite{ref:D90}, including the handling of evolutionary surfaces, has been well established and further extended to evolving surfaces~\cite{ref:DE07} in more complex geometries.

However, in three-dimensional domains, Lagrangian methods face substantial challenges, particularly when dealing with topological changes such as merging or splitting of surfaces. These scenarios are common in practical applications, making the Lagrangian approach less feasible. To address these challenges, Eulerian methods, such as the level set method, have become the preferred choice. The level set method~\cite{ref:OS88} represents the interface implicitly as the zero level set of a higher-dimensional function. This implicit representation allows for natural handling of topological changes, such as merging and splitting, without the need for complex remeshing algorithms.

The level set method has been widely adopted for its robustness and flexibility in handling complex interface dynamics. For researchers interested in the general application and theory of level set methods, or in basic numerical methods and recent research trends, the references~\cite{ref:S99book,ref:OF00,ref:OF01,ref:C05,ref:GFO18,ref:SS20} and the references therein are recommended. In image processing, the level set method has been successfully applied to tasks such as image segmentation, utilizing active contours to detect object boundaries without relying on edge information~\cite{ref:CV01}. In computational fluid dynamics, the method effectively models incompressible two-phase flows, capturing the interface between different fluid phases with high accuracy~\cite{ref:SSO94, ref:SFSO98}

Despite its advantages, implementing accurate boundary conditions in the level set framework remains a significant challenge, particularly on polyhedral meshes in three-dimensional domains.
\pf{Note that a typical choice to represent the initial interface implicitly by a level set function, is to compute the signed distance function by solving the eikonal equation in the whole computational domain \cite{ref:HMFB22,ref:HMF24}.} \pf{To track the dynamic interface, traditional approaches of level set methods} often require reinitialization to maintain the signed distance property, which can introduce errors and perturb the zero level set, reducing accuracy. Reinitialization, as described by a time-dependent formulation~\cite{ref:SSO94}, involves periodically correcting the level set function to ensure it remains a signed distance function. However, this process can inevitably move the zero level set, leading to inaccuracies in the interface location~\cite{ref:SF99,ref:RS00}.

\vv{The banded level set method}~\cite{ref:PMOZK99} provides a computationally efficient way to reinitialize the level set function, ensuring that the distance property is preserved while minimizing perturbations to the zero level set. This method localizes the level set update by only recalculating values near the interface, thus enhancing computational efficiency and maintaining accuracy. \vv{Since the level set function is then kept constant outside a narrow band around the surface~\cite{ref:PMOZK99,ref:AS95}, it may seem that a constant value on $\partial \Omega$ is enough.} However, the localization is challenging in polyhedral cells typically used in complicated computational domains\vv{, where the characteristic lengths within one mesh differ by an order of magnitude. The band moreover carries a boundary of its own, which meets the computational boundary as soon as the surface comes near it, and there the same question returns.}

In this paper, we propose to use the eikonal equation as a natural nonlinear boundary condition for the level set formulation: advective and normal flow equations. The eikonal equation, which characterizes the signed distance function, is solved in the cell-layer adjacent to the boundary, effectively imposing the boundary condition on the level set equation. \vv{What is prescribed at the inflow boundary is then not a value of the level set function, which is in general not known there, but the property that the solution has at the initial time being signed distance function. It is one candidate among the possible choices of boundary conditions, and the numerical examples compare it with those in common use.}

Let us denote $\Omega \subset \mathbb{R}^3$ as Lipschitz, convex, and simply connected computational domain and $T>0$ be the final time of the evolution. The level set formulation uses an implicit function $u : \Omega \times [0,T] \rightarrow \mathbb{R}$ to represent an evolving surface $\Gamma_{t}(u)$ at $t \in [0,T]$ as the zero level set:
\begin{align}\label{eq:zerolevel}
\Gamma_{t}(u) = \{ \mathbf{x} \in \Omega : u(\mathbf{x},t) = 0\}.
\end{align}
Under a given velocity $\mathbf{v}$ on $\Omega$ to evolve the surface,
an equation of evolving surface $\Gamma_{t}(u)$ is described by the level set formulation:
\begin{equation}\label{eq:levelset}
\frac{\partial }{\partial t} u(\mathbf{x},t) + \mathbf{v}(\mathbf{x}, t, \nabla u(\mathbf{x}, t)) \cdot \nabla u(\mathbf{x},t) = 0, \quad (\mathbf{x},t) \in \Omega \times [0,T], 
\end{equation}
where we use \pf{linear combinations of} two choices $\mathbf{v} = \mathbf{v}(\mathbf{x})$ and $\mathbf{v} = \pm \frac{\nabla u}{|\nabla u|}$ to represent advective and normal flows, respectively. An initial condition $u^0(\mathbf{x})$ is given by the signed distance function from the surface $\Gamma_0$:
\begin{equation}\label{eq:initcond}
u(\mathbf{x},0) = u^0(\mathbf{x}), \quad \mathbf{x} \in \Omega.
\end{equation}
According to the convention of the signed distance, if $\Gamma_t(u)$ is a closed surface in $\Omega$, the value of the function $u$ outside the enclosed region is the distance from the surface and the function $u$ takes the negative value of the distance from the surface otherwise. 

To ensure the well-posedness of the problem~\eqref{eq:levelset} governed by the transport equation, it is crucial to prescribe not only an initial condition but also a boundary condition~\cite{ref:A98}. A particular choice of boundary condition should ensure stable computations without artificially affecting the evolution of the surface itself. It is common to choose the zero Neumann boundary condition (ZNBC) when the evolving surface does not touch the boundary.
\begin{equation}\label{eq:zeroNeumann}
\nabla u(\mathbf{x},t) \cdot \mathbf{n}(\mathbf{x}) = 0, \quad (\mathbf{x},t) \in \partial \Omega \times [0,T],
\end{equation}
where $\mathbf{n}$ is the outward normal vector to $\partial \Omega$. It is used together with the local level set method~\cite{ref:PMOZK99,ref:AS95} to make the level set function $u$ constant outside of the local region. If the assumption that the evolving surface is inside the computational domain does not hold, then the zero Neumann boundary condition~\eqref{eq:zeroNeumann} certainly results in a distortion of the evolving surface by forcing the normal of the surface to be orthogonal to the normal of the boundary of the computational domain. 

The Dirichlet boundary condition can also be applied on the inflow boundary~$\partial \Omega^{-}$
\begin{equation}\label{eq:inflow_bdry}
u(\mathbf{x},t) = u_D(\mathbf{x},t), \quad (\mathbf{x},t) \in \partial \Omega^{-} \times [0,T],
\end{equation}
where $\partial \Omega^{-} = \{\mathbf{x} \in \partial \Omega : \mathbf{v}(\mathbf{x},t,\nabla u(\mathbf{x},t)) \cdot \mathbf{n}(\mathbf{x}) < 0 \}$. For example, one can prescribe the values defined by the initial condition with $u_D(\mathbf{x},0) = u^0(\mathbf{x})$, $\mathbf{x} \in \partial \Omega^{-}$, but this choice may negatively impact the accuracy and stability of the evolving surface~\cite{ref:SS20,ref:RL11}. For the purpose of testing the behavior of numerical methods with known exact solutions, the Dirichlet boundary condition is a natural choice. However, when the shape of the evolving surface or the computational domain is complicated in industrial applications, the Dirichlet boundary condition may not be generally known. 

The linearly extended boundary condition (LEBC) is used to solve~\eqref{eq:levelset} on three-dimensional (3D) polyhedral meshes~\cite{ref:HMFB17_1,ref:HMFB17, ref:HMFMB19}. However, if the solution differs significantly from a linear function, there will be a non-negligible error between the linearly extended values and the exact values. The accumulation of such an error can lead to inaccuracy of the numerical solution over long time behavior.

In the numerical solution of level set equations~\eqref{eq:levelset}, a considerable amount of research has focused on the robustness and accuracy of reinitialization procedure \cite{ref:GFO18,ref:SSO94,ref:SFSO98,ref:FHA17,ref:KLK18}, but relatively little attention has been paid to the boundary conditions for the level set function. The purpose of the reinitialization is to recover the property $|\nabla u (\mathbf{x},t)| = 1$ on the whole domain or in a narrow band around the zero level set. The disadvantage of reinitialization is that it can artificially change the position of the evolved surface. Consequently, numerical methods that can avoid or minimize the use of reinitialization steps are an attractive alternative. An approach to achieve it indirectly is to weakly enforce the property $|\nabla u (\mathbf{x},t)| = 1$ by a penalty method~\cite{ref:LXGF05, ref:ZZSZ13} while the evolution equation is solved. A stabilization is further studied in free surface simulation and the stabilized Streamline-Upwind-Petrov-Galerkin (SUPG) method~\cite{ref:TS16} is compared to the SUPG method with Crank-Nicholson method~\cite{ref:RL11,ref:B10} in the finite element method.

In this work, we solve~\eqref{eq:levelset} in~$\Omega$ and the eikonal equation on~$\partial \Omega$ simultaneously, referring as the \textbf{eikonal boundary condition}, to avoid the deviation from $|\nabla u| = 1$ of the evolved level set function due to an inappropriate choice of boundary conditions. Since the initial condition is chosen as a signed distance function of the surface $\Gamma_{0}(u)$, the eikonal boundary condition at $t=0$ is compatible with the initial condition. \pf{Opposite to Dirichlet boundary conditions, which would require for a time evolution of the level set function a knowledge of, e.g., the values of the signed distance function at the boundary, the eikonal boundary conditions can help to preserve the property $|\nabla u| = 1$, at least near the boundary, in a natural way.}

We demonstrate the effectiveness of our method through various numerical experiments, showing that it \nn{is more robust than the linearly extended and the zero Neumann boundary conditions, that it reaches the accuracy of the Dirichlet boundary condition built from the exact solution, which is in general not available, and that it} handles large CFL numbers in simpler cases. We deliberately use simple velocity fields to clearly illustrate the core concept and validate the proposed method. These findings can be applied to more complex velocity fields on hexahedral meshes, incorporating advanced techniques for extending a velocity field~\cite{ref:AS99} defined only on the zero level set. However, addressing the complexities of such extensions on polyhedral meshes lies beyond the scope of this paper. Here, we aim to demonstrate the effectiveness of our approach in simpler, controlled settings, which serves as a foundational step toward more intricate velocity formulations in future work.

The rest of the paper is organized as follows. In the next section, we introduce a basic notation to describe a numerical method on 3D polyhedral meshes and present a numerical algorithm when using the Dirichlet boundary condition. At the end of Section~\ref{sec:EKBC}, we explain the proposed method in detail. Numerical experiments are presented in Section~\ref{sec:NumEx}. Finally, in Section~\ref{sec:conclusion} we conclude the results. Some technical details of the proposed numerical methods are presented in the Appendix.

\section{Eikonal boundary condition}\label{sec:EKBC}

\subsection{Notations}\label{sec:notation}

\begin{figure}
\begin{center}
\begin{tabular}{cc}
	\includegraphics[height=4.5cm]{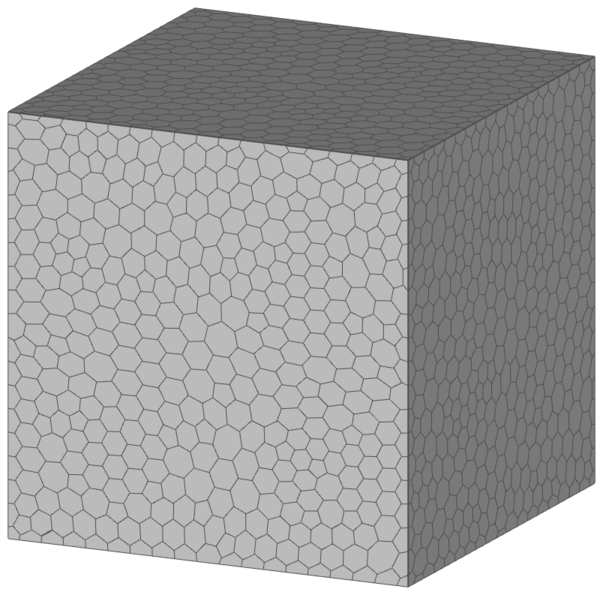} &
	\includegraphics[height=4.5cm]{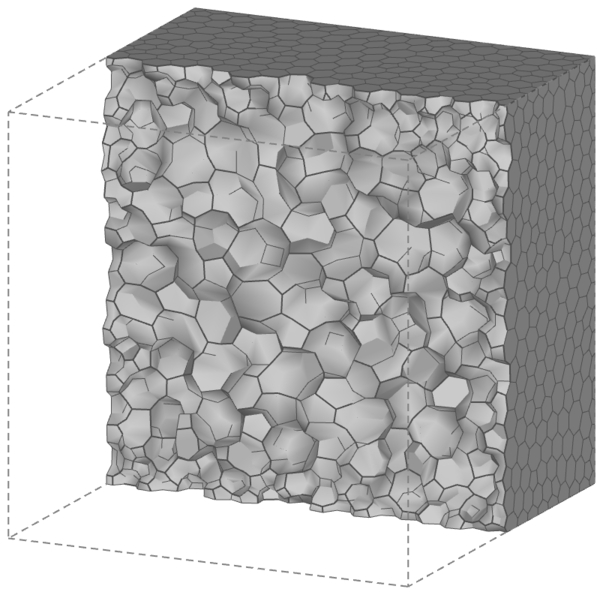} \\
	(a) & (b) \\
	\includegraphics[height=4.5cm]{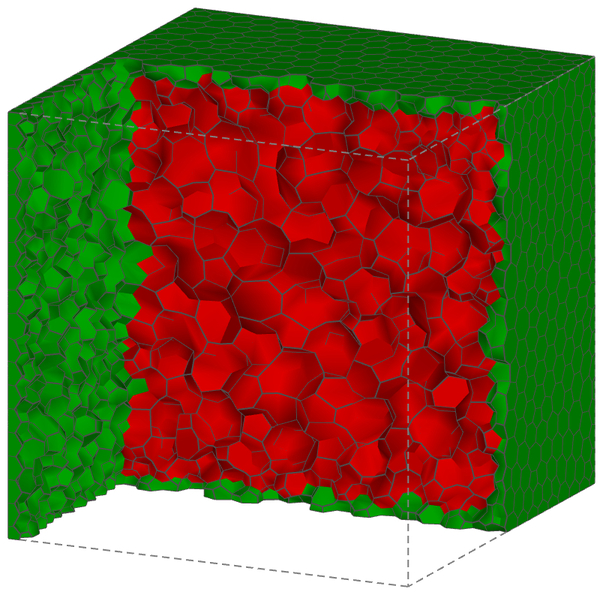} &
	\includegraphics[height=3.2cm]{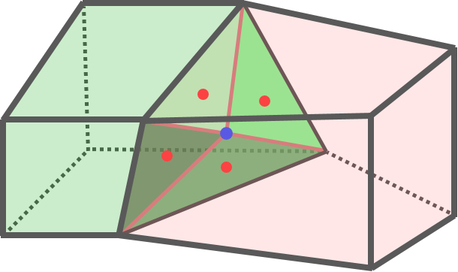} \\
	(c) & (d)
\end{tabular}
\end{center}
\caption{(a) is a polyhedron mesh whose boundary is a cube. (b) shows the inner cross-section of (a). (c) presents a part of the boundary (green) and internal (red) cells. (d) is an illustration of boundary (green) and internal (red) cells with tessellated triangular faces.} \label{fig:cell_dig}
\end{figure}

A computational domain~$\Omega$ is discretized by a union of non-overlapped polyhedral cells with a non-zero volume $|\Omega_p| > 0$:
\begin{equation}\label{eq:dis_dom}
\bar{\Omega} = \bigcup_{p \in \mathcal{I}} \bar{\Omega}_p,
\end{equation}
where $\Omega_p$ is open and $\mathcal{I}$ is a set of the indices of cells. If there is a common surface between two adjacent cells $\Omega_p$ and $\Omega_q$, $p, q \in \mathcal{I}$, we call the surface an \textit{internal face}. If a boundary surface of the cell $\Omega_p$ has a non-empty common area with $\partial \Omega$, we call the surface a \textit{boundary face}. Since a boundary surface of a polyhedron cell in the three-dimensional (3D) domain is mostly not a plane, the surface is always tessellated by triangles; see an example of the surface tessellated by four red triangles in Figure~\ref{fig:cell_dig}-(d); see more details described by \ref{A:step2} in the Appendix or the paper~\cite{ref:HMFMB19}. We denote $\mathcal{F}_p$  and $\mathcal{B}_p$ as indices of all \textit{triangles tessellated from the internal and boundary face} of the cell $\Omega_p$, respectively. We call a cell $\Omega_p$ as an \textit{internal cell} if $p \in \mathcal{I}_{\text{int}} = \{ p \in \mathcal{I} :\mathcal{B}_p = \emptyset \}$. A cell $\Omega_p$, $p \in \mathcal{I}_{\text{bdr}} = \mathcal{I} \setminus \mathcal{I}_{\text{int}}$ is called a \textit{boundary cell}; see red (internal) or green (boundary) cells in Figure~\ref{fig:cell_dig}-(c). For a cell $\Omega_p$, $p \in \mathcal{I}$, we define a set $\mathcal{N}_p$ as the indices of \textit{$1$-ring face neighbor cells}~$\Omega_q$, $q \in \mathcal{I}$, such that the intersection $\partial \Omega_p \cap \partial \Omega_q$ is a face of the non-zero area between two adjacent cells. 

We define a set of the indices of tessellated internal and boundary faces, that is, triangles, as $\mathcal{F}$ and $\mathcal{B}$. In Figure~\ref{fig:cell_dig}-(d), if the left (green) and right (red) cells are $\Omega_p$ and $\Omega_q$, $p,q\in\mathcal{I}$, respectively, we denote $\mathbf{x}_p$ and $\mathbf{x}_q$ as the center of corresponding cells; see the formula~\eqref{eq:app_cell_ct} in the Appendix or in the paper~\cite{ref:HMF24}. Between left and right cells, there is $e_f$, $f \in \mathcal{F}$, an internal triangle. For the $e_f$, $f\in \mathcal{F}_p$, we denote $\mathbf{x}_f$ (red point) as the center of the triangle on $\partial \Omega_p$ and the vector $\mathbf{n}_{pf}$ as the outward unit normal to the triangle and $|e_f|$ as the area fo the triangle. Then, $e_f \subset \partial \Omega_q$ for $q\in\mathcal{N}_p$, we have $\mathbf{n}_{qf} = -\mathbf{n}_{pf}$. For a boundary triangle $e_b$, $b\in\mathcal{B}_p$, the unit vector $\mathbf{n}_{b} = \mathbf{n}_{pb}$ is normal to the triangle which is the outward normal to $\partial \Omega$. We use a directional vector specified by two position vectors $\mathbf{x}_a$ and $\mathbf{x}_b$ as a notation $\mathbf{d}_{ab} = \mathbf{x}_b - \mathbf{x}_a$. In the remainder of the paper, the subscripts $f$, $b$, and $a$ are used to denote an internal triangle $e_f$, a boundary triangle $e_b$, and any of these triangles, respectively, and the subscripts $p$ or $q$ are used to denote a cell unless otherwise noted.

\subsection{Proposed algorithm}\label{sec:prop}
The former part of the subsection is devoted to explaining a concise overview of the cell-centered finite volume method for solving the level set equation~\eqref{eq:levelset} on polyhedral meshes. Although the required information is already presented in the paper~\cite{ref:HMFMB19}, the overview is necessary to clarify the proposed method for simultaneously solving~\eqref{eq:levelset} and the eikonal equation in the latter part of the subsection.

We rewrite the level set equation~\eqref{eq:levelset}:
\begin{equation}
\frac{\partial u}{\partial t}(\mathbf{x}, t) + \nabla \cdot \left(u \mathbf{v} \right) (\mathbf{x}, t) - u(\mathbf{x}, t) \nabla \cdot \mathbf{v}(\mathbf{x}, t) = 0,
\end{equation}
and it can be evaluated at $(\mathbf{x}_p, t) \in \Omega_p \times [0,T]$:
\begin{equation}
\label{eq-lseq-at-point}
\pf{
\frac{\partial}{\partial t}u(\mathbf{x}_p, t) + \nabla \cdot \left(u \mathbf{v} \right) (\mathbf{x}_p, t) - u(\mathbf{x}_p, t)} \nabla \cdot \mathbf{v}(\mathbf{x}_p, t) = 0.
\end{equation}
\pf{For $\mathbf{G} = u \mathbf{v}$ and $\mathbf{G} =\mathbf{v}$, we approximate $\nabla \cdot \mathbf{G}(\mathbf{x}_p,t)$ by averaging} the divergence over the volume of the cell $\Omega_p$ and applying then Gauss's theorem to the integral, that is,
\begin{equation}\label{eq:GaussTh}
\nabla \cdot \mathbf{G}(\mathbf{x}_p,t) \pfb{\approx} \frac{1}{|\Omega_p|} \int_{\Omega_p} \nabla \cdot \mathbf{G}(\mathbf{x},t) dV = \frac{1}{|\Omega_p|} \int_{\partial \Omega_p} \mathbf{G} \cdot \mathbf{n} dS,
\end{equation}
where $\mathbf{n}$ is the outward normal to $\partial \Omega_p$. \pf{Applying these approximations in \eqref{eq-lseq-at-point} and using $u_p(t) \approx u(\mathbf{x}_p,t)$,}  we obtain the following formulation,
\begin{equation}\label{eq:ccFVM_s}
\pf{\frac{d}{dt}}u_p(t) + \frac{1}{|\Omega_p|}\int_{\partial \Omega_p} u \mathbf{v} \cdot \mathbf{n} dS - \frac{u_p(t)}{|\Omega_p|}\int_{\partial \Omega_p} \mathbf{v} \cdot \mathbf{n} dS = 0.
\end{equation}
For each face $e_a$, $a \in \mathcal{F}_p \cup \mathcal{B}_p$, we define the flux,
\begin{equation}\label{eq:flux_mu}
\mu_{pa}(t) = \int_{e_{a}} \mathbf{v} \cdot \mathbf{n}_{pa} dA \approx \mathbf{v} (\mathbf{x}_a, t, \nabla u(\mathbf{x}_a, t)) \cdot \mathbf{n}_{pa} |e_a|,
\end{equation}
where $\mathbf{x}_a$ is the center of the face and $\nabla u(\mathbf{x}_a,t)$ is a representative gradient on the triangle~$e_a$, that we explain how to compute in the Appendix. Finally, we have the spatial discretization:
\begin{equation}\label{eq:ccFVM}
|\Omega_p|\pf{\frac{d}{dt}} u_p(t) +  \sum_{a \in \mathcal{F}_{p} \cup \mathcal{B}_{p}} \left( u_{pa}(t) - u_p(t) \right) \mu_{pa}(t) = 0,
\end{equation}
where the value $u_{pa}(t)$ is determined by the sign of the flux $\mu_{pa}(t)$ and we explain the evaluation of $u_{pa}(t)$ below.

Considering a fixed time step $\triangle t$ and $t^n = n \triangle t$, we would like to find a solution $u_p^n = u(\mathbf{x}_p,t^n)$ for all $p \in \mathcal{I}$ and $n\in\mathbb{N}$ from known values $u_p^{n-1}$. Let us denote the index sets of triangles depending on the sign of the flux at $t= t^{n-1}$:
\begin{equation}\label{eq:flux_sign_mu}
\begin{alignedat}{3}
\mathcal{F}^{\mu^{-}}_{p} &= \{ f \in \mathcal{F}_{p} \: | \: \mu_{pf}^{n-1} < 0 \}, \quad &&\mathcal{F}^{\mu^{+}}_{p} = \mathcal{F}_{p} \setminus \mathcal{F}^{\mu^{-}}_{p},\\
\mathcal{B}^{\mu^{-}}_{p} &= \{ b \in \mathcal{B}_{p} \: | \: \mu_{pb}^{n-1}  < 0 \}, \quad &&\mathcal{B}^{\mu^{+}}_{p} = \mathcal{B}_{p} \setminus \mathcal{B}^{\mu^{-}}_{p},	
\end{alignedat}
\end{equation}
where $\mu_{pa}^{n-1} = \mu_{pa}(t^{n-1})$, $a \in \mathcal{F} \cup \mathcal{B}$. The upwind principle in the evaluation of~$u_{pa}$ in~\eqref{eq:ccFVM} is used depending on the sign of $\mu_{pa}^{n-1}$ and we apply the inflow-implicit outflow-explicit~(IIOE) method~\cite{ref:MO11,ref:MOU14}. The core concept of the IIOE method is to handle the outflow from a cell explicitly while treating the inflow implicitly, ensuring stability by constructing an M-matrix from the inflow fluxes, which provides favorable solvability properties. Furthermore, we denote a representative gradient, $\mathcal{D} u_p^n$ on the cell $\Omega_p$, which is calculated distinctly for internal and boundary cells. The method varies based on the applied boundary conditions, with details on the computation process available in the Appendix. Now, we have the discretization scheme of~\eqref{eq:levelset} on internal cells~$\Omega_p$,~$p \in \mathcal{I}_{\text{int}}$:
\begin{equation}\label{eq:algo_int}
\begin{split}
&\frac{|\Omega_p|}{\pfb{\triangle} t} \left( u_p^{n} - u_p^{n-1} \right) + \sum_{f \in \mathcal{F}^{\mu^{-}}_{p}} \left( u_q^{n} + \mathcal{D} u_{q}^{n} \cdot \mathbf{d}_{qf}  - u_p^{n} \right) \mu_{pf}^{n-1} \\
&  + \sum_{a \in \mathcal{F}^{\mu^{+}}_{p} } \left( \mathcal{D} u_{p}^{n-1} \cdot \mathbf{d}_{pa} \right) \mu_{pa}^{n-1} = 0.
\end{split}
\end{equation}
\pfb{For an internal triangle $e_f$, $f\in\mathcal{F}_p$, the index $q=q(f)\in\mathcal{N}_p$ denotes the unique neighboring cell whose shared face contains $e_f$, that is, $e_f\subset\partial\Omega_p\cap\partial\Omega_q$; this index~$q$ is also used in the remaining discretization schemes.} If we consider the Dirichlet boundary condition~\cite{ref:HMFMB19}, the discretization scheme needs to include the values on boundary cells~$\Omega_p$,~$p \in \mathcal{I}_{\text{bdr}}$:
\begin{equation}\label{eq:algo_DBC}
\begin{split}
&\frac{|\Omega_p|}{\pfb{\triangle} t} \left( u_p^{n} - u_p^{n-1} \right) + \sum_{f \in \mathcal{F}^{\mu^{-}}_{p}} \left( u_q^{n} + \mathcal{D} u_{q}^{n} \cdot \mathbf{d}_{qf}  - u_p^{n} \right) \mu_{pf}^{n-1} \\
& + \sum_{b \in \mathcal{B}^{\mu^{-}}_{p}} \left( u_b^{n} - u_p^{n} \right) \mu_{pb}^{n-1}  + \sum_{a \in \mathcal{F}^{\mu^{+}}_{p} \cup \mathcal{B}^{\mu^{+}}_{p}} \left( \mathcal{D} u_{p}^{n-1} \cdot \mathbf{d}_{pa} \right) \mu_{pa}^{n-1} = 0,
\end{split}
\end{equation}
where $u_b^n = u_D(\mathbf{x}_b, t^n)$ is given.

In contrast to solving~\eqref{eq:algo_DBC} on the whole computational cells, the central concept of the proposed algorithm is to make the coupling to solve numerically the governing equation~\eqref{eq:levelset} on internal cell~$\Omega_p$,~$p \in \mathcal{I}_{\text{int}}$ and the eikonal equation $|\nabla u| = 1$ on boundary cells~$\Omega_b$,~$b \in \mathcal{I}_{\text{bdr}}$, simultaneously. In the finite volume method, a similar concept is used to compute the oblique derivative boundary value problems~\cite{ref:MCMM15}. In the finite element method, the nonlinear satellite-fixed geodetic boundary value problem~\cite{ref:MMCM23} used the eikonal equation on the boundary.

For the internal cells $\Omega_p$, $p \in \mathcal{I}_{\text{int}}$, we already have the discretization scheme~\eqref{eq:algo_int} to solve~\eqref{eq:levelset}. Therefore, to complete the proposed method, it is sufficient to explain the discretization of the eikonal equation on the boundary cells. In order to deal with the nonlinearity of the eikonal equation, a linearization by a semi-implicit approach is used at $t=t^n$ by the known solution $u(\mathbf{x},t^{n-1})$ from the previous time step:
\begin{equation}\label{eq:lin_ek}
\frac{\nabla u(\mathbf{x},t^{n-1})}{|\nabla u(\mathbf{x},t^{n-1})|} \cdot \nabla u(\mathbf{x},t^n) = 1, \quad \mathbf{x} \in \Omega_p, \:\: p \in \mathcal{I}_{\text{bdr}}\pfb{.}
\end{equation}
Applying the regularized absolute value $|\mathbf{x}|_\delta = \sqrt{|\mathbf{x}|^2 + \delta^2}$ with a small constant $\delta = 10^{-12}$,
we evaluate $u$ and its gradient  at $\mathbf{x}_p \in \Omega_p$, $p \in \mathcal{I}_{\text{bdr}}$ in the above equation:
\begin{equation}\label{eq:lin_ek_vec_cal}
\mathbf{w}(\nabla u_p(t^{n-1})) \cdot \nabla u_p(t^n) = \nabla \cdot (u \mathbf{w})(\nabla u_p(t^{n-1})) - u_p(t^{n-1}) \nabla \cdot \mathbf{w}(\nabla u_p(t^{n-1})) = 1,
\end{equation}
where $\nabla u_p(t^{n-1}) = \nabla u(\mathbf{x}_p, t^{n-1})$ and 
\begin{align}
\mathbf{w}(\nabla u(\mathbf{x},  t^{n-1})) = \frac{\nabla u(\mathbf{x}, t^{n-1})}{|\nabla u(\mathbf{x}, t^{n-1})|_\delta}, \quad \mathbf{x} \in \Omega_p, \:\: p \in \mathcal{I}_{\text{bdr}}.
\end{align}
We use the same method~\eqref{eq:GaussTh} for the equation~\eqref{eq:lin_ek_vec_cal} to derive
\begin{align}
\frac{1}{|\Omega_p|} \int_{\partial \Omega_p} u\mathbf{w} \cdot \mathbf{n}dS - \frac{u_p(t^{n-1})}{|\Omega_p|} \int_{\partial \Omega_p} \mathbf{w} \cdot \mathbf{n} dS = 1.
\end{align}
Defining the regularized flux on the face $e_a$, $a \in \mathcal{F} \cup \mathcal{B}$,
\begin{align}\label{eq:flux_nu}
\nu_{pa}(t^{n-1}) = \int_{e_a} \mathbf{w}(\nabla u(\mathbf{x}, t^{n-1})) \cdot \mathbf{n}_{pa} dA \approx \mathbf{w}(\nabla u(\mathbf{x}_a, t^{n-1})) \cdot \mathbf{n}_{pa} |e_a|,
\end{align}
where $\nabla u(\mathbf{x}_a,t^{n-1})$ is a representative gradient on the triangle~$e_a$, which we explain how to compute in the Appendix, the discretization scheme is obtained:
\begin{align}\label{eq:ccFVM_bdry}
\sum_{a \in \mathcal{F}_{p} \cup \mathcal{B}_{p}} \left( u_{pa}(t^{n}) - u_p(t^{n}) \right) \nu_{pa}(t^{n-1}) = |\Omega_p|,
\end{align}
where the value $u_{pa}(t^{n})$ is determined by the sign of the flux $\nu_{pa}(t^{n-1})$ and we explain the evaluation of~$u_{pa}(t^{n})$ below.

Now, using the upwind principle in the evaluation of~$u_{pa}(t^n)$ in~\eqref{eq:ccFVM_bdry} and the index sets of triangles depending on the sign of the flux~$\nu_{pa}^{n-1} = \nu_{pa}(t^{n-1})$,
\begin{equation}\label{eq:flux_sign_nu}
\begin{alignedat}{3}	
\mathcal{F}_{p}^{\nu^{-}} &= \{ f \in \mathcal{F}_p : \nu_{pf}^{n-1} < 0 \}, \quad \mathcal{F}_{p}^{\nu^{+}} &&= \mathcal{F}_p \setminus \mathcal{F}_{p}^{\nu^{-}},  \\
\mathcal{B}_{p}^{\nu^{-}} &= \{ b \in \mathcal{B}_p : \nu_{pb}^{n-1} < 0 \}, \quad \mathcal{B}_{p}^{\nu^{+}} &&= \pfb{\mathcal{B}}_p \setminus \mathcal{B}_{p}^{\nu^{-}},
\end{alignedat}
\end{equation}
we have the discretization scheme of the eikonal equation on the boundary cells $\Omega_p$, $p \in \mathcal{I}_{\text{bdr}}$:
\begin{equation}\label{eq:algo_bdr}
\begin{split}
&\sum_{f \in \mathcal{F}^{\nu^{-}}_{p}} \left( u_q^{n} + \mathcal{D} u_q^{n} \cdot \mathbf{d}_{qf} - u_p^{n} \right) \nu_{pf}^{n-1} + \sum_{f \in \mathcal{F}^{\nu^{+}}_{p}} \left( \mathcal{D} u_p^{n} \cdot \mathbf{d}_{pf} \right) \nu_{pf}^{n-1} \\
&+ \sum_{b \in \mathcal{B}^{\nu^{+}}_{p} } \left( \mathcal{D} u_p^{n} \cdot \mathbf{d}_{pb} \right) \nu_{pb}^{n-1} =  |\Omega_p|.
\end{split}	
\end{equation}
It is crucial to note that, unlike in equation~\eqref{eq:algo_DBC}, there are no terms on inflow boundary triangles, $e_b$, $b \in \mathcal{B}_p^{\nu^{-}}$, because including them would violate the Soner boundary condition.
\nn{It enforces the constraint $\mathbf{n} (\mathbf{x}) \cdot \nabla u(\mathbf{x}) \geq 0$ on the boundary, which excludes the inadmissible control directions and so selects the viscosity solution of the eikonal equation in the boundary cells;} see more details in the papers~\cite{ref:HMF24,ref:S86,ref:CDL90,ref:HMFB22}. If we solve~\eqref{eq:levelset} also on the boundary cells, we clearly see that the discretization of the inflow boundary is necessary and certain boundary values must be specified. However, in this paper, we propose to solve the eikonal equation on the boundary cells, which eliminates discretization on inflow boundary triangles.

In order to solve the proposed method~\eqref{eq:algo_int} and~\eqref{eq:algo_bdr} practically on polyhedral meshes by parallel computing using domain decomposition with 1-ring face neighborhood, we apply a deferred correction method~\cite{ref:BHS84}. Let us denote~$u^{0}$ the initial condition~\eqref{eq:initcond} and~$u^{n,0} = u^{n-1}$,~$n \in \mathbb{N}$. Then, we would like to find the unknown values~$u_p^{n,k}$ from the known values~$u_p^{n,k-1}$ for~$k \in \mathbb{N}$ and all~$p \in \mathcal{I}$ from the couple equations: for an internal cell~$\Omega_p$,~$p \in \mathcal{I}_{\text{int}}$, 
\begin{equation}\label{eq:deffered_algo_int}
\begin{split}
&\frac{|\Omega_p|}{\pfb{\triangle} t} \left( u_p^{n,k} - u_p^{n-1} \right) + \sum_{f \in \mathcal{F}^{\mu^{-}}_{p}} \left( u_q^{n,k} + \mathcal{D} u_q^{n,k-1} \cdot \mathbf{d}_{qf}  - u_p^{n,k} \right) \mu_{pf}^{n-1} \\ &+ \sum_{f \in \mathcal{F}^{\mu^{+}}_{p} } \left( \mathcal{D} u_p^{n-1} \cdot \mathbf{d}_{pf} \right) \mu_{pf}^{n-1} = 0,
\end{split}
\end{equation}
and for a boundary cell~$\Omega_p$, ~$p \in \mathcal{I}_{\text{bdr}}$, 
\begin{equation}\label{eq:deffered_algo_bdr}
\begin{split}
&\sum_{f \in \mathcal{F}^{\nu^{-}}_{p}} \left( u_q^{n,k} + \mathcal{D} u_q^{n,k-1} \cdot \mathbf{d}_{qf} - u_p^{n,k} \right) \nu_{pf}^{n-1} + \sum_{f \in \mathcal{F}^{\nu^{+}}_{p}} \left( \mathcal{D} u_p^{n,k-1} \cdot \mathbf{d}_{pf} \right) \nu_{pf}^{n-1} \\
&+ \sum_{b \in \mathcal{B}^{\nu^{+}}_{p} } \left( \mathcal{D} u_p^{n,k-1} \cdot \mathbf{d}_{pb} \right) \nu_{pb}^{n-1} =  |\Omega_p|.
\end{split}
\end{equation}
Rewriting~\eqref{eq:deffered_algo_int} and~\eqref{eq:deffered_algo_bdr} formally as a matrix equation, 
\begin{equation}\label{eq:mateq}
\mathbf{A}^{n-1} u^{n,k} = \mathbf{f}(u^{n,k-1}),
\end{equation}
it is solved by an algebraic multigrid method for each $k^{\text{th}}$ iteration. Then the $k$-iteration can be stopped at the smallest integer to make the residual error to be smaller than a chosen bound $\eta=10^{-12}$:
\begin{equation}\label{eq:residual}
K = \min\left\{ k \in \mathbb{N} : \frac{1}{|\mathcal{I}|} \sum_{p \in \mathcal{I}} \left| \left(\mathbf{A}^{n-1}u^{n,k} - \mathbf{f}(u^{n,k})\right)_p \right| < \eta\right\},
\end{equation}
where the parenthesis with a subscript $\left( \mathbf{r} \right)_p$ denotes the $p^{\text{th}}$ component of the vector $\mathbf{r}$. Finally, we select the solution of~$u$ at~$t=t^n$ as $u^{n,K}$ and proceed to compute the solution at the next time step.

When the Dirichlet boundary condition is used, the same deferred correction method is applied to~\eqref{eq:algo_int} and~\eqref{eq:algo_DBC}:
\begin{equation}\label{eq:deffered_algo_DBC}
\begin{split}
&\frac{|\Omega_p|}{\pfb{\triangle} t} \left( u_p^{n,k} - u_p^{n-1} \right) + \sum_{f \in \mathcal{F}^{\mu^{-}}_{p}} \left( u_q^{n,k} + \mathcal{D} u_{q}^{n,k-1} \cdot \mathbf{d}_{qf}  - u_p^{n,k} \right) \mu_{pf}^{n-1} \\
&+ \sum_{b \in \mathcal{B}^{\mu^{-}}_{p}} \left( u_b^{n} - u_p^{n,k} \right) \mu_{pb}^{n-1} + \sum_{a \in \mathcal{B}^{\mu^{+}}_{p} \cup \mathcal{F}^{\mu^{+}}_{p}} \left( \mathcal{D} u_{p}^{n-1} \cdot \mathbf{d}_{pa} \right) \mu_{pa}^{n-1} = 0,
\end{split}
\end{equation}
where $u_b^{n} = u_D(\mathbf{x}_b,t^n)$ is given by the Dirichlet boundary condition.

\begin{algorithm}[t] 
\caption{A solution procedure per time step} \label{alg:EKBC}
\begin{algorithmic}
\Procedure{Eikonal boundary condition for level set method}{} 
\State Initialization of signed distance function $u^0(\mathbf{x}) = u(\mathbf{x},0)$ by~\eqref{eq:initcond}.
\State Set $n=1$.
\While{$n \pfb{\triangle} t \leq T$} 
\State Compute $\mu_{pf}^{n-1}$ and $\mathcal{D}^{n-1}_p u$ in~\eqref{eq:deffered_algo_int} and  $\nu_{pf}^{n-1}$ in~\eqref{eq:deffered_algo_bdr}.
\While{$k \leq K$}
\State Compute $\mathcal{D}^{n,k-1}_p u$
\State Solve the matrix equation~\eqref{eq:mateq} to find $u^{n}$.
\EndWhile  
\State $n \leftarrow n + 1$.
\EndWhile  
\EndProcedure
\end{algorithmic}
\end{algorithm}

\begin{remark}
In the standard cell-centered FVM, the unknown values are located at the centers of the cells. Therefore, in a standard scheme, it is required to approximate the discretization on the boundary using the values at the centers of the boundary cells. By solving the eikonal equation in the cell-layer adjacent to the boundary, that is, boundary cells in Figure~\ref{fig:cell_dig}-(c), we effectively impose a boundary condition on the level set equation. This approach ensures that the boundary conditions are accurately reflected in the computational domain, maintaining the integrity of the method.
\end{remark}

\begin{remark}
It is a common approach to add ``ghost cells'' outside the domain to handle boundary conditions, such as Neumann boundary conditions, and this practice is widely accepted. The only difference in our approach is that we create this layer inside the domain because adding it outside would be challenging when the boundary of the computational domain is a complex shape. This layer is refined as we refine the mesh, ensuring convergence to the boundary.
\end{remark}

\section{Numerical results}\label{sec:NumEx}

\begin{table}
\centering	
\begin{tabular}{
	c
	S[table-format=1]
	S[table-format=7]
	S[table-format=1.2,table-figures-exponent=2,table-sign-mantissa,table-sign-exponent]
	S[table-format=1.2,table-figures-exponent=2,table-sign-mantissa,table-sign-exponent]
	S[table-format=1.2,table-figures-exponent=2,table-sign-mantissa,table-sign-exponent]}
\toprule
mesh & \multicolumn{1}{c}{M} & \multicolumn{1}{c}{$|\mathcal{I}_{\text{M}}|$ } & \multicolumn{1}{c}{$h_{\text{M}}^{\text{min}}$} & \multicolumn{1}{c}{$h_{\text{M}}^{\text{ave}}$} & \multicolumn{1}{c}{$h_{\text{M}}^{\text{max}}$} \\
\midrule
\multirow{4}{*}{$\mathcal{M}_{\text{M}}$} & 1 & 4079 & 6.56E-02 & 1.90E-01 & 4.21E-01 \\
& 2 & 32004 & 2.29E-02 & 9.52E-02 & 2.00E-01 \\
& 3 & 252433 & 1.36E-02 & 4.76E-02 & 9.46E-02 \\
& 4 & 2024478 & 5.51E-03 & 2.48E-02 & 4.48E-02 \\
\bottomrule
\end{tabular}
\caption{The numbers of polyhedral cells $\mathcal{I}_{\text{M}}$ and the characteristic lengths $h_{\text{M}}^{\text{min}}$, $h_{\text{M}}^{\text{ave}}$, and $h_{\text{M}}^{\text{max}}$~\eqref{eq:char_len} of the meshes are presented; see the shape of the computational domains at the level $\text{M}=1$ in Figures~\ref{fig:cell_dig}-(a) and (b).}\label{tab:meshes}
\end{table}

Various examples are presented to show the numerical properties of using the eikonal boundary condition. The polyhedral mesh discretizing the cube:
\begin{equation}
\Omega = [-1.25,1.25]^3 \subset \mathbb{R}^3
\end{equation}
generated by \texttt{AVL FIRE}\textsuperscript{\texttt{TM}} is illustrated in Figures~\ref{fig:cell_dig}-(a) and (b) and the number of polyhedral cells $|\mathcal{I}_{\text{M}}|$ and the characteristic length $h_{\text{M}}^{\text{ave}}$ of the meshes are presented for four levels of meshes, $\text{M} \in \{1,2,3,4\}$, in Table~\ref{tab:meshes}, where
\begin{equation}\label{eq:char_len}
h^{\text{min}}_{\text{M}} = \min_{p \in \mathcal{I}_{\text{M}}}  |\Omega_p|_{B}^{\frac{1}{3}}, \quad		
h_{\text{M}}^{\text{ave}} = \frac{1}{|\mathcal{I}_{\text{M}}|} \sum_{p \in \mathcal{I}_{\text{M}}} |\Omega_p|_{B}^{\frac{1}{3}},\quad	
h^{\text{max}}_{\text{M}} = \max_{p \in \mathcal{I}_{\text{M}}}  |\Omega_p|_{B}^{\frac{1}{3}}.
\end{equation}
where $|\Omega_p|_{B}$ is the volume of the box whose diagonal is a vector $\mathbf{x}_M - \mathbf{x}_m$, where $\mathbf{x}_m$ and $\mathbf{x}_M$ are componentwise minimum and maximum of all vertices of the cell $\Omega_p$, respectively. The $\text{M}$ means the level of mesh refinement. When $\text{M}$ increases, finer cells are generated and we roughly have $h_{{\text{M}} + 1}^{\text{ave}} \approx \frac{1}{2} h_{\text{M}}^{\text{ave}}$ to check the experimental order of convergence ($EOC$).

For the quantitative comparison of the numerical results, an experimental order of convergence ($EOC$) is checked. We use the final time $T>0$ and the time step 
\begin{equation}\label{eq:time_step}
\triangle t_{\text{M}} = \frac{0.1}{2^{\text{M}-1}}, \quad \text{M} \in \{1,2,3,4\},
\end{equation}
is used on the corresponding level of the mesh $\mathcal{M}_{\text{M}}$ in Table~\ref{tab:meshes}. Since exact solutions for all examples are explicitly known, the following errors are reported for the test cases in Table~\ref{tab:test_case}. Let $u_e$ be the exact solution.
\begin{itemize}
\item The error $E^{1,\mathcal{Z}}$ measures the average deviation of the numerical solution from the exact solution at the zero level set over the entire simulation period. It is crucial for assessing how accurately the method captures the interface dynamics.
\begin{equation}\label{eq:err_L1_zero}
\begin{split}
	E^{1,\mathcal{Z}}_{\text{M}} &= \frac{1}{T \sum_{p \in \mathcal{Z}_\text{M}} |\Omega_p| } \sum_{n=1}^{N} \sum_{p \in \mathcal{Z}_\text{M}} |u(\mathbf{x}_p,t^{n}) - u_e(\mathbf{x}_p,t^{n})| |\Omega_p| \pfb{\triangle t_{\text{M}}},
\end{split}
\end{equation}
where $\mathcal{Z}_\text{M}$ is the set of the index of cells to contain the zero level set of the exact solution.
\item The error $E^{\infty,\mathcal{Z}}$ represents the maximum pointwise deviation at the zero level set throughout the simulation. It highlights the worst-case scenario, providing insight into the method's robustness against peak errors.
\begin{equation}\label{eq:err_Linf_zero}
\begin{split}
	E^{\infty,\mathcal{Z}}_{\text{M}} &= \max_{n \in \{1,2,\ldots,N\}} \max_{p \in \mathcal{Z}_\text{M}} |u(\mathbf{x}_p,t^n) - u_e(\mathbf{x}_p,t^n)|.
\end{split}
\end{equation}
\item The error $E^{\text{v}}$ measures the difference in the volume enclosed by the evolving surface between the numerical and exact solutions. It is essential for evaluating how well the method preserves the volume of the shape, which is critical in applications such as fluid dynamics and material science.
\begin{equation}\label{eq:err_vol}
\begin{split}
	E^{\text{v}}_{\text{M}} &= \frac{1}{N} \sum_{n=1}^{N} |V(\Gamma_{t^{n}}(u)) - V(\Gamma_{t^{n}}(u_e))|,
\end{split}
\end{equation}	
where $V(\Gamma)$ is the volume of the \pfb{region} enclosed by $\Gamma$.
\item The error $E^1$ is the $\pfb{L^1(\Omega\times(0,T))}$ norm of the error over the entire computational domain. It provides a global measure of the numerical solution's accuracy, indicating how well the method performs across the entire domain.
\begin{equation}\label{eq:err_L1}
\begin{split}
	E^1_{\text{M}} &= \frac{1}{T \sum_{p \in \mathcal{I}_\text{M}} |\Omega_p|} \sum_{n=1}^{N} \sum_{p \in \mathcal{I}_\text{M}} |u(\mathbf{x}_p,t^n) - u_e(\mathbf{x}_p,t^n)| |\Omega_p\pfb{|} \pfb{\triangle t_{\text{M}}}.
\end{split}
\end{equation}
\item The error $E^{1,g}$ is the $\pfb{L^1(\Omega\times(0,T))}$ norm of the error to check the deviation of \pfb{the representative gradient} $|\nabla u|$ \pfb{defined in}~\eqref{eq:grad} from $1$ over the entire computational domain.
\begin{equation}\label{eq:err_L1_g}
\begin{split}
	E^{1,g}_{\text{M}} &= \frac{1}{T \sum_{p \in \mathcal{I}_\text{M}} |\Omega_p|} \sum_{n=1}^{N} \sum_{p \in \mathcal{I}_\text{M}} \big| |\nabla u(\mathbf{x}_p,t^n)| - 1 \big| |\Omega_p| \pfb{\triangle t_{\text{M}}},
\end{split}
\end{equation}
\end{itemize}
Note that the volume in~\eqref{eq:err_vol} is computed by the extracting the interface $\Gamma$ using the marching tetrahedron algorithm. This algorithm is applied to a tetrahedral cell whose base is a face $e_f$ and whose apex is the vertex $\mathbf{x}_p$, for $\forall f \in \mathcal{F}_p$, $\forall p \in \mathcal{I}$. The values at vertices of the tetrahedron are computed by~\eqref{eq:ver_val} in the Appendix.

Then, for an error $E_{\text{M}} \in \{E^{1,\mathcal{Z}}_{\text{M}},\: E^{\infty,\mathcal{Z}}_{\text{M}}, \:E^{\text{v}}_{\text{M}},\: E^1_{\text{M}},\:E^{1,g}_{\text{M}}\}$ on the mesh $\mathcal{M}_{\text{M}}$, the corresponding $EOC_{E_\text{M}}$ is calculated by
\begin{equation}\label{eq:EOC}
EOC_{E_\text{M}} = \pfb{\log\!\left(\frac{E_\text{M+1}}{E_\text{M}}\right) \Big/ \log\!\left(\frac{h_\text{M+1}^{\text{ave}}}{h_\text{M}^{\text{ave}}}\right)}, \quad \text{M} \in \{1,\:2,\:3\},
\end{equation}
where $h_\text{M}^{\text{ave}}$ is the characteristic length~\eqref{eq:char_len} in Table~\ref{tab:meshes}.  

\subsection{Comparison to Dirichlet boundary condition}\label{sec:comp_EDBC}

In this subsection, we numerically compare the quantitative results of using the eikonal and Dirichlet boundary conditions. For all test cases, we need to specify the end time $T$, the velocity field $\mathbf{v}$ in~\eqref{eq:levelset}, and the initial surface $\Gamma_{0}(u)$ that is the zero level set of the signed distance function $u^0$ in~\eqref{eq:initcond}. Using equations of sphere and cube,
\begin{equation}\label{eq:sphere_cube}
\begin{split}
\Lambda_1(\mathbf{c},r) &= \{ \mathbf{x} \in \mathbb{R}^3 : |\mathbf{x} - \mathbf{c}| = r\}, \\
\Lambda_2(\mathbf{c},r) &= \{ \mathbf{x} \in \mathbb{R}^3 : \max_{i} |x_i - c_i| = r\},
\end{split}
\end{equation}
one can describe the \pf{exact} solutions of the test cases explicitly as follows.

\pf{
The initial signed distance functions are defined as follows. For a sphere given by $\Lambda_1(\mathbf{c}, r)$:
\begin{equation}
u_{\text{sph}}^0(\mathbf{x}; \mathbf{c}, r) = |\mathbf{x} - \mathbf{c}| - r.
\end{equation}
For a cube given by $\Lambda_2(\mathbf{c}, r)$, let $d_i = |x_i - c_i| - r$ for $i \in \{1,2,3\}$, then the initial signed distance function is:
\begin{equation}
u_{\text{cub}}^0(\mathbf{x}; \mathbf{c}, r) = \min\left(\max_{i} d_i, 0\right) + \sqrt{\sum_{i=1}^3 \max(d_i, 0)^2}.
\end{equation}
}

\pf{
Next, to account for the external advective velocity field $\mathbf{v}$ in \eqref{eq:levelset}, we use for the translation $\mathbf{v}=(2,2,2)$ the backward-traced coordinates $\tilde{\mathbf{x}}(t) = \mathbf{x} - t \mathbf{v}$, and $\tilde{\mathbf{x}}(t) = (x_1 \cos(\omega t) \pfb{-} x_2 \sin(\omega t), \pfb{+}x_1 \sin(\omega t) + x_2 \cos(\omega t), x_3)^\top$ for the rotation $\mathbf{v}=(\omega x_2, -\omega x_1, 0)$, $\omega \in \mathbb{R}$. Note that for the advection in normal direction we add  to $\mathbf{v}$ the term $\delta \frac{\nabla u}{|\nabla u|}$ with $\delta \in \mathbb{R}$.
}

\pf{
Let $u_*^0$ denote either $u_{\text{sph}}^0$ or $u_{\text{cub}}^0$. The exact solutions $u_e(\mathbf{x}, t)$ of \eqref{eq:levelset} defined for $\mathbf{x} \in \mathbb{R}^3$ and $t>0$ are given as follows:\\
}

\begin{itemize}
\item \pf{
\textbf{Pure Advection} ($\delta = 0$): \ref{TS}, \ref{RS}, \ref{TC}, \ref{RC}
\begin{equation}
	\label{exfirst}
	u_e(\mathbf{x}, t) = u_*^0(\tilde{\mathbf{x}}(t); \mathbf{c}, r).
\end{equation}
}
\item \pf{
\textbf{Pure Shrinking} ($\delta < 0, \mathbf{v}=\mathbf{0}$): \ref{SS}, \ref{SC}
}
\pf{
\begin{equation}
	u_e(\mathbf{x}, t) = u_*^0(\mathbf{x}; \mathbf{c}, r) - \delta t.
\end{equation}
}
\item \pf{
\textbf{Combined Advection and Shrinking} ($\delta < 0, \mathbf{v} \neq \mathbf{0}$): \ref{RSS}, \ref{RSC}
\begin{equation}
	u_e(\mathbf{x}, t) = u_*^0(\tilde{\mathbf{x}}(t); \mathbf{c}, r) - \delta t.
\end{equation}
}
\item \pf{
\textbf{Pure Expansion} ($\delta > 0, \mathbf{v}=\mathbf{0}$): \ref{ES}
\begin{equation}
	\label{exlast}
	u_e(\mathbf{x}, t) = 
	\begin{cases}
		-r & \text{if } u_*^0(\mathbf{x}; \mathbf{c}, r) \le \delta t - r, \\
		u_*^0(\mathbf{x}; \mathbf{c}, r) - \delta t & \text{otherwise.}
	\end{cases}
\end{equation}
For the expanding cube \ref{EC} the exact solution $u_e$ is derived in Appendix~\ref{sec:app_EC}.
}
\end{itemize}

\pfb{Note that the exact solution $u_e(\mathbf{x},t)$ in~\eqref{exfirst}--\eqref{exlast} is the viscosity solution of~\eqref{eq:levelset}.} \nn{The exact solution is defined without a boundary condition: it is the solution
in the whole space. If restricted to $\Omega$, the dynamic interface (the zero level set) is not reaching
$\partial\Omega$. The Dirichlet condition uses the
value of exact solution on the inflow boundary, which has to be known there,
\begin{equation}
u_D(\mathbf{x},t) = u_e(\mathbf{x},t), \quad (\mathbf{x},t) \in \partial \Omega^{-} \times [0,T]
\end{equation}
and the eikonal boundary condition uses the property $|\nabla u|=1$ that it has there, which requires no values of the solution to be known.}

\begin{table}
\makebox[\linewidth][c]{%
\begin{minipage}{1.25\textwidth}
	\centering
	\begin{tabular}{lllc}
		\toprule
		Test case & $\Gamma_0(u)$ & $\mathbf{v}$ & $T$ \\
		\midrule
		\textbf{TS} \namelabel{TS}: A translation of the sphere &  $\Lambda_1((-0.5,-0.5,-0.5),0.5)$ & $(2,2,2)$ & $\pfb{0.5}$ \label{ex_TS} \\
		\cmidrule(lr){1-1} \cmidrule(lr){2-4}
		\textbf{RS} \namelabel{RS}: A rotation of the sphere & $\Lambda_1((0.625,0,0),0.5)$ & $(\pi x_2, -\pi x_1, 0)$ & $2$ \label{ex_RS} \\ 
		\cmidrule(lr){1-1} \cmidrule(lr){2-4}
		\textbf{ES} \namelabel{ES}: An expansion of the sphere & $\Lambda_1((0,0,0),0.5)$ & $\displaystyle \frac{\nabla u}{|\nabla u|}$ & $0.5$ \\
		\cmidrule(lr){1-1} \cmidrule(lr){2-4}
		\textbf{SS} \namelabel{SS}: A shrinking of the sphere & $\Lambda_1((0,0,0),1)$ & $\displaystyle -\frac{\nabla u}{|\nabla u|}$ & $0.5$ \\
		\cmidrule(lr){1-1} \cmidrule(lr){2-4}
		\textbf{TC} \namelabel{TC}: A translation of the cube & $\Lambda_2((-0.5,-0.5,-0.5),0.5)$ & $(2,2,2)$ & $0.5$ \\
		\cmidrule(lr){1-1} \cmidrule(lr){2-4}
		\textbf{RC} \namelabel{RC}: A rotation of the cube & $\Lambda_2((0.625,0,0),0.5)$ & $(\pi x_2, -\pi x_1, 0)$ & $2$ \\ 
		\cmidrule(lr){1-1} \cmidrule(lr){2-4}
		\textbf{EC} \namelabel{EC}: An expansion of the cube & $\Lambda_2((0,0,0),0.5)$ & $\displaystyle  \frac{\nabla u}{|\nabla u|}$ & $0.5$ \\ 
		\cmidrule(lr){1-1} \cmidrule(lr){2-4}
		\textbf{SC} \namelabel{SC}: A shrinking of the cube & $\Lambda_2((0,0,0),1)$ & $\displaystyle  -\frac{\nabla u}{|\nabla u|}$ & $0.5$ \\
		\cmidrule(lr){1-1} \cmidrule(lr){2-4}
		\textbf{RSS} \namelabel{RSS}: A rotation and shrinking of the sphere & $\Lambda_1((0.625,0,0),0.5)$ & $\displaystyle (\pi x_2, -\pi x_1, 0) -  0.1\frac{\nabla u}{|\nabla u|} $ & $1$ \\ 
		\cmidrule(lr){1-1} \cmidrule(lr){2-4}
		\textbf{RSC} \namelabel{RSC}:A rotation and shrinking of the cube &  $\Lambda_2((0.625,0,0),0.5)$ & $\displaystyle (\pi x_2, -\pi x_1, 0)-0.1\frac{\nabla u}{|\nabla u|}$ & $1$ \\
		\bottomrule
	\end{tabular}
\end{minipage}}
\vspace{0.5em}
\begin{center}
\caption{The $8$ examples from the top are tested in Section~\ref{sec:comp_EDBC} for the comparison to Dirichlet boundary condition. The last two examples are used in Section~\ref{sec:comp_other} for the comparison to linearly extended and zero Neumann boundary condition.}\label{tab:test_case}	
\end{center}	
\end{table}

\begin{figure}
\makebox[\linewidth][c]{%
\begin{minipage}{1.25\textwidth}
	\begin{center}
		\begin{tabular}{cccc}
			\includegraphics[height=3.7cm]{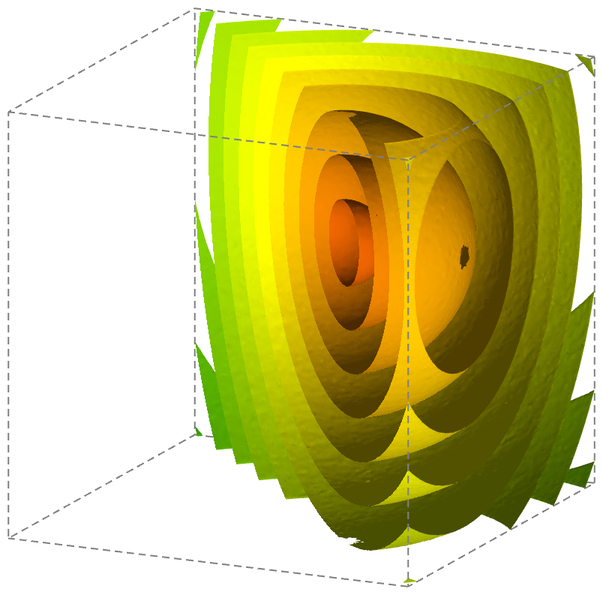} &
			\includegraphics[height=3.7cm]{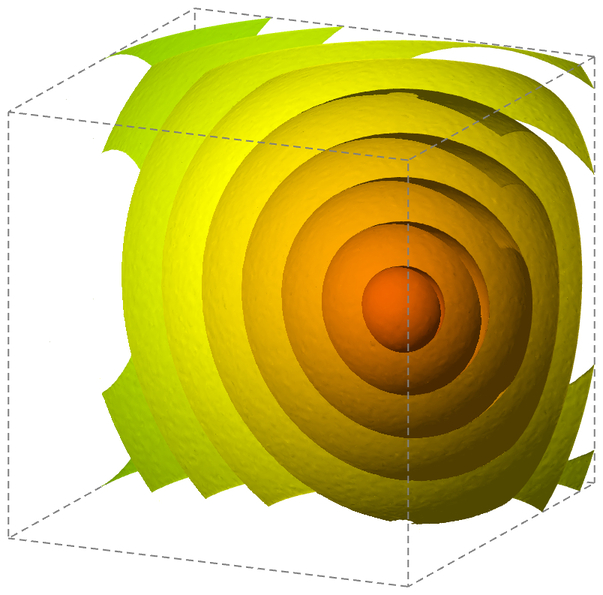} &			
			\includegraphics[height=3.7cm]{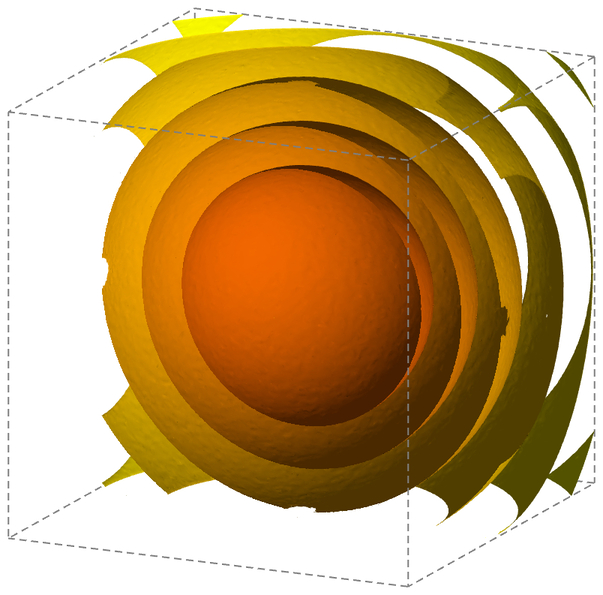} &
			\includegraphics[height=3.7cm]{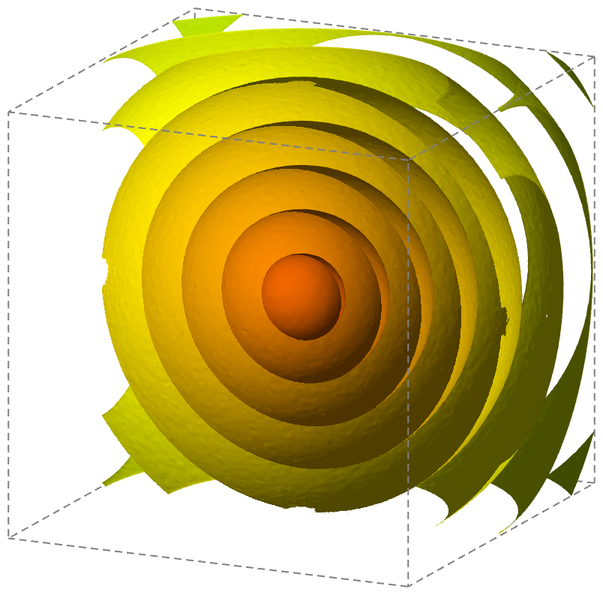} \\
			\ref{TS} & \ref{RS} & \ref{ES} & \ref{SS} \\
			\includegraphics[height=3.7cm]{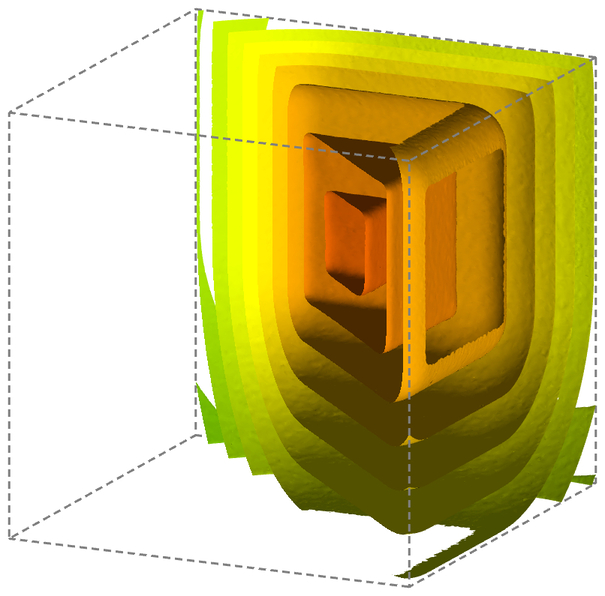} &
			\includegraphics[height=3.7cm]{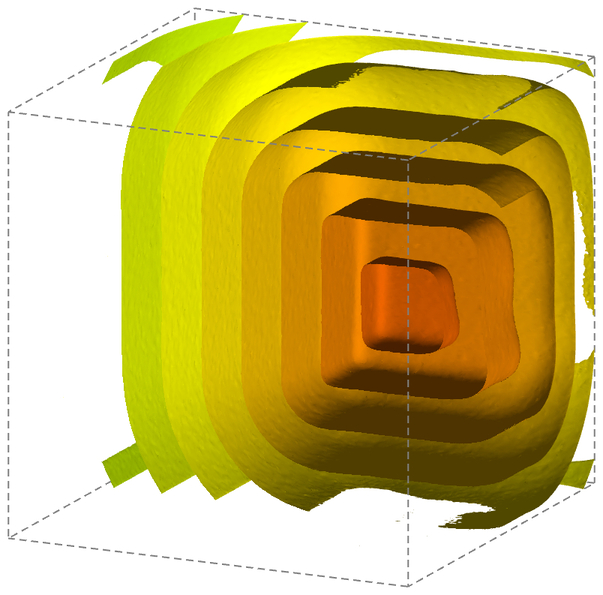} &
			\includegraphics[height=3.7cm]{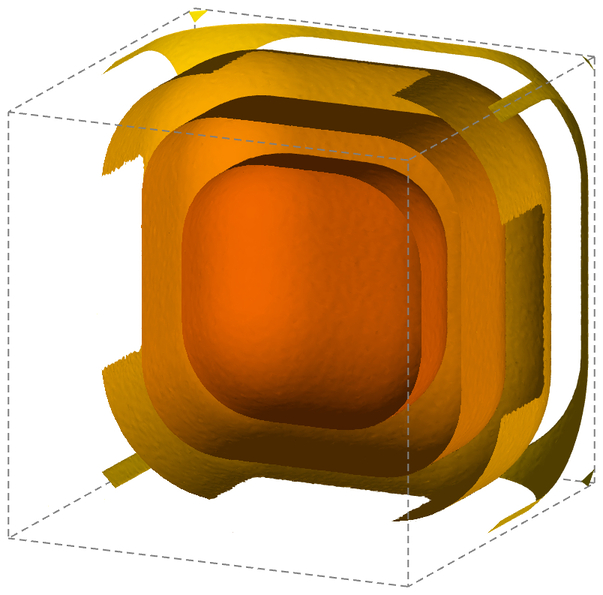} &
			\includegraphics[height=3.7cm]{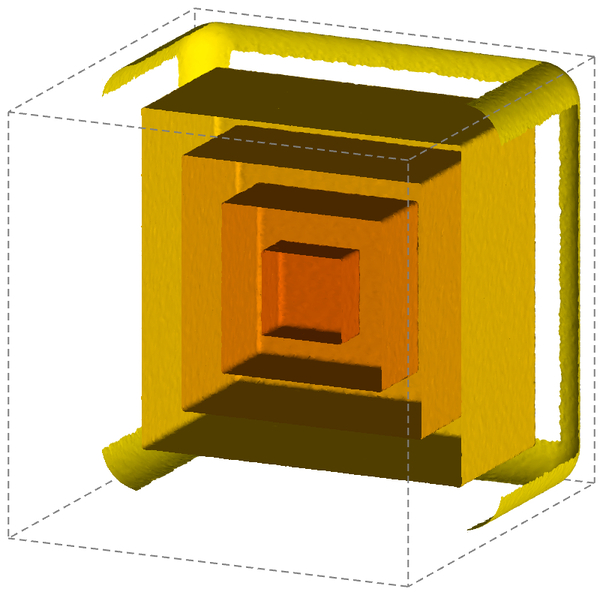} \\
			\ref{TC} & \ref{RC} & \ref{EC} & \ref{SC} 
		\end{tabular}
	\end{center}
\end{minipage}}
\begin{center}
\caption{Using the eikonal boundary condition, the \nn{equally spaced} isosurfaces~\eqref{eq:eqdist_iso} of test cases in Table~\ref{tab:test_case} at the final time $T$ are presented on the mesh $\mathcal{M}_{\text{M}}$ in Table~\ref{tab:meshes}. The second smallest surface is the evolved surface $\Gamma_{T}(u)$~\eqref{eq:zerolevel}. In \ref{TS} and \ref{TC}, we diagonally cut to present the isosurfaces for proper visualization.} \label{fig:EKBC}
\end{center}
\end{figure}
In Figure~\ref{fig:EKBC}, for the test cases above, we present the numerical solutions at the final time $T$ of using the eikonal boundary condition as the isosurfaces of equally spaced values:
\begin{equation}\label{eq:eqdist_iso}
\{\mathbf{x}:u(\mathbf{x},T) = 0.25(l-2), \ l = 1,2,\ldots\}
\end{equation}
The second smallest surface is the result of surface evolution $\Gamma_{T}(u)$ in~\eqref{eq:zerolevel} from $\Gamma_{0}(u)$. All isosurfaces have almost no distortion qualitatively.

\begin{figure}
\makebox[\linewidth][c]{%
\begin{minipage}{1.25\textwidth}
	\begin{center}
		\begin{tabular}{cccc}
			\includegraphics[width=3.8cm]{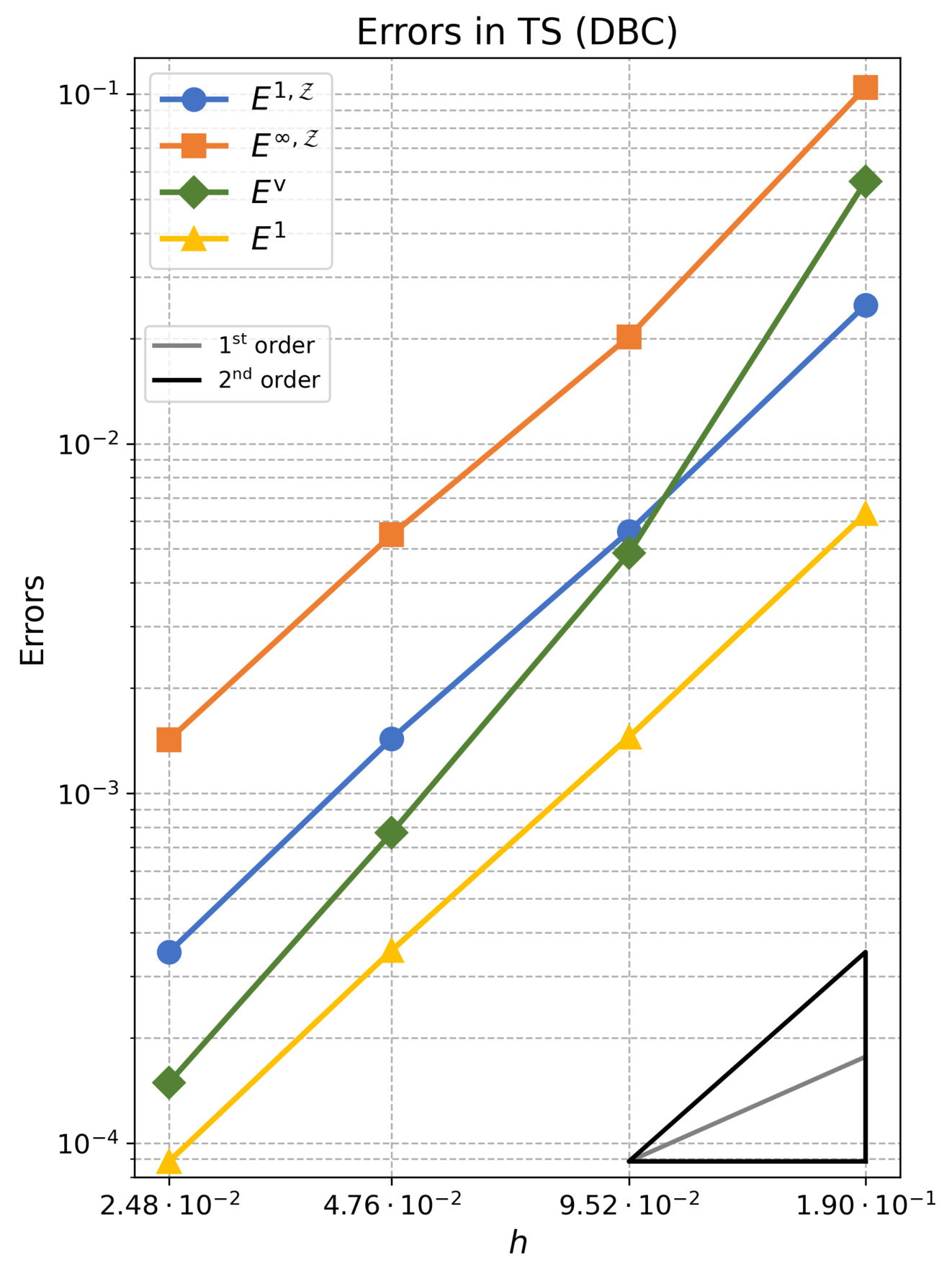} & 			
			\includegraphics[width=3.8cm]{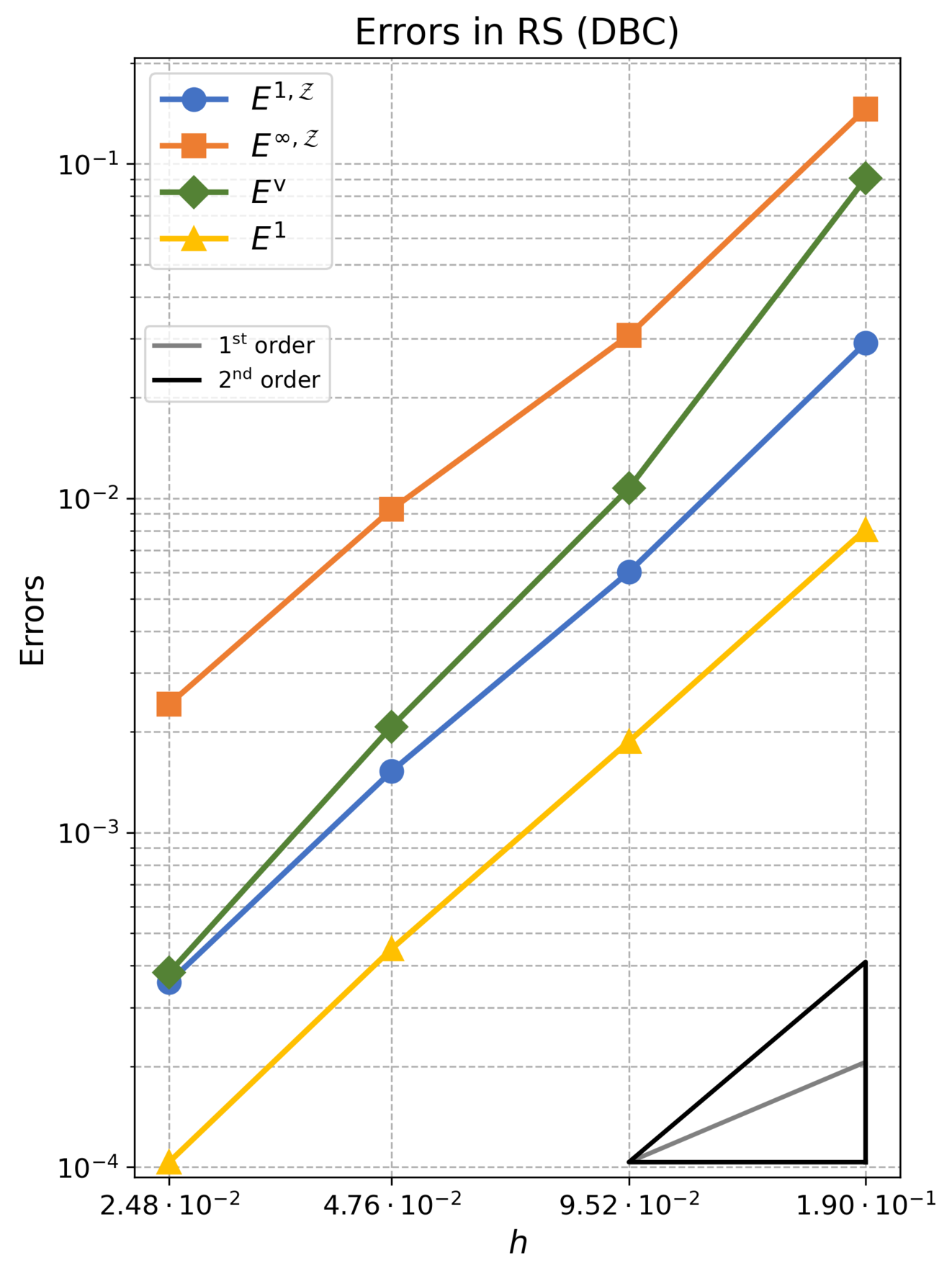} &
			\includegraphics[width=3.8cm]{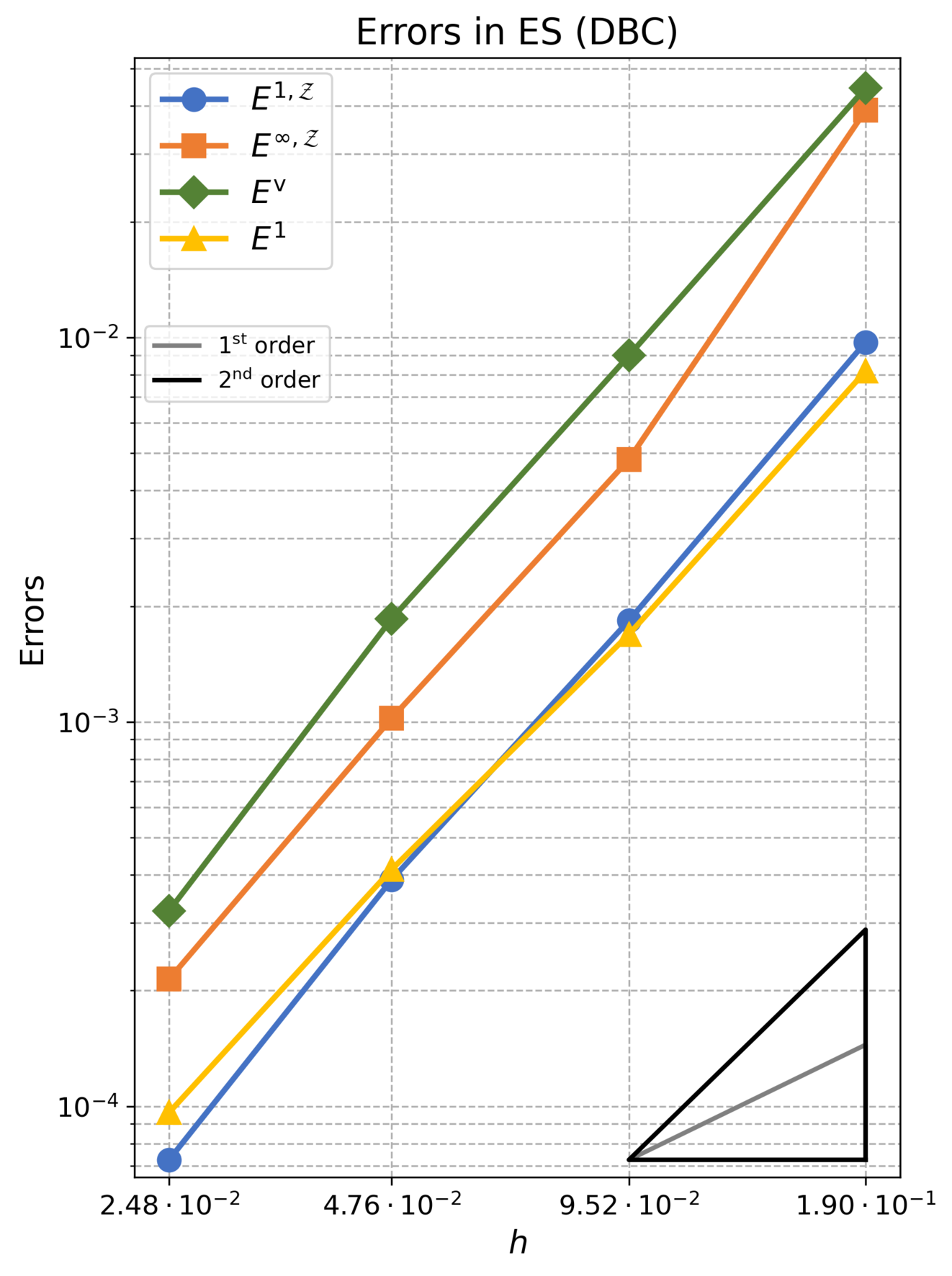} &
			\includegraphics[width=3.8cm]{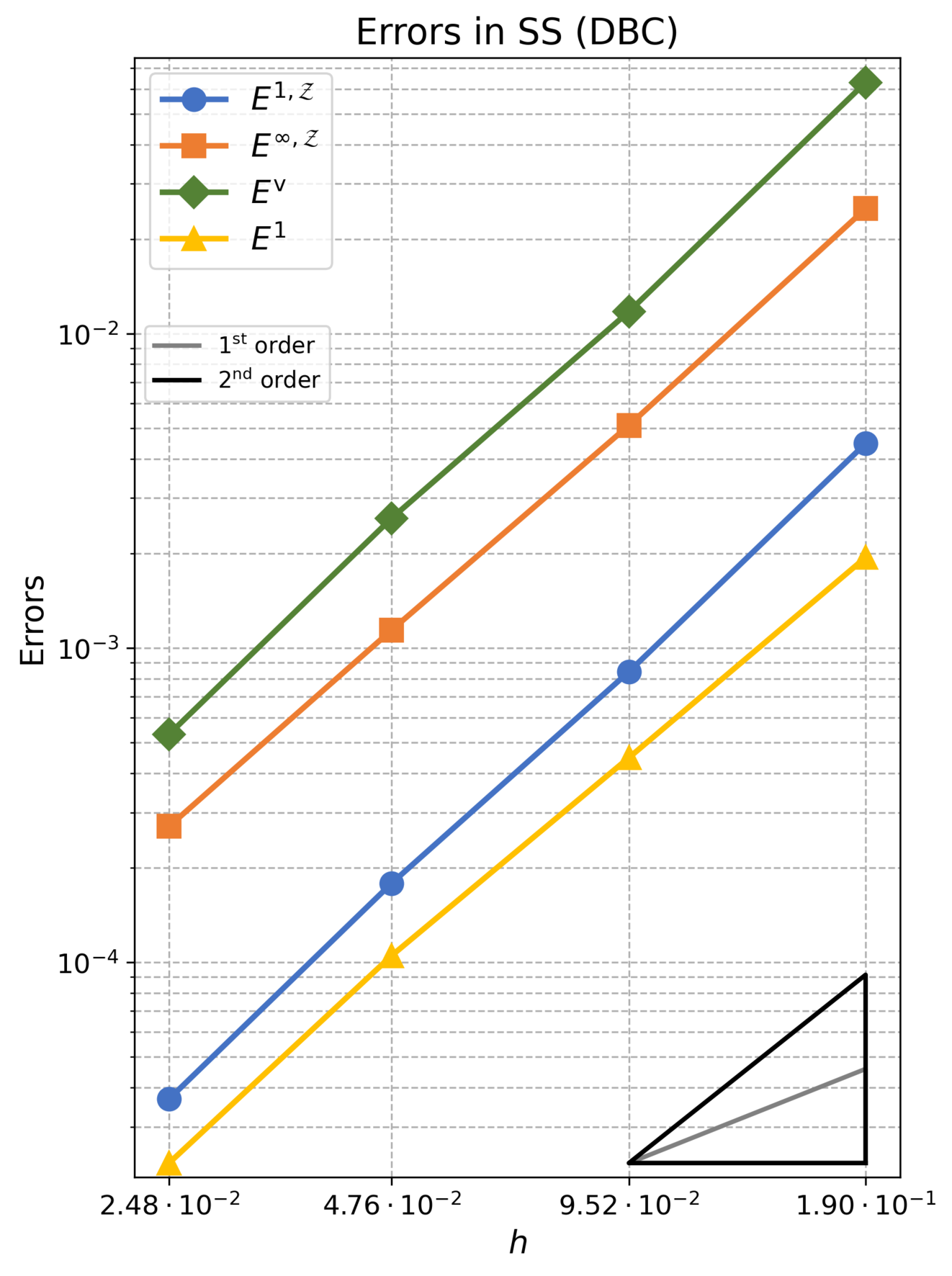} \\
			\includegraphics[width=3.8cm]{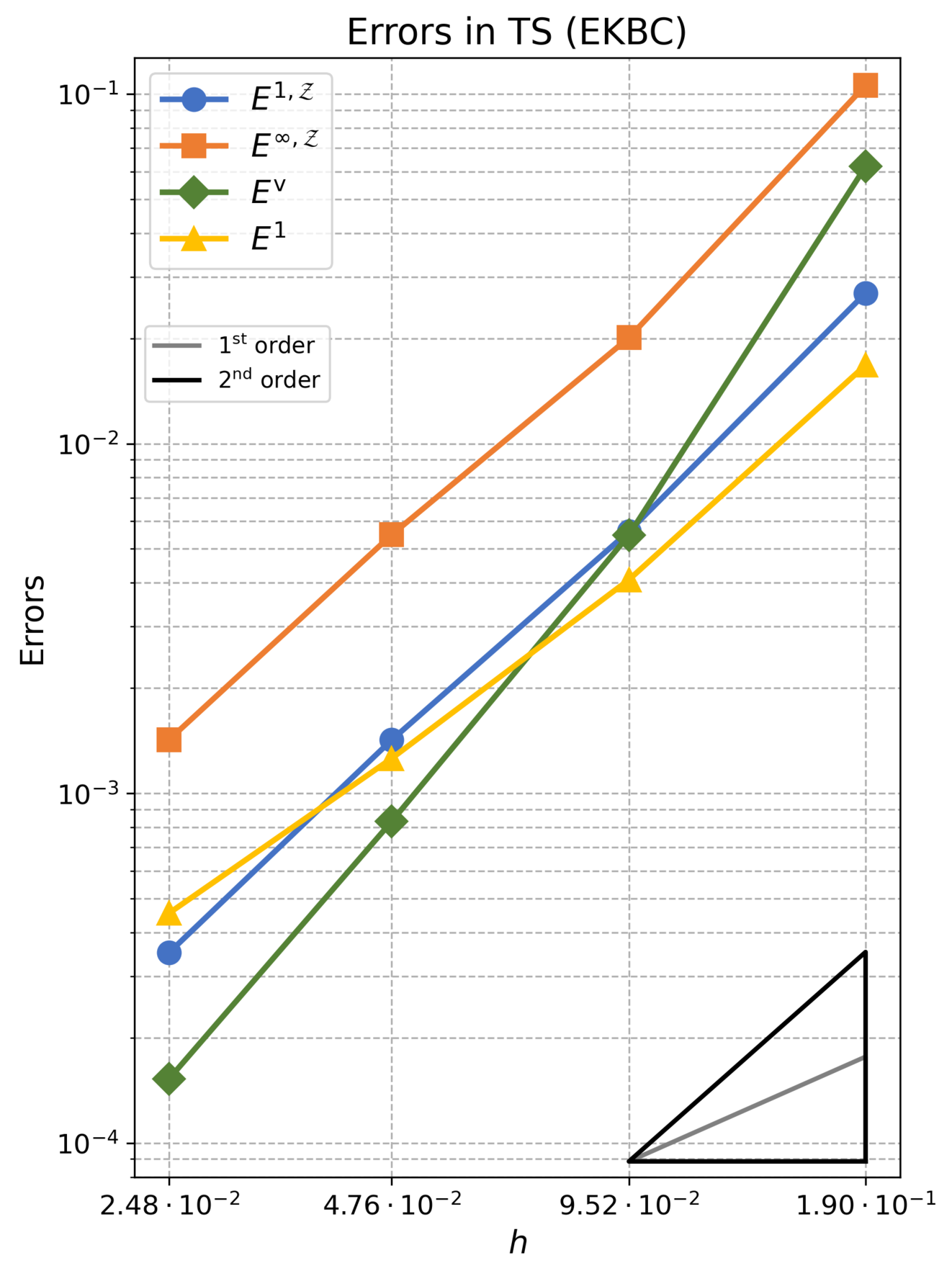} &
			\includegraphics[width=3.8cm]{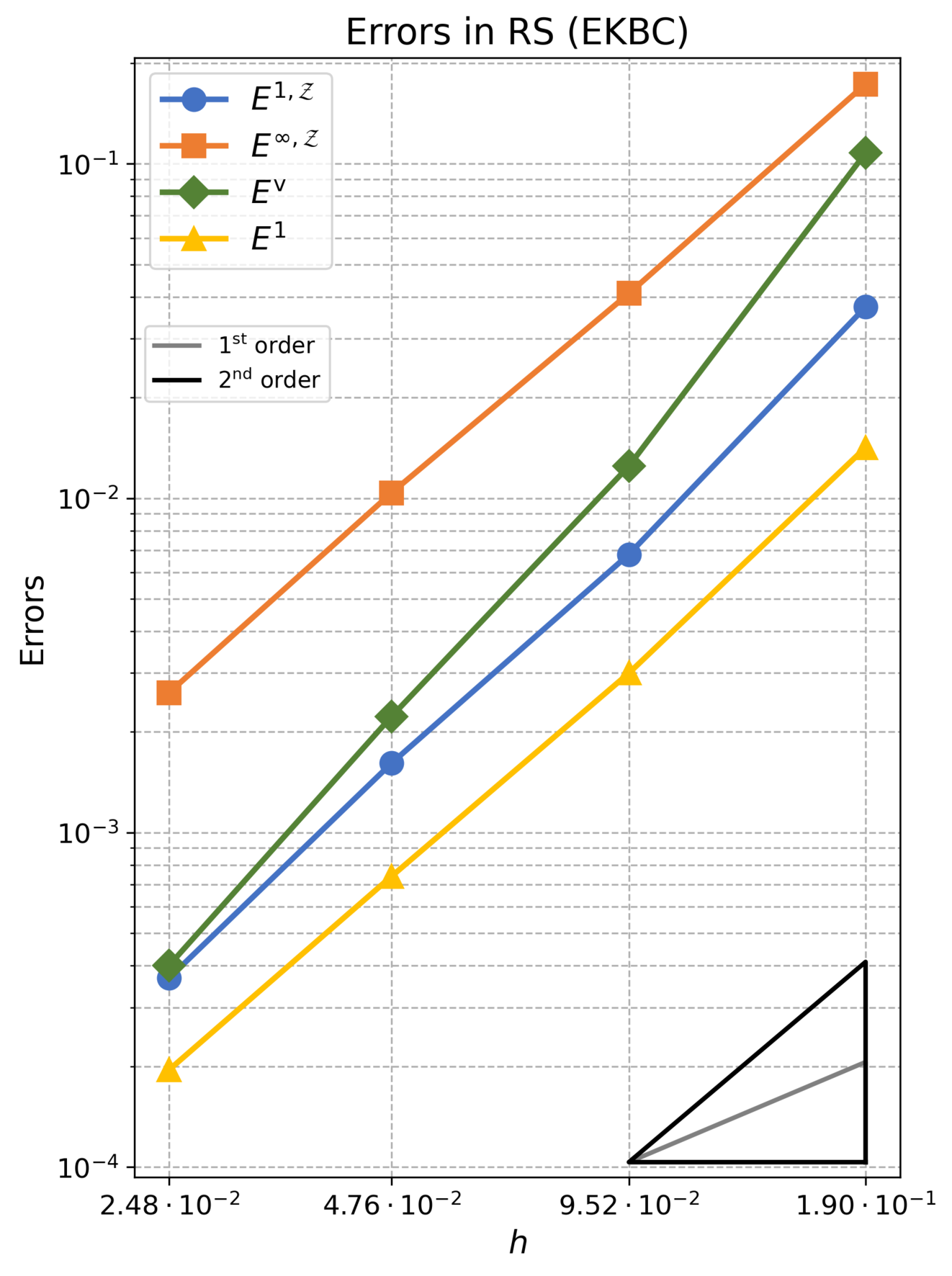} &
			\includegraphics[width=3.8cm]{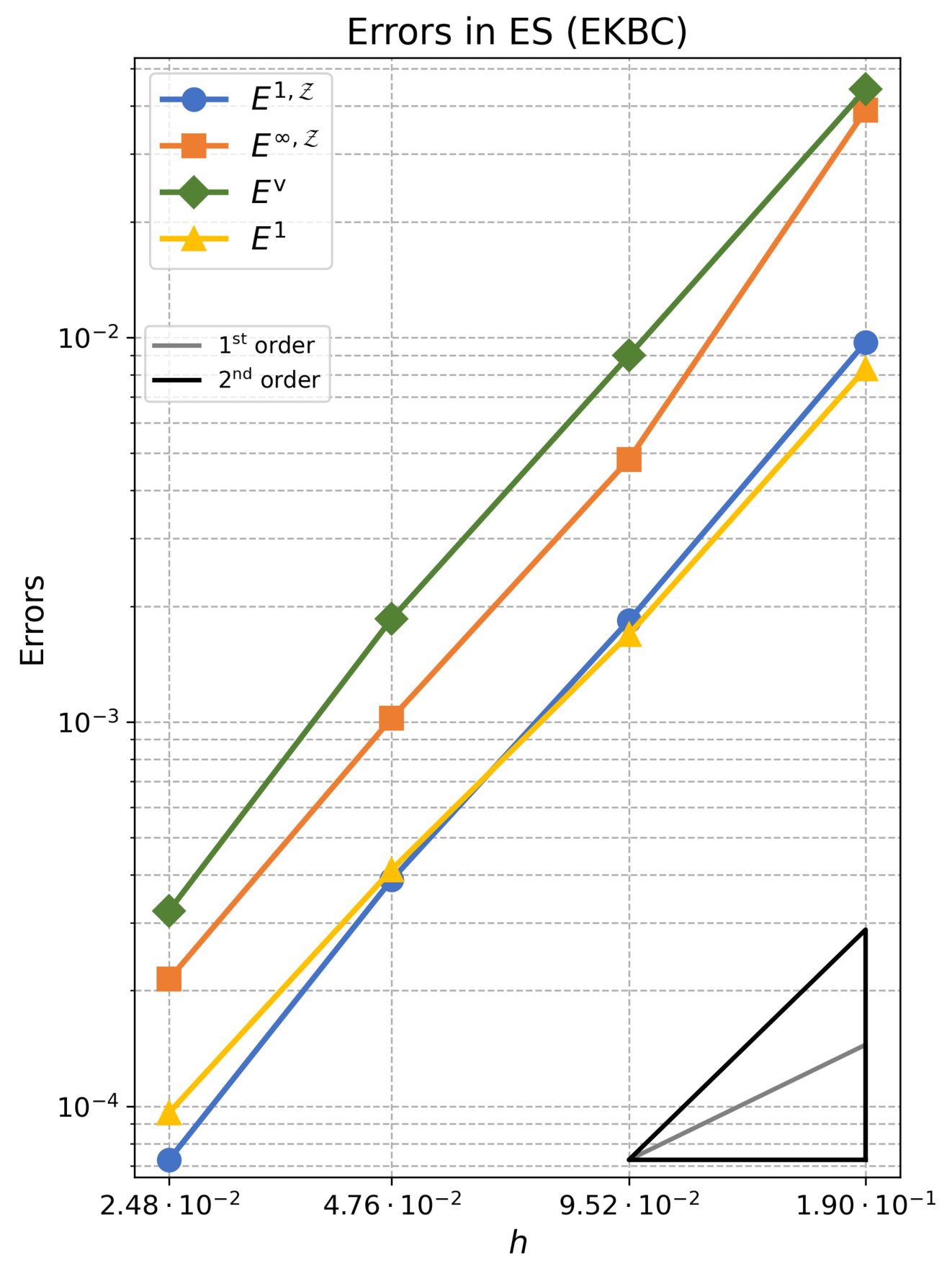} &			
			\includegraphics[width=3.8cm]{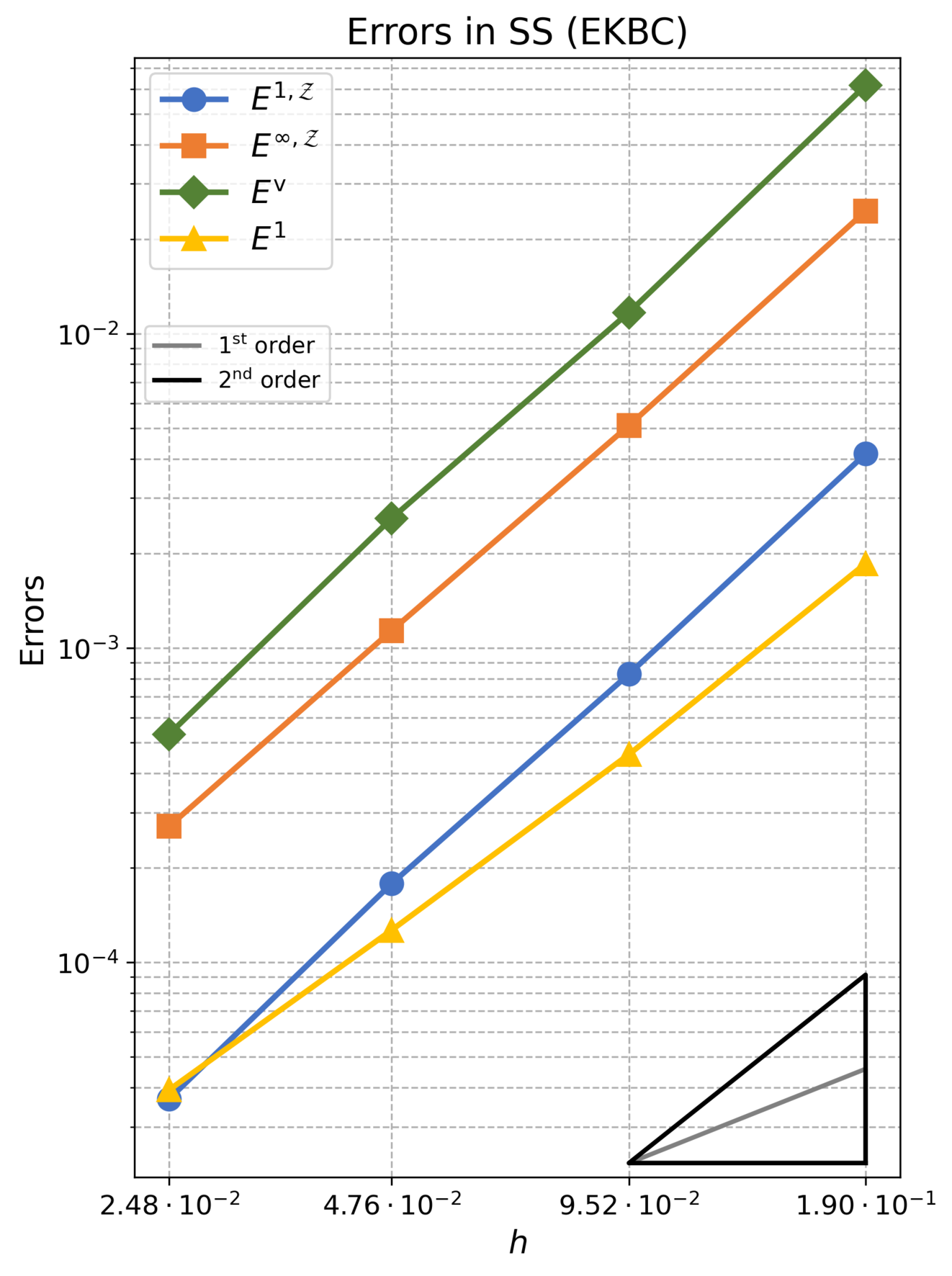} 
		\end{tabular}
	\end{center}
\end{minipage}}
\begin{center}
\caption{Applying $\triangle t_{\text{M}}$~\eqref{eq:time_step}, the characteristic length versus errors $E^{1,\mathcal{Z}}$~\eqref{eq:err_L1_zero}, $E^{\infty,\mathcal{Z}}$~\eqref{eq:err_Linf_zero}, $E^{\text{v}}$~\eqref{eq:err_vol}, $E^1$~\eqref{eq:err_L1} for test cases~\ref{TS},~\ref{RS},~\ref{ES}, and~\ref{SS} are presented by using the Dirichlet (upper row) and the eikonal boundary condition (lower row).} \label{fig:conv_order_S}	
\end{center}
\end{figure}

\begin{figure}
\makebox[\linewidth][c]{%
\begin{minipage}{1.25\textwidth}
	\begin{center}
		\begin{tabular}{cccc}
			\includegraphics[width=3.8cm]{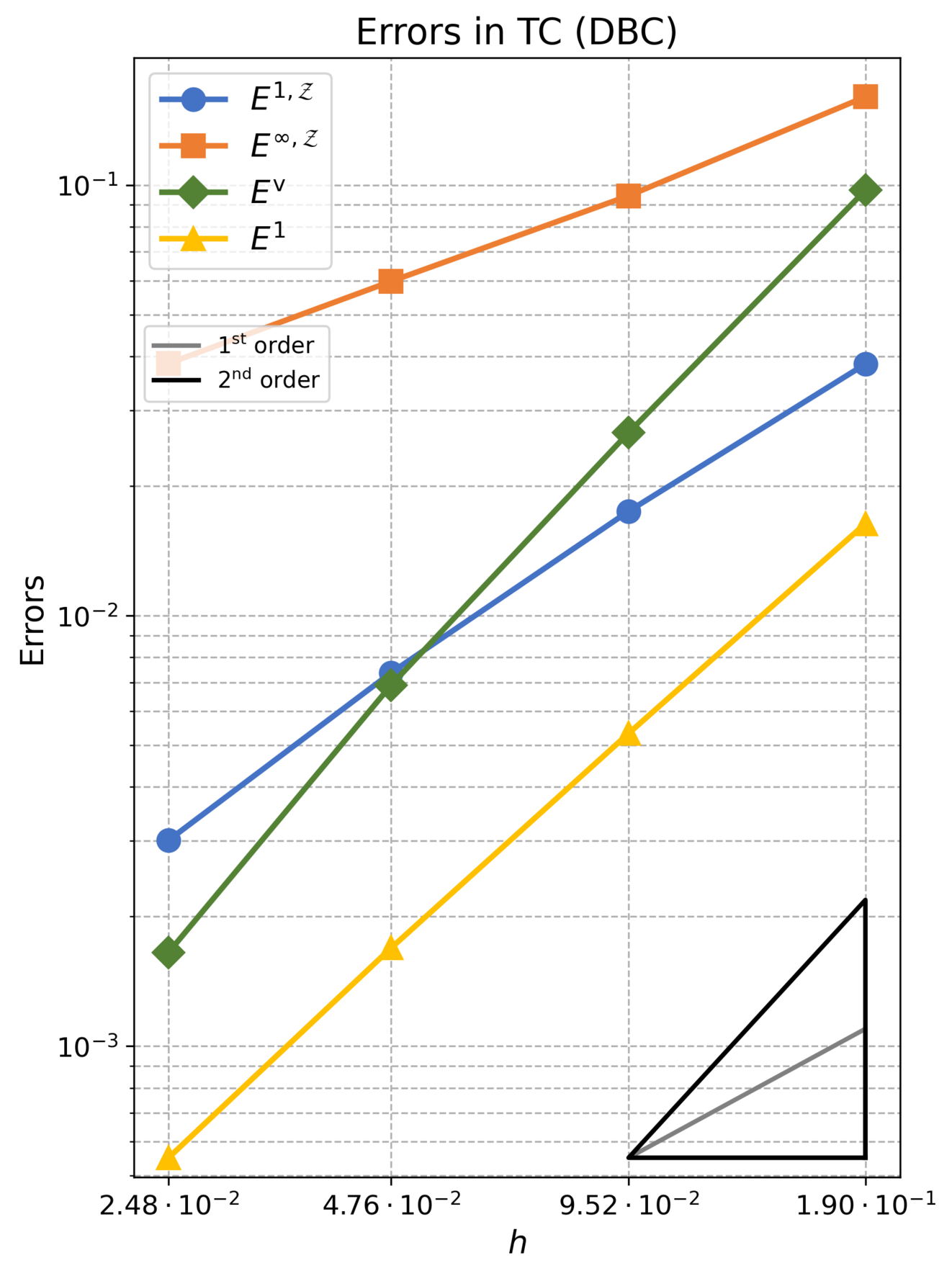} & 
			\includegraphics[width=3.8cm]{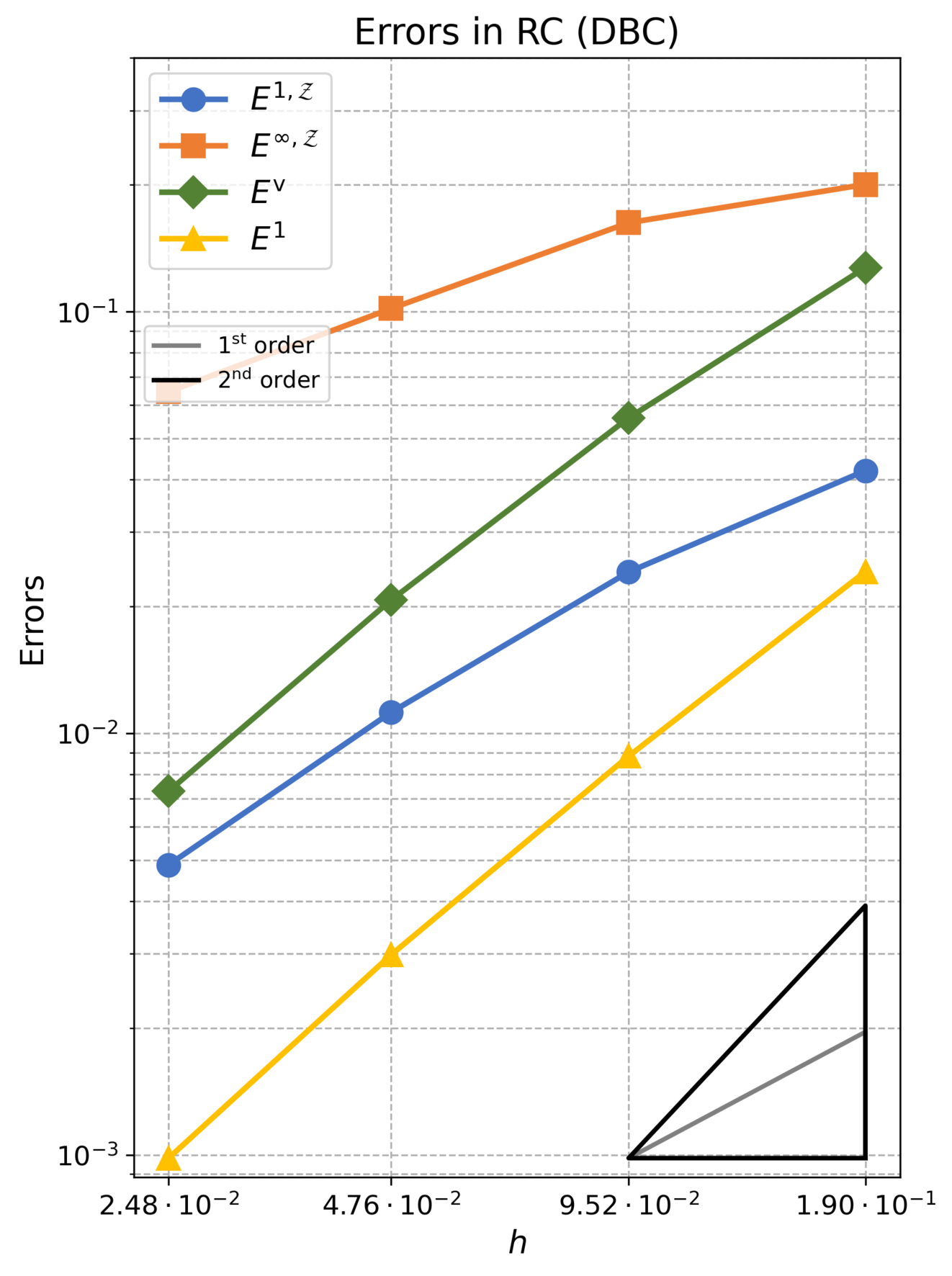} &
			\includegraphics[width=3.8cm]{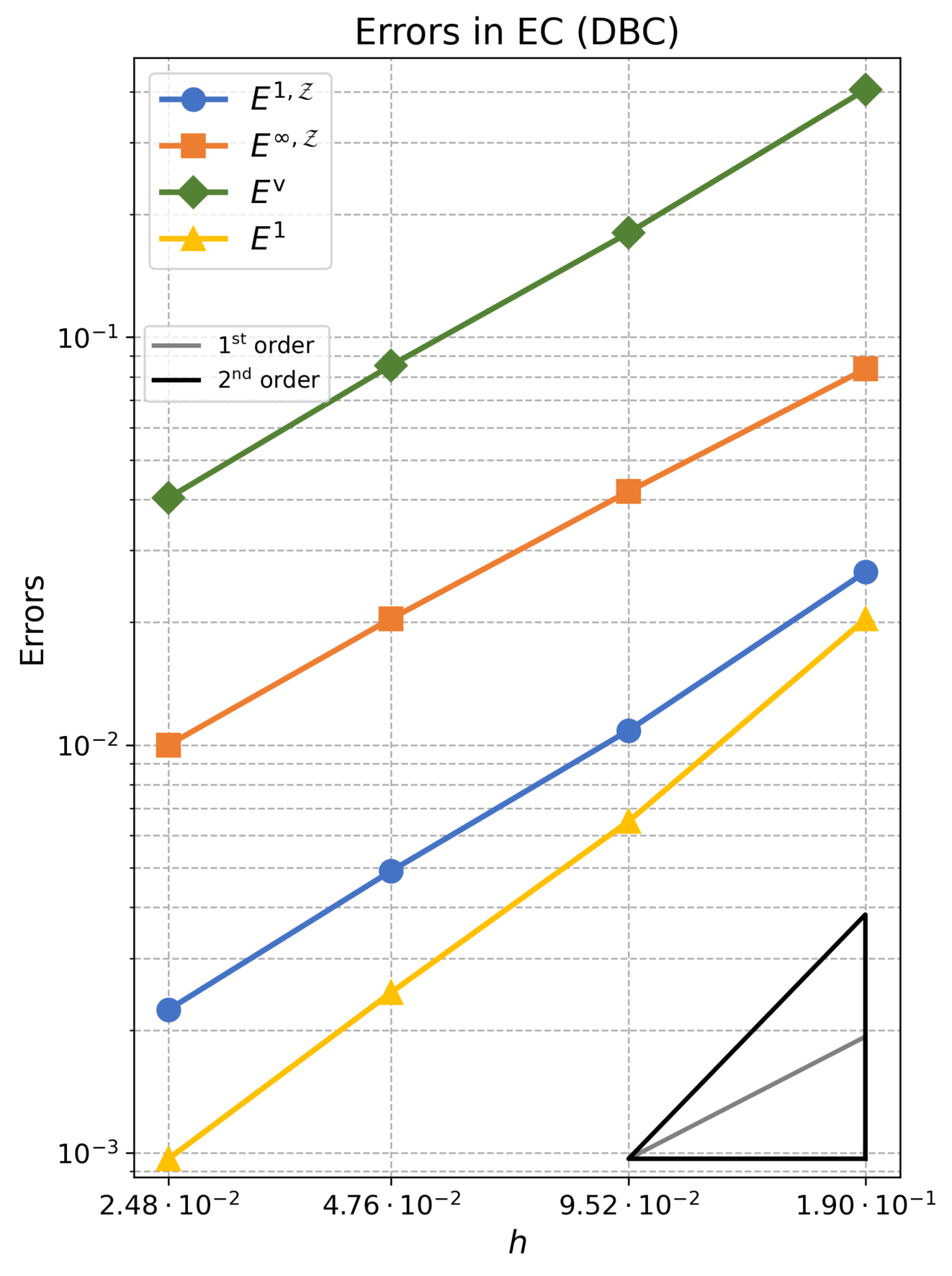} &
			\includegraphics[width=3.8cm]{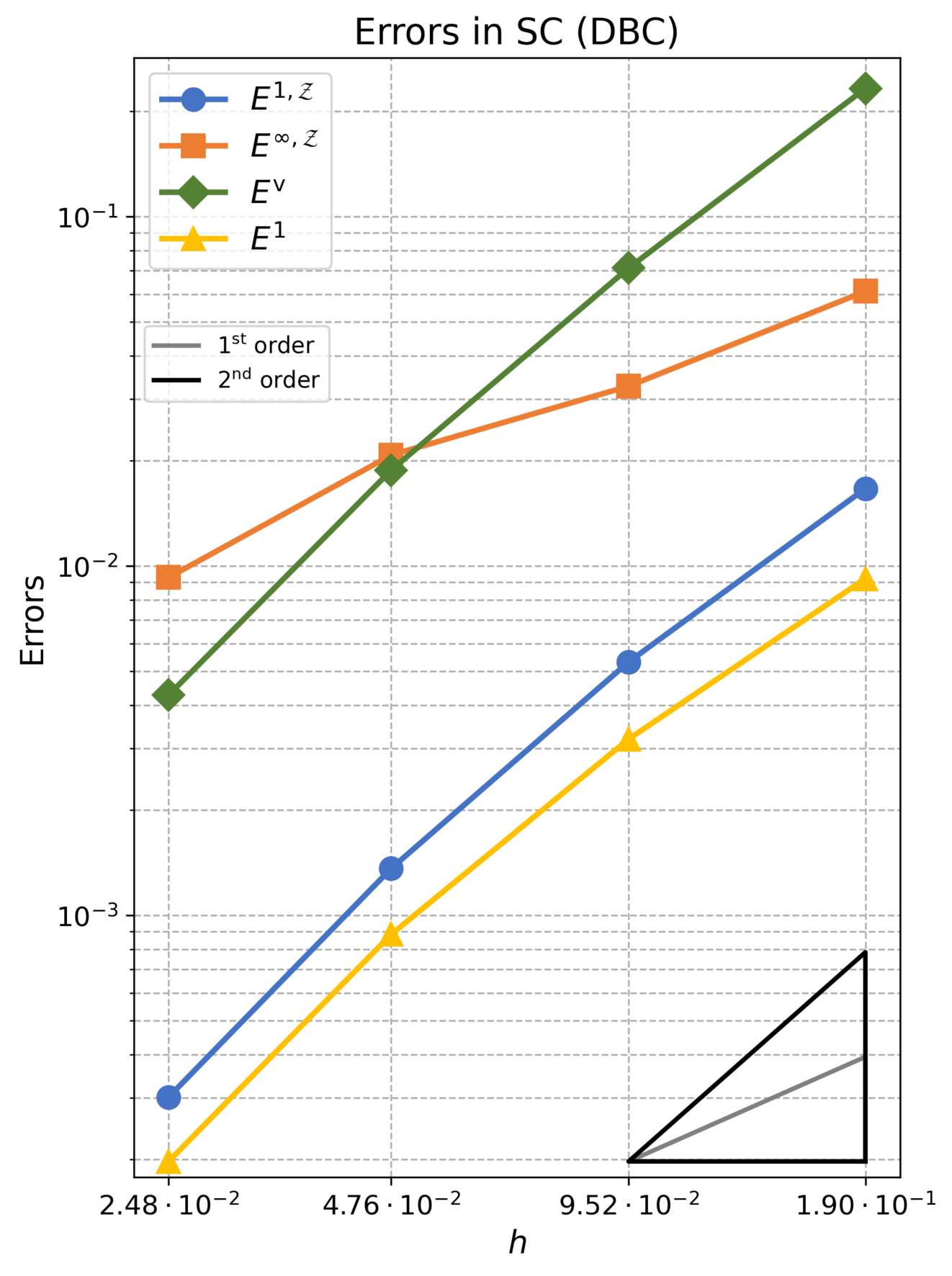} \\
			\includegraphics[width=3.8cm]{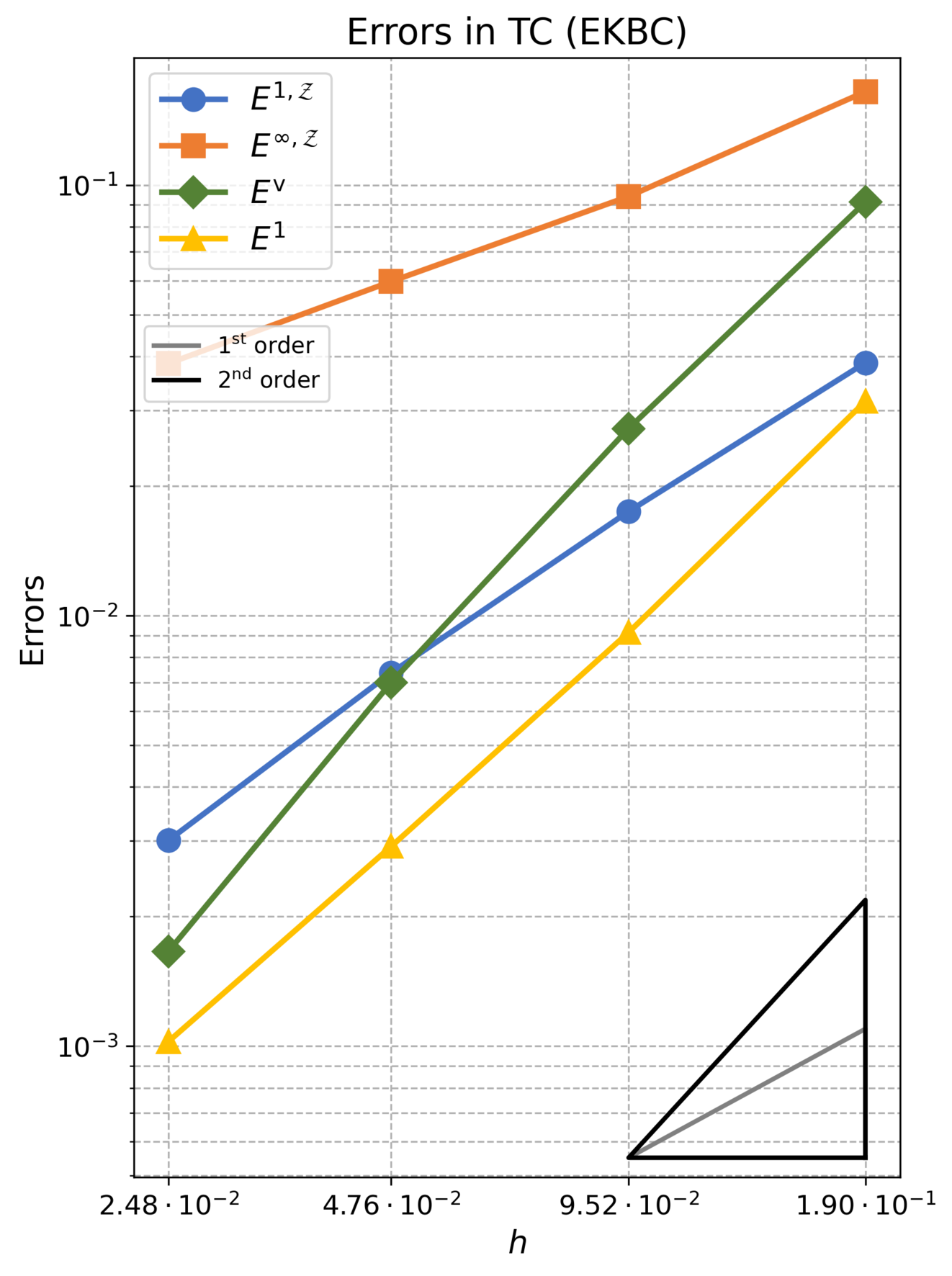} &
			\includegraphics[width=3.8cm]{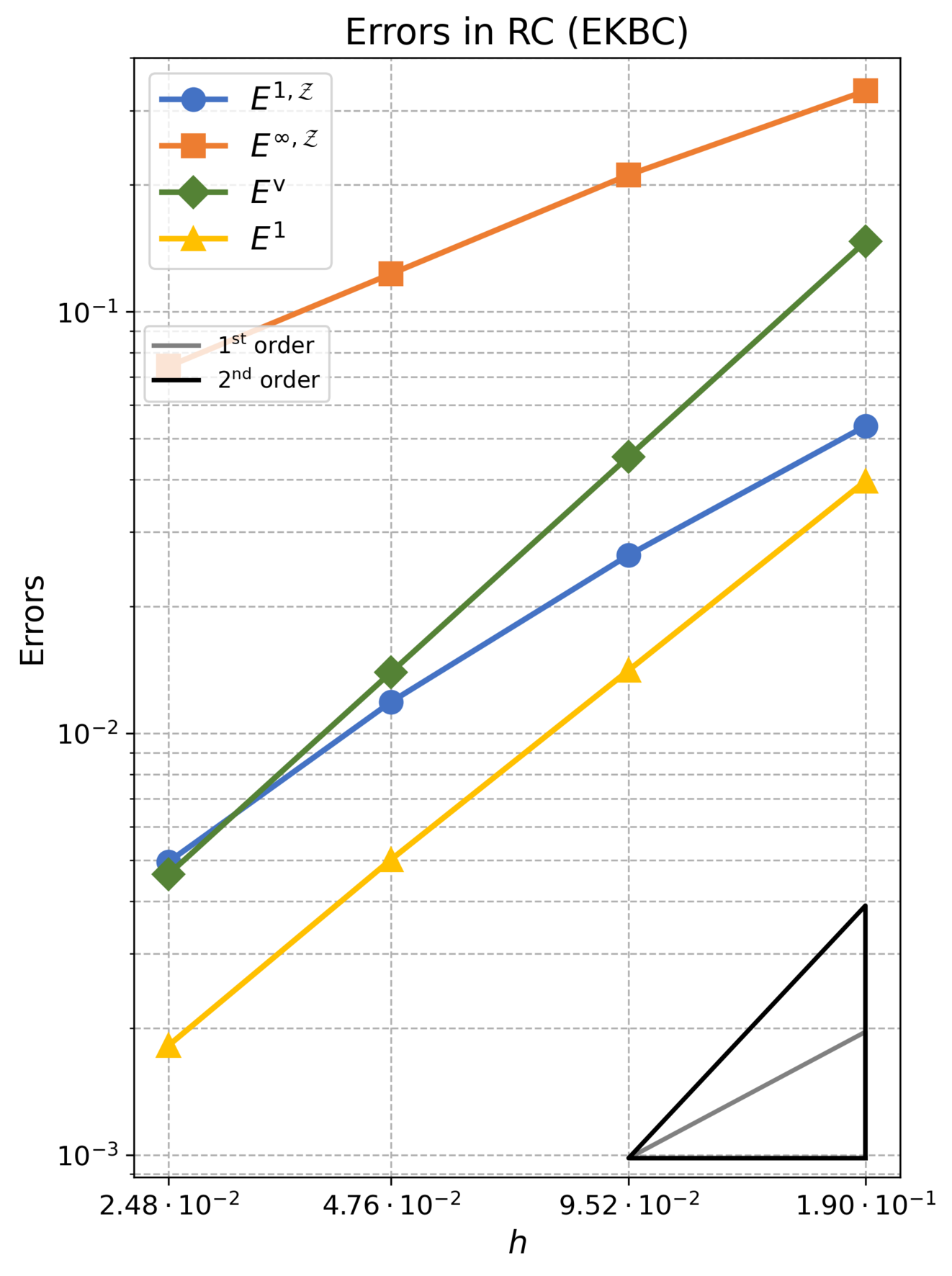} &
			\includegraphics[width=3.8cm]{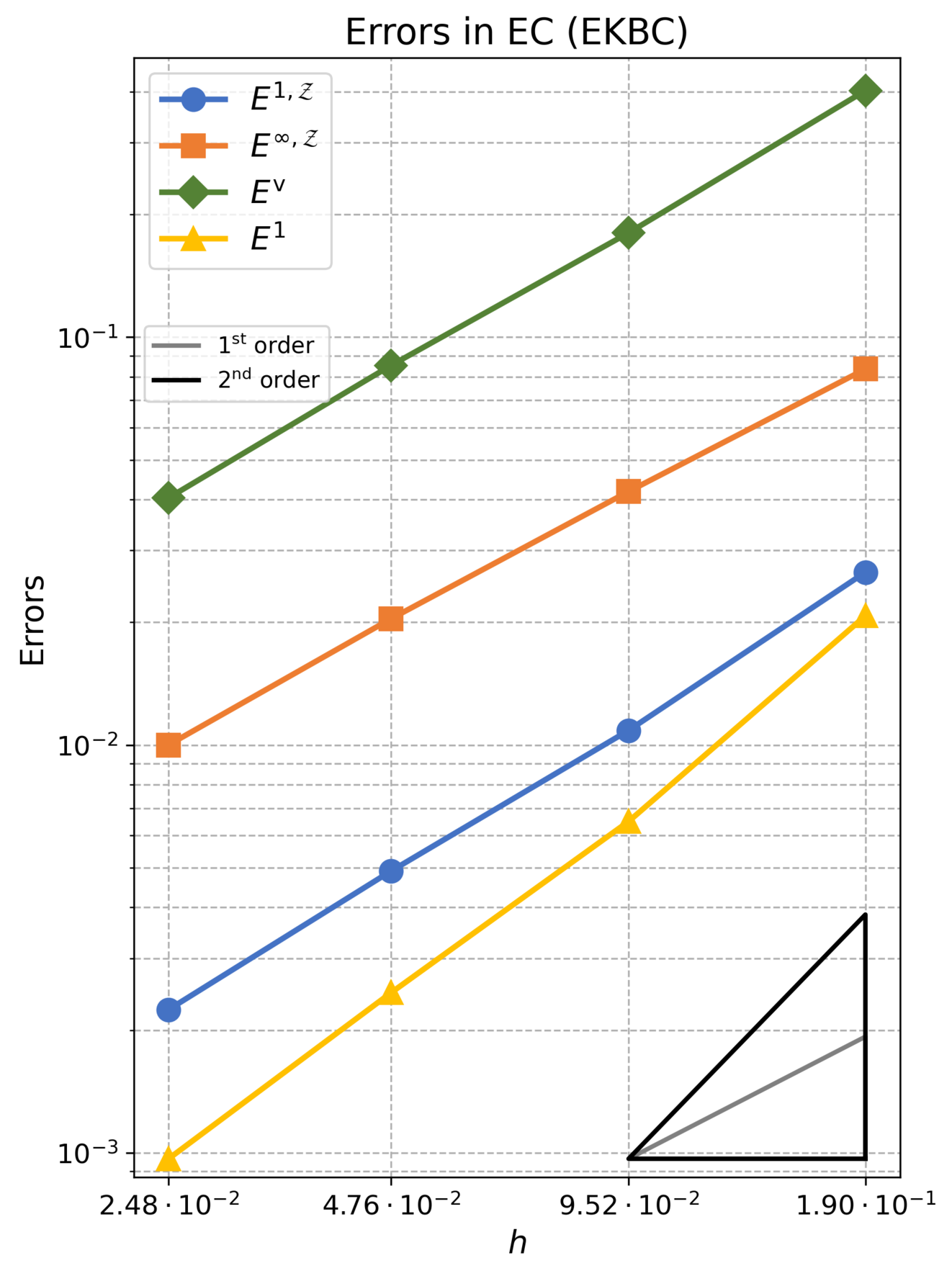} &
			\includegraphics[width=3.8cm]{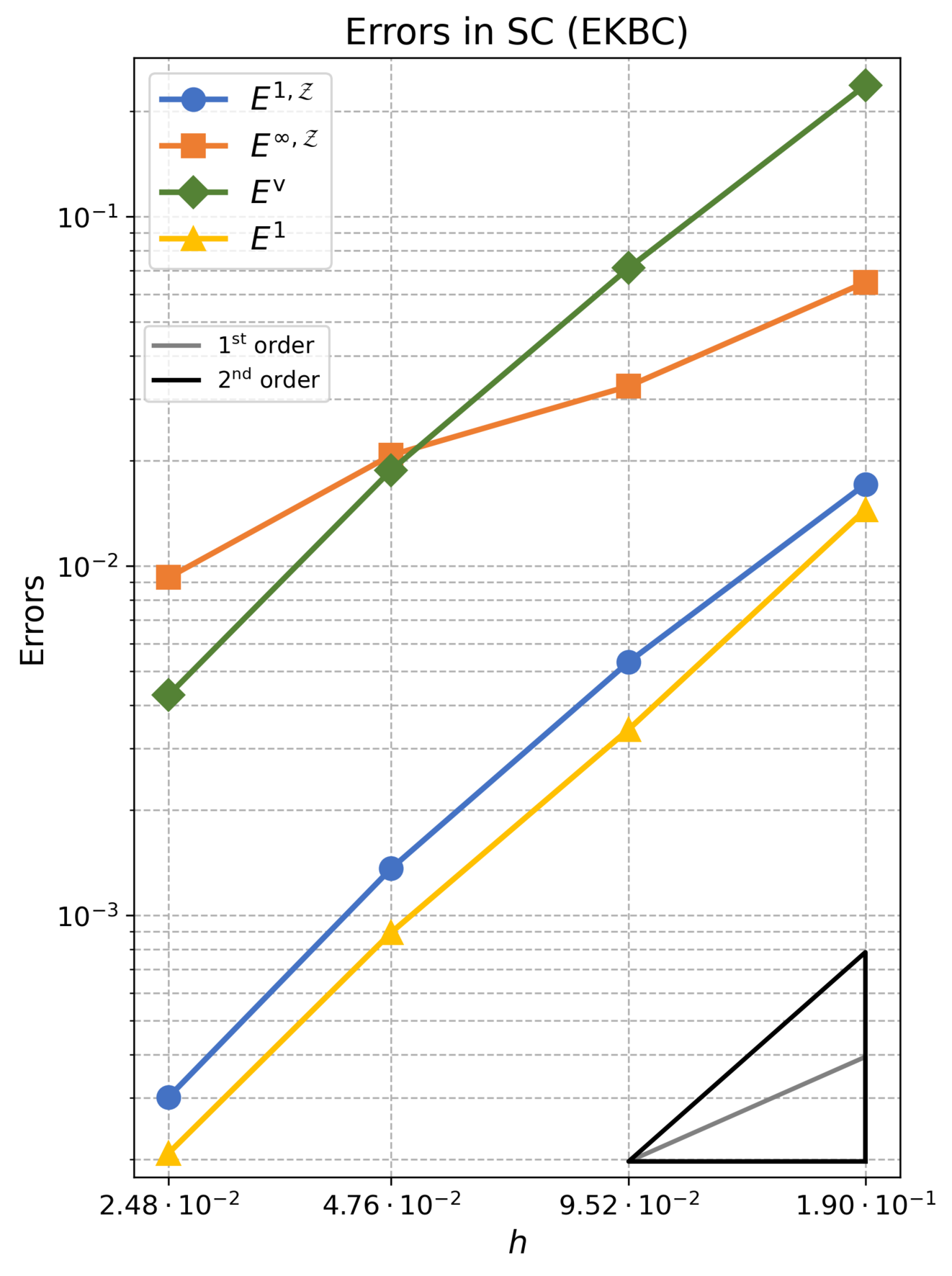}
		\end{tabular}
	\end{center}
\end{minipage}}
\begin{center}
\caption{Applying $\triangle t_{\text{M}}$~\eqref{eq:time_step}, the characteristic length versus errors $E^{1,\mathcal{Z}}$~\eqref{eq:err_L1_zero}, $E^{\infty,\mathcal{Z}}$~\eqref{eq:err_Linf_zero}, $E^{\text{v}}$~\eqref{eq:err_vol}, $E^1$~\eqref{eq:err_L1} for test cases~\ref{TC},~\ref{RC},~\ref{EC}, and~\ref{SC} are presented by using the Dirichlet (upper row) and the eikonal boundary condition (lower row).} \label{fig:conv_order_C}
\end{center}
\end{figure}

The characteristic length~\eqref{eq:char_len} versus errors, $E^{1,\mathcal{Z}}$~\eqref{eq:err_L1_zero}, $E^{\infty,\mathcal{Z}}$~\eqref{eq:err_Linf_zero}, and $E^{\text{v}}$~\eqref{eq:err_vol}, $E^1$~\eqref{eq:err_L1}, are presented by the log-log scaled graph in Figures~\ref{fig:conv_order_S} and~\ref{fig:conv_order_C}. The slopes in the graphs are the experimental order of convergence~\eqref{eq:EOC}; see also values in Table~\ref{tab:comp_EDBC}. For a visual comparison, the first and second orders are shown by the slopes of the hypotenuse of the triangles in all figures. All graphs on the upper and lower rows are the results of using the Dirichlet and eikonal boundary condition, respectively. The main point of graphs is that the proposed method using the eikonal boundary condition can achieve competitive results using the Dirichlet boundary condition \pf{given by exact solutions $u_e(\mathbf{x},t)$ defined in \eqref{exfirst} - \eqref{exlast}}\nn{, and in Appendix~\ref{sec:app_EC} for the expanding cube}. The behavior of the graphs on the upper rows is similar to that of the graphs on the lower rows, particularly in terms of their slops. More specifically, the $EOC$ from the errors, $E^{1,\mathcal{Z}}$, $E^{\infty,\mathcal{Z}}$, and $E^{\text{v}}$, are nearly the same. The proposed method shows lower values of $EOC_{E^{1}_{\text{M=3}}}$ in the test case~\ref{TS}.

\begin{table}
\makebox[\linewidth][c]{%
\begin{minipage}{1.25\textwidth}
	\centering
	\begin{tabular}{cS[table-format=1]
			S[table-format=2.4]
			S[table-format=2.4]
			S[table-format=2.4]}
		\toprule
		\multicolumn{1}{c}{} & \multicolumn{1}{c}{M} & \multicolumn{1}{c}{$CFL^{\text{min}}_{\text{M}}$}    & \multicolumn{1}{c}{$CFL^{\text{ave}}_{\text{M}}$} & \multicolumn{1}{c}{$CFL^{\text{max}}_{\text{M}}$} \\
		\midrule
		\multirow{4}{*}{A} & 1 & 0.8227 & 2.0317 & 5.2842 \\ 
		& 2 & 0.8651 & 2.0915 & 7.5764 \\ 
		& 3 & 0.9156 & 2.1796 & 6.3730 \\ 
		& 4 & 0.9664 & 2.0240 & 7.8618 \\
		\midrule
		\multirow{4}{*}{B} & 1 & 0.0129 & 2.1984 & 7.8992 \\ 
		& 2 & 0.0092 & 2.2319 & 12.0948 \\ 
		& 3 & 0.0067 & 2.3052 & 9.0246 \\ 
		& 4 & 0.0001 & 2.0281 & 11.0649 \\ 
		\midrule
		\multirow{4}{*}{C} & 1 & 0.2375 & 0.5865 & 1.5254 \\ 
		& 2 & 0.2497 & 0.6038 & 2.1871 \\ 
		& 3 & 0.2643 & 0.6292 & 1.8397 \\ 
		& 4 & 0.2790 & 0.5843 & 2.2695 \\ 
		\bottomrule
	\end{tabular}
\end{minipage}}
\vspace{0.5em}
\begin{center}
\caption{The $CFL$s~\eqref{eq:CFL} with $\triangle t_{\text{M}}$~\eqref{eq:time_step} are presented for the test cases A=~\ref{TS} or~\ref{TC}, B=~\ref{RS} or~\ref{RC}, and C=~\ref{ES},~\ref{SS},~\ref{EC}, or~\ref{SC}. For all test cases on the mesh $\mathcal{M}_{\text{M}}$, the same time step $\triangle t_{\text{M}}$ is used, while $CFL$ varies due to differences in velocity and characteristic length.\label{tab:CFL}}
\end{center}
\end{table}

The $CFL$ number in Table~\ref{tab:CFL} is presented to show that the proposed method is robust to the size of the time step to solve advective or normal flow equations in~\eqref{eq:levelset}. In Table~\ref{tab:meshes}, the difference between the smallest and the largest characteristic length in the cube domain is roughly the order of $10$ regardless of the level $\text{M}$. It means that the corresponding $CFL$s can be the similar order if a constant velocity is used:
\begin{align}\label{eq:CFL}
CFL^{\text{min}}_{\text{M}} &= \min_{p \in \mathcal{I}_{\text{M}}}  \frac{|\mathbf{v}(\mathbf{x}_p)|\triangle t_{\text{M}}}{|\Omega_p|_{B}^{\frac{1}{3}}}, \\	
CFL_{\text{M}}^{\text{ave}} &= \frac{1}{|\mathcal{I}_{\text{M}}|} \sum_{p \in \mathcal{I}_{\text{M}}} \frac{|\mathbf{v}(\mathbf{x}_p)|\triangle t_{\text{M}}}{|\Omega_p|_{B}^{\frac{1}{3}}}, \\
CFL^{\text{max}}_{\text{M}} &= \max_{p \in \mathcal{I}_{\text{M}}}  \frac{|\mathbf{v}(\mathbf{x}_p)|\triangle t_{\text{M}}}{|\Omega_p|_{B}^{\frac{1}{3}}}.
\end{align}
In Table~\ref{tab:CFL}, the corresponding $CFL$s are shown for the mentioned test cases when $\triangle t_{\text{M}}$~\eqref{eq:time_step} is used.

\subsection{Comparison to linearly extended and zero Neumann boundary conditions}\label{sec:comp_other}

\begin{figure}
\begin{center}
\begin{tabular}{ccc}
	\includegraphics[width=3.8cm]{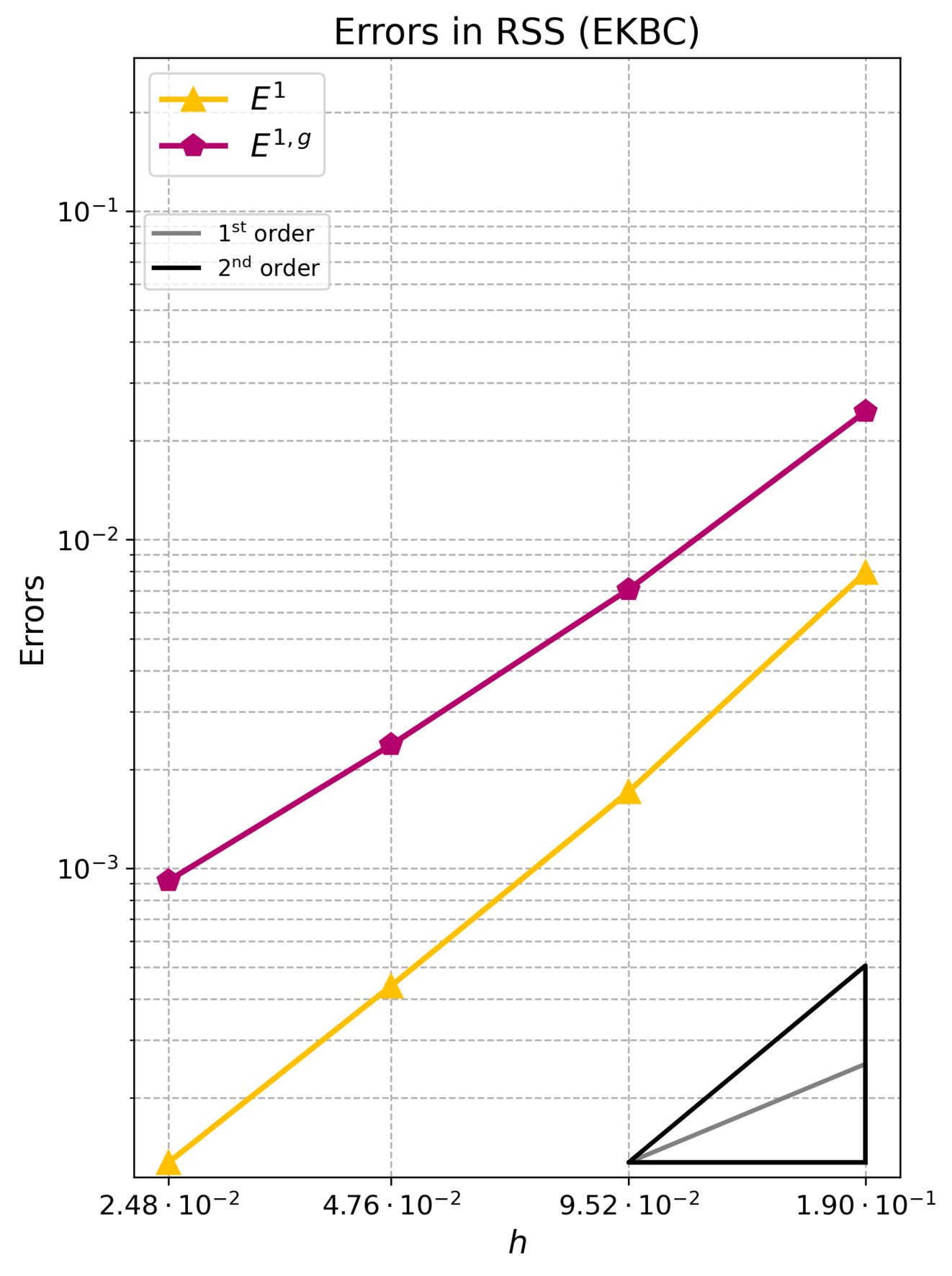} & 
	\includegraphics[width=3.8cm]{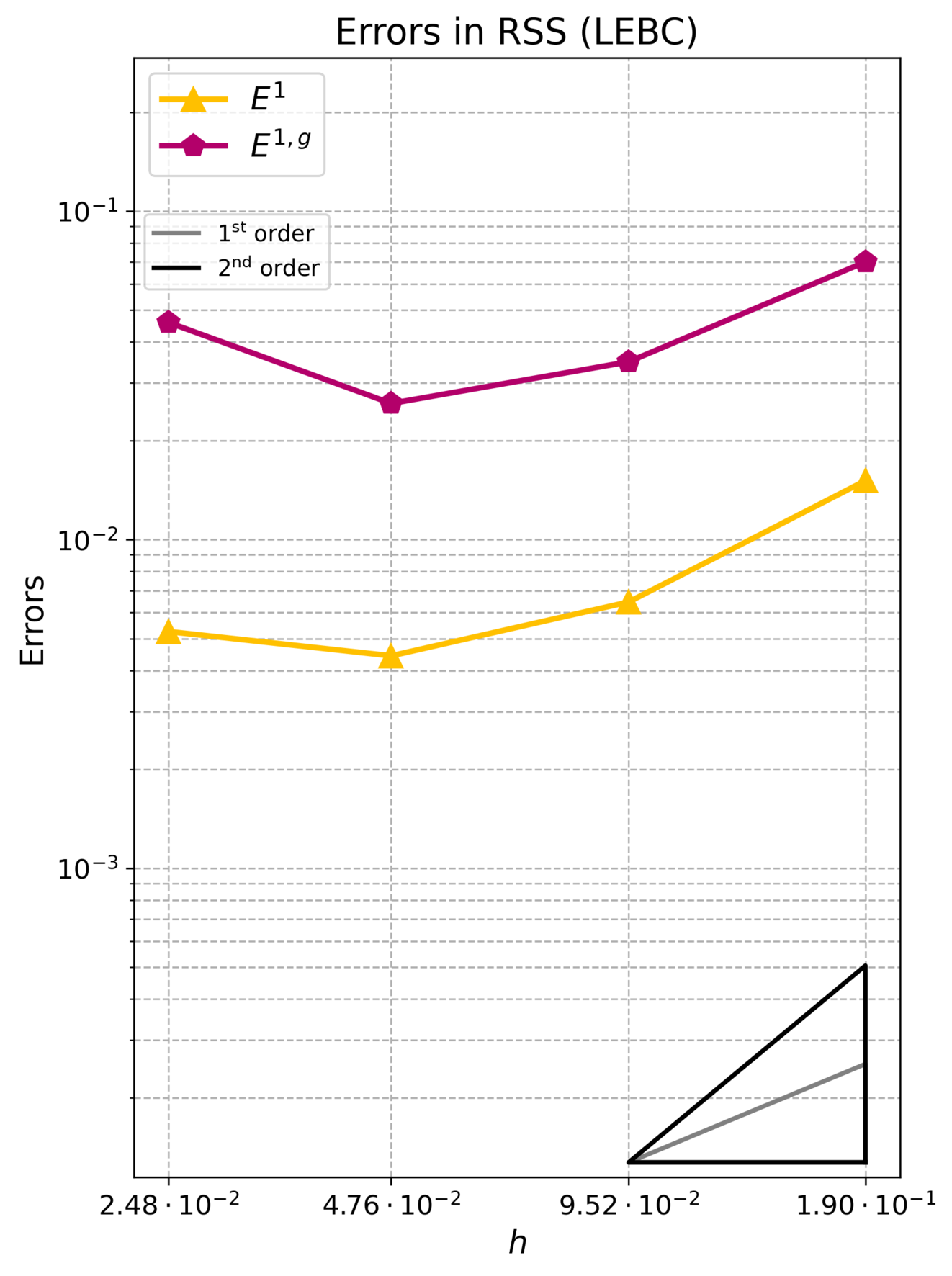} &
	\includegraphics[width=3.8cm]{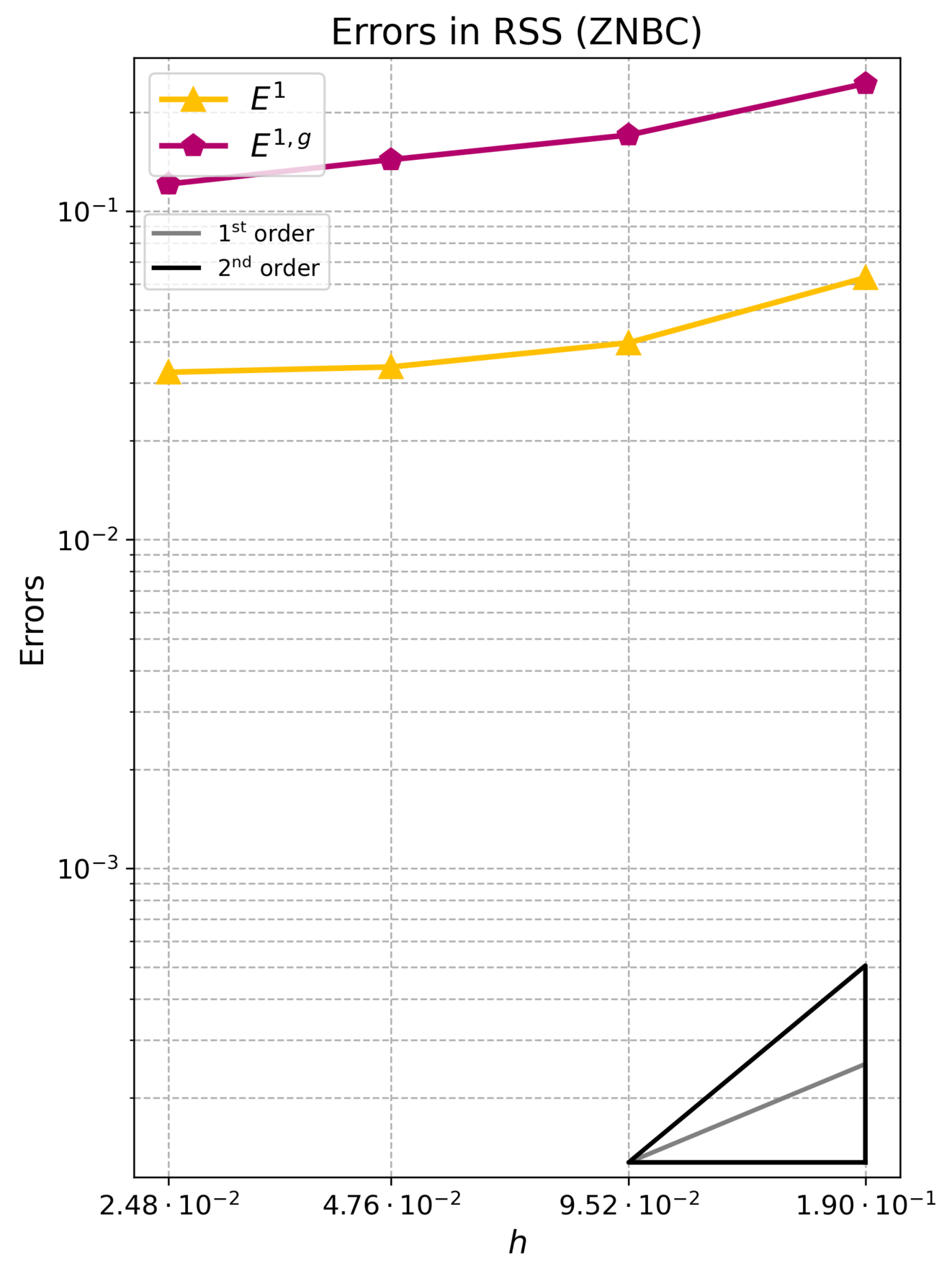} \\
	\includegraphics[width=3.8cm]{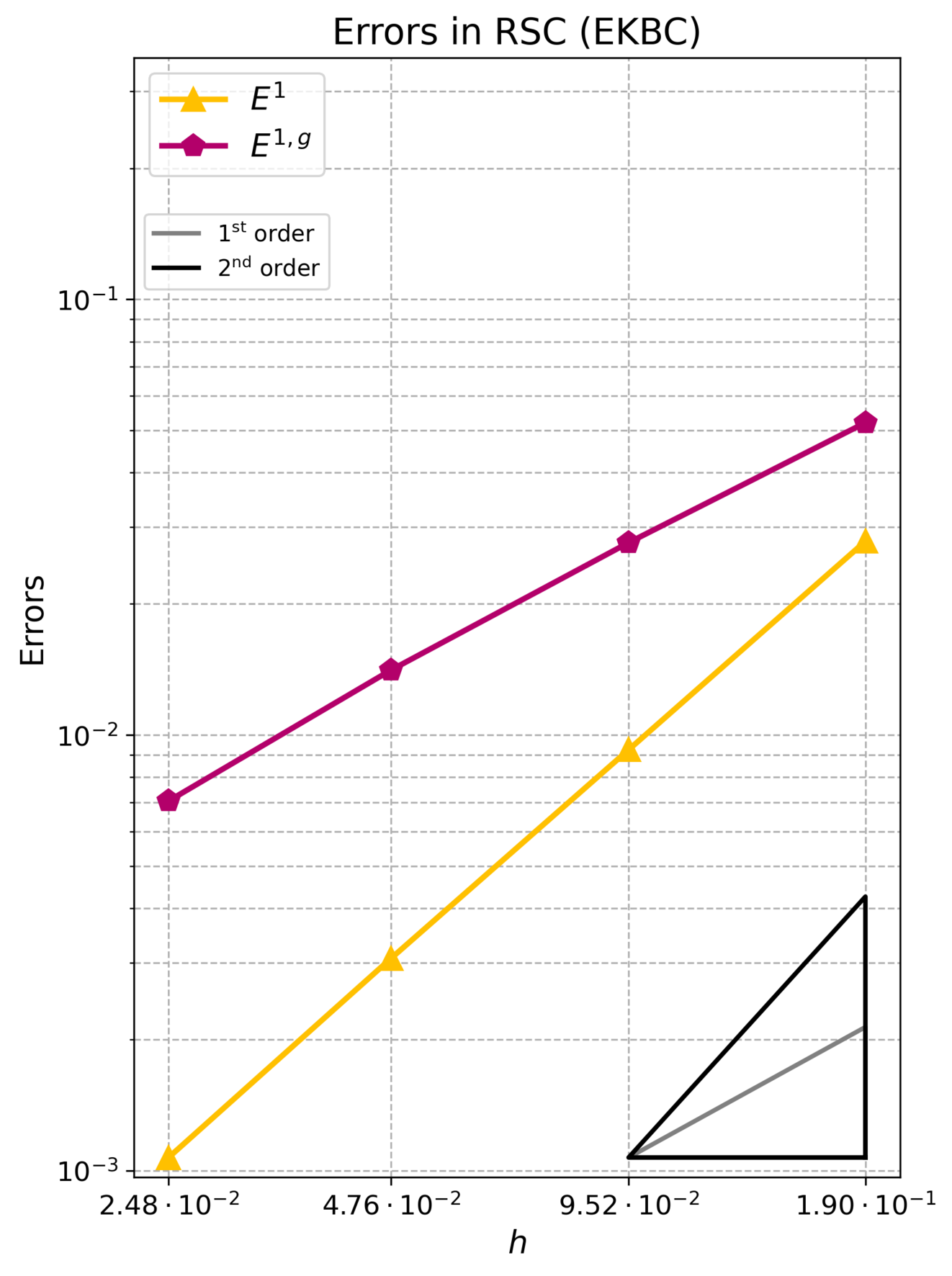} & 
	\includegraphics[width=3.8cm]{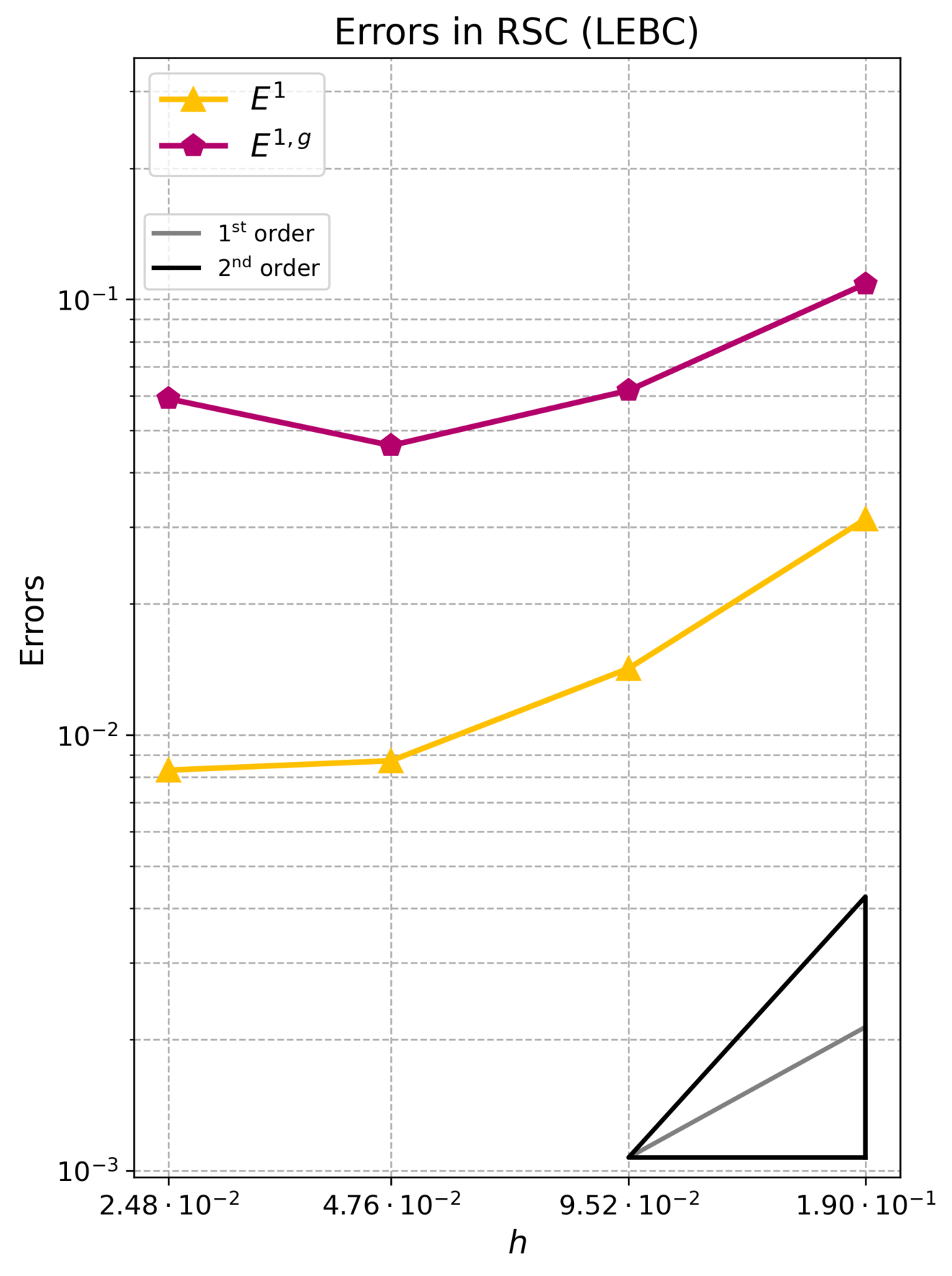} &
	\includegraphics[width=3.8cm]{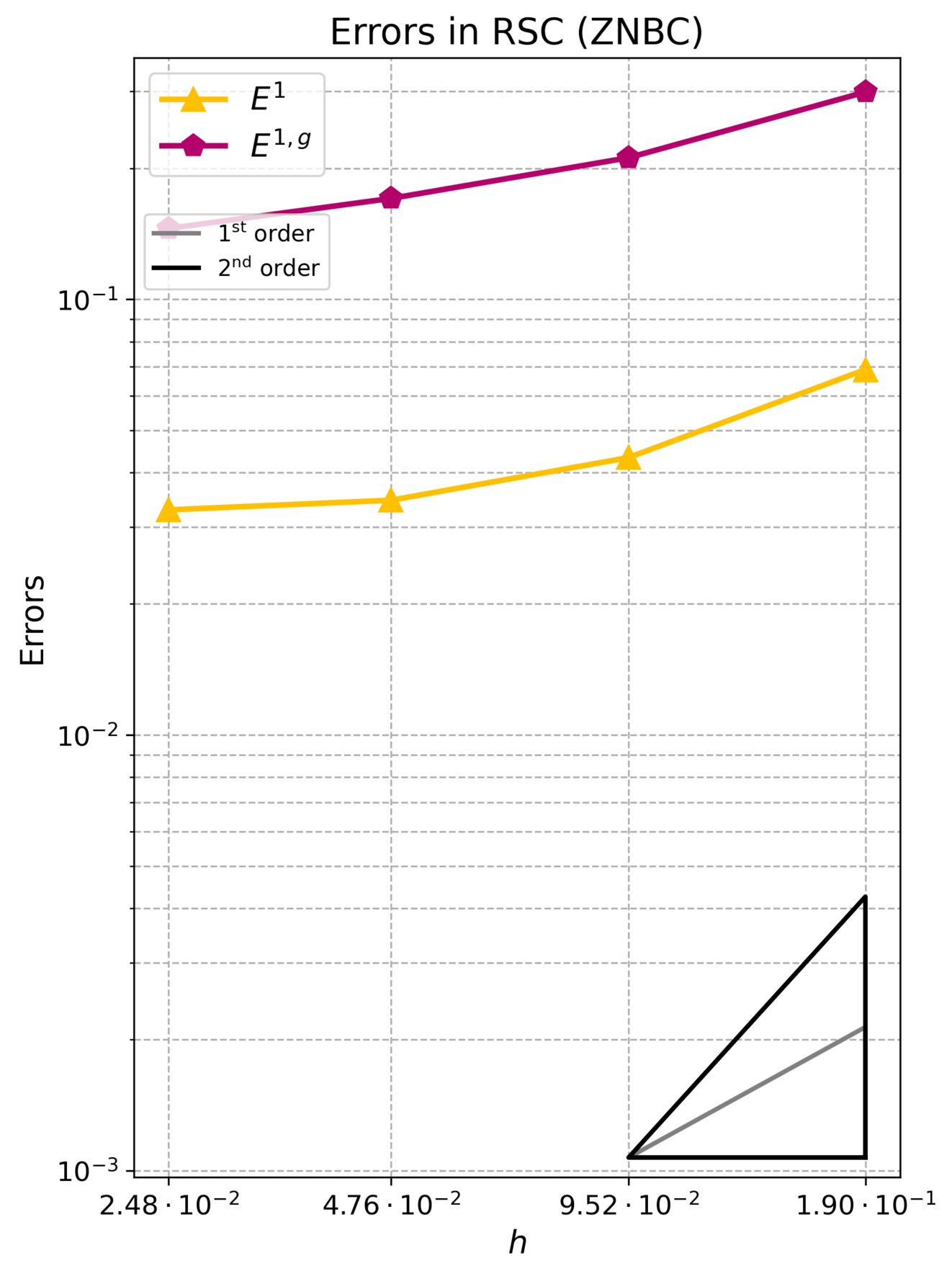} \\
\end{tabular}
\end{center}
\caption{Applying $\triangle t_{\text{M}}$~\eqref{eq:time_step}, the characteristic length versus errors $E^1_{\text{M}}$~\eqref{eq:err_L1} and $E^{1,g}_{\text{M}}$~\eqref{eq:err_L1_g} for test cases~\ref{RSS} and~\ref{RSC} in Table~\ref{tab:test_case} are presented by using EKBC (left column), LEBC (middle column), and ZNBC (right column).} \label{fig:conv_grad_order}	
\end{figure}

In this subsection, we compare the numerical solutions and their gradient results of using EKBC, LEBC, and ZNBC. For two test cases, we use a combination of rotation and shrinking velocity until the end time $T=1$ and the initial surface $\Gamma_{0}(u)$ that is the zero level set of the signed distance function $u^0$. Using equations of sphere and cube~\eqref{eq:sphere_cube}, we describe the analytical solutions of the test cases in Table~\ref{tab:test_case}.

In Figure~\ref{fig:conv_grad_order}, the characteristic length versus errors $E^1$~\eqref{eq:err_L1} and $E^{1,g}$~\eqref{eq:err_L1_g} are presented for the numerical results using EKBC (left column), LEBC (middle column), and ZNBC (right column). The graphs clearly explain that the results of using LEBC or ZNBC do not converge in $\pfb{L^1(\Omega\times(0,T))}$ norm of the error between exact and numerical solutions. More interestingly, the graph with the pentagon symbols shows the convergence behavior of absolute value of the gradient under $\pfb{L^1(\Omega\times(0,T))}$ norm; see~\eqref{eq:err_L1_g}. The numerical convergence of the absolute value of the gradient to be $1$ \nn{confirms that the computation reproduces the property $|\nabla u|=1$ that the exact solution of these test cases has almost everywhere.} \nn{\ The measurement is meaningful here because the exact solution of a rotation combined with a shrinking satisfies $|\nabla u|=1$ almost everywhere. For LEBC and ZNBC the deviation grows from the boundary inwards: the values they supply at the inflow boundary are not compatible with $|\nabla u|=1$, and the inflow carries that into the domain.}

\subsection{General closed surface}\label{sec:gen_surf}

To test an evolving general closed surface, we use the Stanford Bunny from the Stanford 3D scanning repository\footnote{ \texttt{http://graphics.stanford.edu/data/3Dscanrep}.} to rotate or translate the surface. As a Dirichlet boundary condition is not available for such a general surface, numerical results using EKBC are compared with those using ZNBC~\eqref{eq:zeroNeumann} and LEBC~\cite{ref:HMFB17}.

\begin{table}
\makebox[\linewidth][c]{%
\begin{minipage}{1.25\textwidth}
	\centering
	\centering
	\begin{tabular}{
			S[table-format=1]
			S[table-format=7]
			S[table-format=1.2,table-figures-exponent=2,table-sign-mantissa,table-sign-exponent]
			S[table-format=1.2,table-figures-exponent=2,table-sign-mantissa,table-sign-exponent]
			S[table-format=1.2,table-figures-exponent=0,table-sign-mantissa,table-sign-exponent]
			S[table-format=1.2,table-figures-exponent=2,table-sign-mantissa,table-sign-exponent]
			S[table-format=1.2,table-figures-exponent=2,table-sign-mantissa,table-sign-exponent]
			S[table-format=1.2,table-figures-exponent=0,table-sign-mantissa,table-sign-exponent]
			S[table-format=1.2,table-figures-exponent=2,table-sign-mantissa,table-sign-exponent]}
		\toprule
		\multicolumn{3}{c}{$\Sigma$} & \multicolumn{3}{c}{\ref{TB}} &\multicolumn{3}{c}{\ref{RB}}  \\
		\cmidrule(lr){1-3} \cmidrule(lr){4-6} \cmidrule(lr){7-9}  
		\multicolumn{1}{c}{M} & \multicolumn{1}{c}{$\mathcal{J}_{\text{M}}$} & \multicolumn{1}{c}{$h^{\text{ave}}_{\text{M}}$} & \multicolumn{1}{c}{$CFL^{\text{min}}_{\text{M}}$} & \multicolumn{1}{c}{$CFL^{\text{ave}}_{\text{M}}$} & \multicolumn{1}{c}{$CFL^{\text{max}}_{\text{M}}$} & \multicolumn{1}{c}{$CFL^{\text{min}}_{\text{M}}$} & \multicolumn{1}{c}{$CFL^{\text{ave}}_{\text{M}}$} & \multicolumn{1}{c}{$CFL^{\text{max}}_{\text{M}}$} \\
		\midrule
		1 & 737183  & 2.30E-03 & 8.13E-01 & 2.85E+00 & 1.47E+01 & 1.29E-03 & 2.37E+00 & 1.08E+01 \\
		\cmidrule(lr){1-3} \cmidrule(lr){4-6} \cmidrule(lr){7-9} 
		2 & 2391873 & 1.71E-03 & 6.82E-01 & 1.47E+00 & 9.69E+00 & 1.76E-04 & 1.14E+00 & 9.16E+00 \\
		\bottomrule	
	\end{tabular}
\end{minipage}}
\vspace{0.5em}
\begin{center}
\caption{Two levels of discretizing the computational domain, $\Sigma = [-0.0682,0.1078]\times[-0.092,0.094]\times[-0.1483,0.0497] \subset \mathbb{R}^3$, used in~\ref{TB} are presented. We also show the corresponding characteristic length~\eqref{eq:char_len} and $CFL$s in~\ref{TB} and~\ref{RB}.\label{tab:dom_bunny}}
\end{center}
\end{table}

\begin{figure}
\begin{center}
\begin{tabular}{cc}
	\includegraphics[height=5cm]{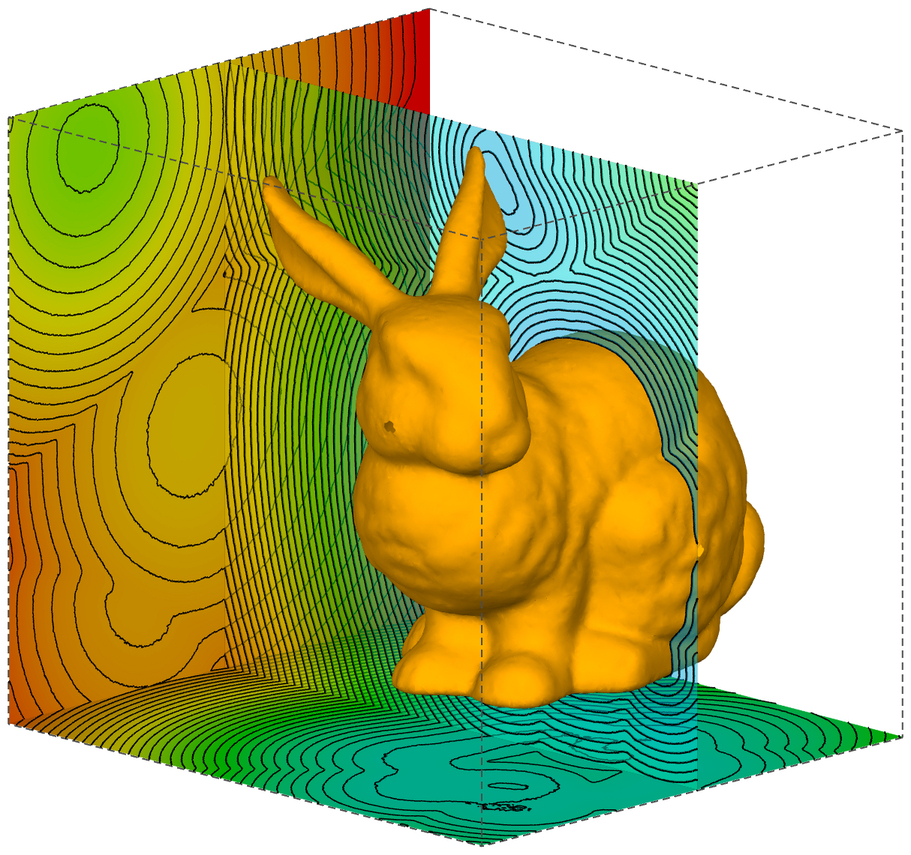} &
	\includegraphics[height=5cm]{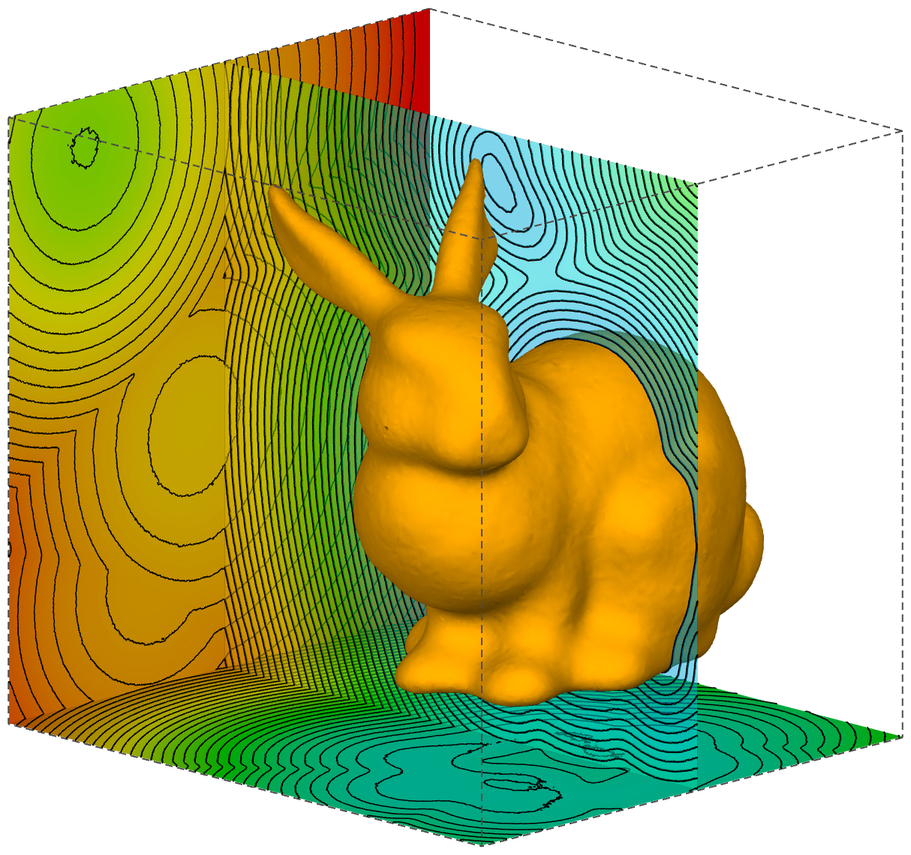} \\
	(a) $T=0$ & (b) $T=2$ \\
	\includegraphics[height=5cm]{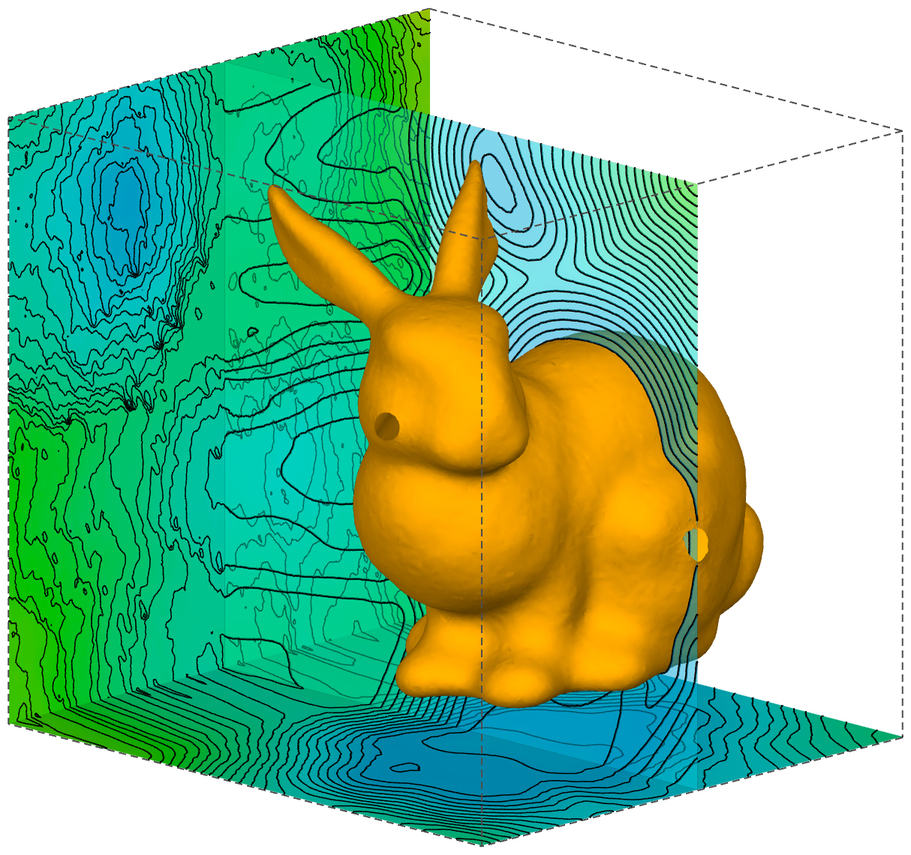} &
	\includegraphics[height=5cm]{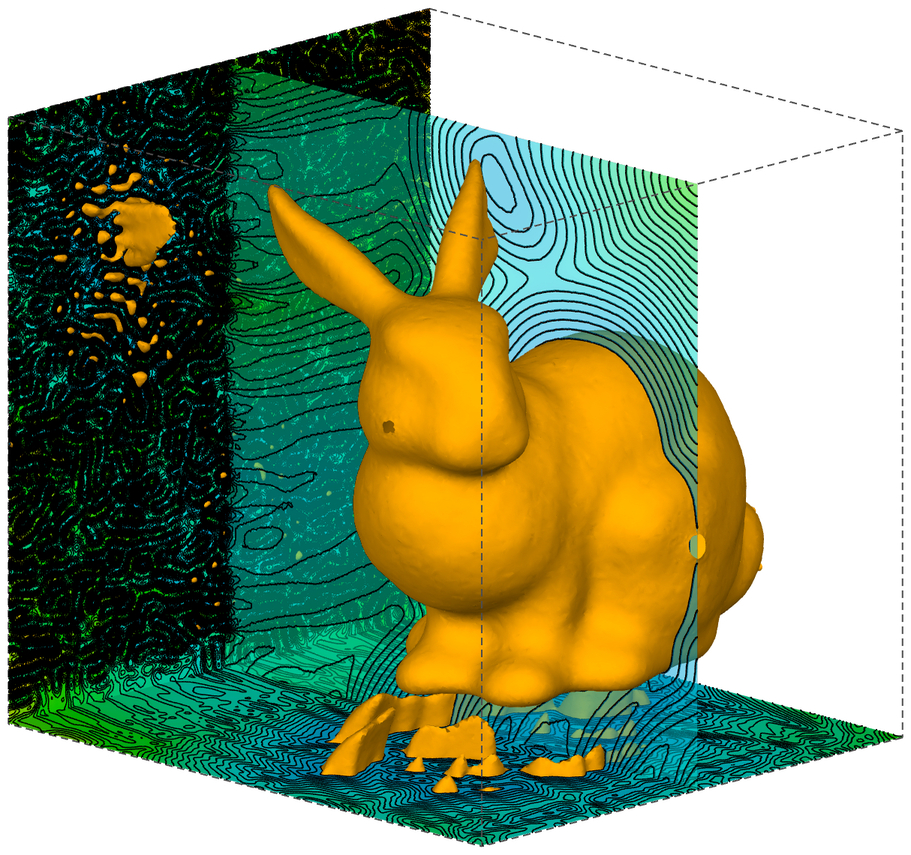} \\
	(c) $T=2$ & (d) $T=2$ \\
\end{tabular}
\end{center}
\caption{(a) The initial surface of~\ref{TB} is presented with the $60$ levels of isocurves of $u^0$ on a middle plane and two sides of boundary of the computational domain. The results of~\ref{TB} are presented at $T=2$ using the EKBC, ZNBC, and LEBC in (b), (c), and (d), respectively.} \label{fig:TranBunny}
\end{figure}

\begin{figure}
\begin{center}
\begin{tabular}{cc}
	\includegraphics[height=5cm]{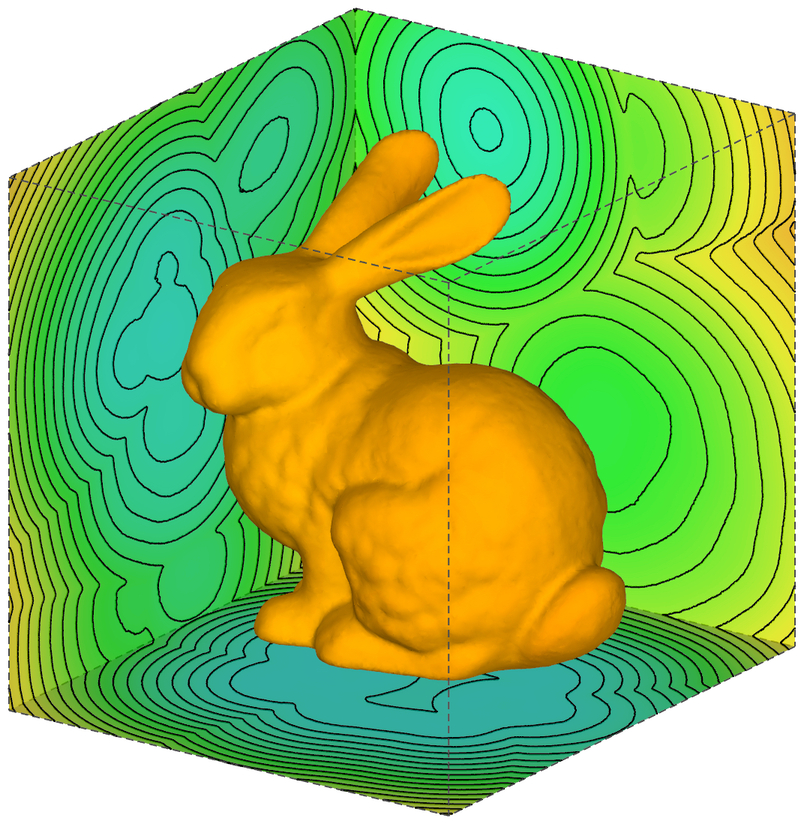} &
	\includegraphics[height=5cm]{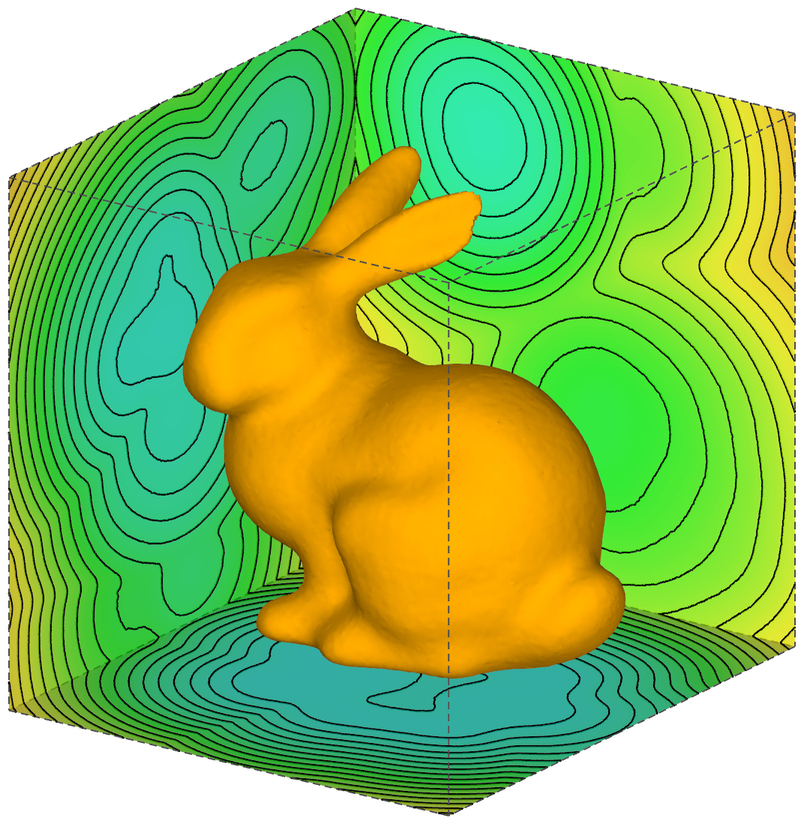} \\
	(a) $T=0$ & (b) $T=2$ \\
	\includegraphics[height=5cm]{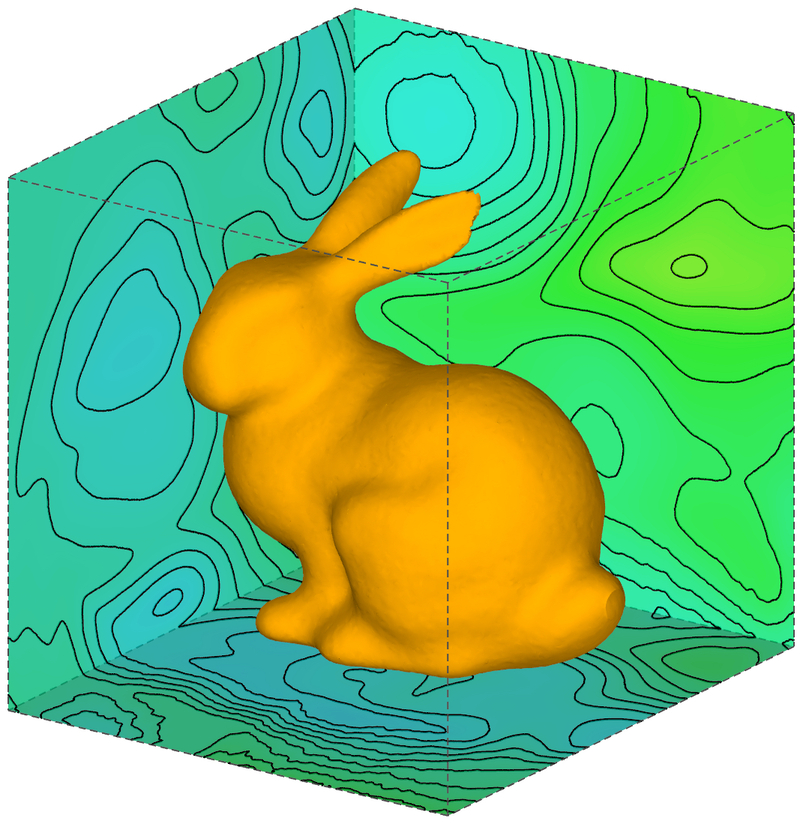} &
	\includegraphics[height=5cm]{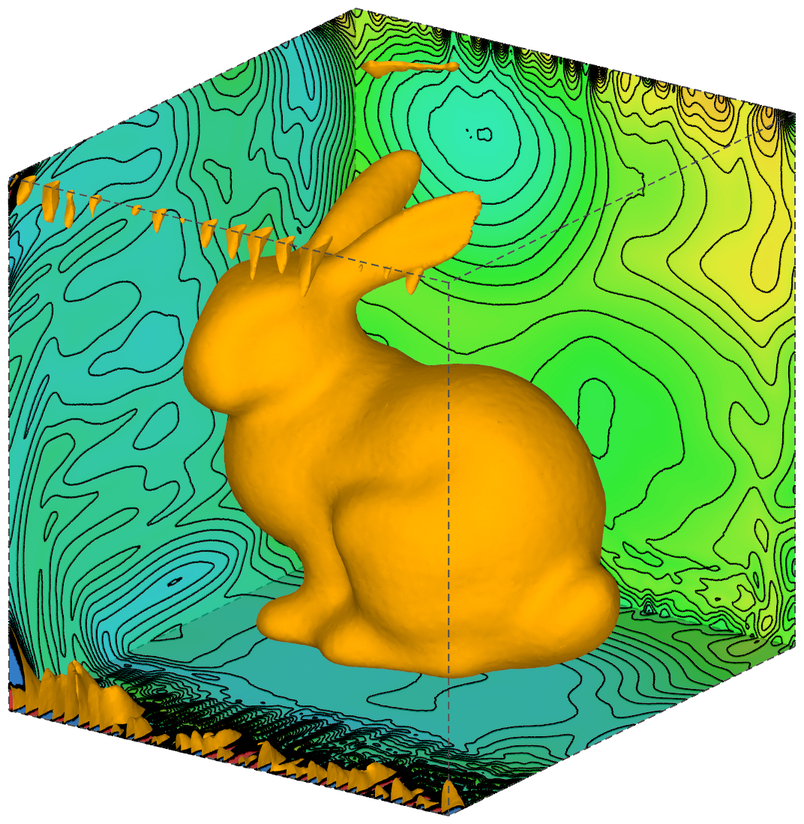} \\
	(c) $T=2$ & (d) $T=2$
\end{tabular}
\end{center}
\caption{(a) The initial surface of~\ref{RB} is presented with the $40$ levels of isocurves of $u^0$ on three sides of boundary of the computational domain.  The results of~\ref{RB} are presented at $T=2$ using the EKBC, ZNBC, and LEBC in (b), (c), and (d), respectively.} \label{fig:RotBunny}
\end{figure}

The initial signed distance function~\eqref{eq:initcond} is obtained by Laplacian regularized eikonal equation~\cite{ref:HMF24}. In Figures~\ref{fig:TranBunny}-(a) and~\ref{fig:RotBunny}-(a), we present the initial surface (Stanford Bunny) and the isocurves on some planes in the computational domain. The details of the test cases are explained below:
\begin{itemize}
\item \textbf{TB}: \namelabel{TB} A translation of the Stanford Bunny in Figure~\ref{fig:TranBunny}-(a) until $T=2$ with the velocity:
\begin{align}\label{eq:vel_tran_bunny}
\mathbf{v} =
\begin{cases}
	(0.0198, -0.031, -0.0768), & 0.0 \leq T \leq 1.0,\\
	(-0.0198, 0.031, 0.0768), &  1.0 < T \leq 2.0.
\end{cases}
\end{align}
The computational domain $\Sigma$ in Table~\ref{tab:dom_bunny} is deliberately selected so that the boundary of the domain is in close proximity to the three sides of the initial surface. We use the time step $\pfb{\triangle} \tau_{\text{M}} = 5.0\cdot10^{-2} \cdot 2^{-(M-1)}$ for $M \in \{1,2\}$ and the final time $T=2$.
\item \textbf{RB}: \namelabel{RB} A rotation of the Stanford Bunny in Figure~\ref{fig:RotBunny}-(a) until $T=2$ with the velocity:
\begin{align}\label{eq:vel_rot_bunny}
\mathbf{v} =
\begin{cases}
	(\pi x_3, 0, -\pi x_1 ), & 0.0 \leq T \leq 0.5,\\
	(0, -\pi x_3, \pi x_2 ), & 0.5 < T \leq 1.0,\\
	(\pi x_2, -\pi x_1, 0), &  1.0 < T \leq 1.5, \\
	(0, \pi x_3, -\pi x_2), & 1.5 < T \leq 2.0.
\end{cases}
\end{align}
In this test case, the computational domain $\Sigma$ in Table~\ref{tab:dom_bunny} is shifted by a vector $(-0.0198,0.016,0.0503)$ to include rotated bunny in the domain. We use the time step $ \pfb{\triangle} \tau_{\text{M}} = 1.25\cdot10^{-2} \cdot 2^{-(M-1)}$ for $M \in \{1,2\}$ and the final time $T=2$.
\end{itemize}
The final time and the velocity functions~\eqref{eq:vel_tran_bunny} and~\eqref{eq:vel_rot_bunny} are selected to ensure that the closed surface remains within the computational domain. Moreover, the velocity fields chosen in the above examples return the evolving surface to the location of the initial surface when the final time is reached.

In Figure~\ref{fig:TranBunny}, the results at $T=2$ of~\ref{TB} using the EKBC, ZNBC, and LEBC are presented in (b), (c), and (d), respectively. Clearly, comparing the result of the proposed algorithm with the initial isocurves in Figure~\ref{fig:TranBunny}-(a), the isocurves are well preserved in the case of constant advective velocity~\eqref{eq:vel_tran_bunny}. However, the results of solving~\eqref{eq:levelset} with the ZNBC or LEBC show a severe distortion near the inflow boundaries. After a long period of evolving the surface, the error in the domain accumulates more and more and causes significant inaccuracy or instability on the evolving surface. In Figure~\ref{fig:RotBunny}, the same phenomenon is observed in the case of using the rotational velocity. Since the velocity function~\eqref{eq:vel_rot_bunny} returns the evolving surface to the initial surface when the final time reaches, the results at $T=2$ should be similar to the initial surface and the isocurves in Figure~\ref{fig:RotBunny}-(a). The result of the proposed method shows to achieve the mentioned objective, but the results using the ZNBC or LEBC cannot obtain the same accuracy.

\begin{table}
\centering
\renewcommand{\arraystretch}{1.1}
\begin{tabular}{cS[table-format=1]
	S[table-format=1.2,table-figures-exponent=2,table-sign-mantissa,table-sign-exponent]
	S[table-format=1.4]
	S[table-format=1.2,table-figures-exponent=2,table-sign-mantissa,table-sign-exponent]
	S[table-format=1.4]
	S[table-format=1.2,table-figures-exponent=2,table-sign-mantissa,table-sign-exponent]
	S[table-format=1.4]
	S[table-format=1.2,table-figures-exponent=2,table-sign-mantissa,table-sign-exponent]
	S[table-format=1.4]
	S[table-format=1.2,table-figures-exponent=2,table-sign-mantissa,table-sign-exponent]
	S[table-format=1.4]
	S[table-format=1.2,table-figures-exponent=2,table-sign-mantissa,table-sign-exponent]
	S[table-format=1.4]}
\toprule
& & \multicolumn{2}{c}{EKBC} & \multicolumn{2}{c}{LEBC} & \multicolumn{2}{c}{ZNBC} \\
\cmidrule(lr){1-2} \cmidrule(lr){3-4} \cmidrule(lr){5-6} \cmidrule(lr){7-8}
& \multicolumn{1}{c}{M} & \multicolumn{1}{c}{$e^1_{\text{M}}$} & \multicolumn{1}{c}{$ROE_{e^1_{\text{M}}}$} & \multicolumn{1}{c}{$e^1_{\text{M}}$} &  \multicolumn{1}{c}{$ROE_{e^1_{\text{M}}}$} & \multicolumn{1}{c}{$e^1_{\text{M}}$} & \multicolumn{1}{c}{$ROE_{e^1_{\text{M}}}$}
\\
\cmidrule(lr){1-2} \cmidrule(lr){3-4} \cmidrule(lr){5-6} \cmidrule(lr){7-8}
\multirow{2}{*}{\ref{TB}} & 1 & 1.05E-03 & 5.1957 & 1.52E-02 & 0.1196   & 1.50E-02 & -0.0021 \\
& 2 & 2.27E-04 &        & 1.47E-02 &          & 1.50E-02 &         \\
\cmidrule(lr){1-2} \cmidrule(lr){3-4} \cmidrule(lr){5-6} \cmidrule(lr){7-8}
\multirow{2}{*}{\ref{RB}} & 1 & 1.00E-03 & 4.3227 & 7.60E-03 & 2.5613   & 6.31E-03 & -0.1204 \\
& 2 & 2.80E-04 &        & 3.57E-03 &          & 6.53E-03 &         \\
\midrule
& \multicolumn{1}{c}{M} & \multicolumn{1}{c}{$e^\infty_{\text{M}}$} & \multicolumn{1}{c}{$ROE_{e^\infty_{\text{M}}}$} & \multicolumn{1}{c}{$e^\infty_{\text{M}}$} &  \multicolumn{1}{c}{$ROE_{e^\infty_{\text{M}}}$} & \multicolumn{1}{c}{$e^\infty_{\text{M}}$} & \multicolumn{1}{c}{$ROE_{e^\infty_{\text{M}}}$}
\\
\cmidrule(lr){1-2} \cmidrule(lr){3-4} \cmidrule(lr){5-6} \cmidrule(lr){7-8}
\multirow{2}{*}{\ref{TB}} & 1 & 2.15E-02 & 5.2533 & 8.51E-02 & -1.4906  & 8.51E-02 & -0.2715 \\
&  2 & 4.57E-03 &        & 1.32E-01 &          & 9.21E-02 &         \\
\cmidrule(lr){1-2} \cmidrule(lr){3-4} \cmidrule(lr){5-6} \cmidrule(lr){7-8}
\multirow{2}{*}{\ref{RB}} & 1 & 2.15E-02 & 4.8894 & 6.48E-02 & -13.0739 & 7.75E-02 & 0.2955  \\
& 2 & 5.09E-03 &        & 3.05E+00 &          & 7.11E-02 &   \\
\midrule
& \multicolumn{1}{c}{M} & \multicolumn{1}{c}{$e^{\text{v}}_{\text{M}}$} & \multicolumn{1}{c}{$ROE_{e^{\text{v}}_{\text{M}}}$} & \multicolumn{1}{c}{$e^{\text{v}}_{\text{M}}$} &  \multicolumn{1}{c}{$ROE_{e^{\text{v}}_{\text{M}}}$} & \multicolumn{1}{c}{$e^{\text{v}}_{\text{M}}$} & \multicolumn{1}{c}{$ROE_{e^{\text{v}}_{\text{M}}}$}
\\
\cmidrule(lr){1-2} \cmidrule(lr){3-4} \cmidrule(lr){5-6} \cmidrule(lr){7-8}
\multirow{2}{*}{\ref{TB}} & 1 & 8.77E-06 & 5.6631 & \multicolumn{1}{c}{-} & \multicolumn{1}{c}{-} & 5.04E-06 & 4.5583 \\
& 2 & 1.65E-06 &  & \multicolumn{1}{c}{-} & & 1.32E-06 & \\
\cmidrule(lr){1-2} \cmidrule(lr){3-4} \cmidrule(lr){5-6} \cmidrule(lr){7-8}
\multirow{2}{*}{\ref{RB}} & 1 & 9.69E-06 & 4.0944 & \multicolumn{1}{c}{-} & \multicolumn{1}{c}{-} & 1.01E-05 & 4.0571 \\
& 2 & 2.90E-06 & & \multicolumn{1}{c}{-} & & 3.05E-06 &  \\
\midrule
& \multicolumn{1}{c}{M} & \multicolumn{1}{c}{$E^{1,g}_{\text{M}}$} & \multicolumn{1}{c}{$EOC_{E^{1,g}_{\text{M}}}$} & \multicolumn{1}{c}{$E^{1,g}_{\text{M}}$} &  \multicolumn{1}{c}{$EOC_{E^{1,g}_{\text{M}}}$} & \multicolumn{1}{c}{$E^{1,g}_{\text{M}}$} & \multicolumn{1}{c}{$EOC_{E^{1,g}_{\text{M}}}$}
\\
\cmidrule(lr){1-2} \cmidrule(lr){3-4} \cmidrule(lr){5-6} \cmidrule(lr){7-8}
\multirow{2}{*}{\ref{TB}}  & 1 & 4.91E-02 & 2.8324 & 2.31E-01 & 0.8104 & 2.09E-01 & 0.3273  \\
& 2 & 2.13E-02 &        & 1.82E-01 &        & 1.90E-01 &         \\
\cmidrule(lr){1-2} \cmidrule(lr){3-4} \cmidrule(lr){5-6} \cmidrule(lr){7-8}
\multirow{2}{*}{\ref{RB}}  & 1 & 4.71E-02 & 2.4147 & 1.98E-01 & 0.3496 & 1.66E-01 & -0.5766 \\
& 2 & 2.31E-02 &        & 1.79E-01 &        & 1.97E-01 &  \\
\bottomrule
\end{tabular}
\vspace{0.5em}
\caption{The reduction of error ($ROE$) is computed by the same method in $EOC$~\eqref{eq:EOC} with the errors~\eqref{eq:err_1_bunny},~\eqref{eq:err_inf_bunny}, and~\eqref{eq:err_v_bunny} for~\ref{TB} and~\ref{RB}. The errors $e^{\text{v}}_{\text{M}}$ in case of LEBC are omitted due to fluctuating volumes caused by the creation and disappearance of surfaces over time; see the text for details.\label{tab:err_bunny}}
\end{table}

Since the exact solution of test cases~\ref{TB} and~\ref{RB} is unknown, it is not feasible to compute the same errors introduced at the beginning of Section~\ref{sec:NumEx}. Instead, we access the following quantitative errors below at the final time using the initial condition $u^0(\mathbf{x})$, as the given velocity returns the evolving surface to their original shape:
\begin{align}
e^1_{\text{M}} &=  \frac{1}{\sum_{p \in \mathcal{J}_\text{M}} |\Omega_p| } \sum_{p \in \mathcal{J}_\text{M}} |u(\mathbf{x}_p,T) - u^0_{\text{M}}(\mathbf{x}_p)| |\Omega_p|,  \label{eq:err_1_bunny} \\ 
e^\infty_{\text{M}} &=  \max_{p \in \mathcal{J}_\text{M}} |u(\mathbf{x}_p,T) - u^0_{\text{M}}(\mathbf{x}_p)|. \label{eq:err_inf_bunny} 
\end{align}
For $t \in [0,T]$, the average of volume difference between evolving surface $V(\Gamma_{t}(u))$ and the initial surface $V(\Gamma_0(u^{0}))$ can be measured: 
\begin{align}
e^{\text{v}}_{\text{M}} &= \frac{1}{N} \sum_{n=1}^{N} |V(\Gamma_{t^{n}}(u)) - V(\Gamma_0(u^{0}_{\text{M}}))|, \label{eq:err_v_bunny}
\end{align}
When using LEBC, the volume difference is not meaningful, as surfaces closed to the boundary appears and disappear over time, causing the overall volume $V(\Gamma_{t}(u))$ to fluctuate. While the computed volume difference with respect to the initial volume may occasionally reach zero, this does not imply the exactness of keeping the volume; see Figures~\ref{fig:TranBunny}-(d) and~\ref{fig:RotBunny}-(d). The deviation of $|\nabla u|$ from $1$ is measured by~\eqref{eq:err_L1_g} and the gradient $\nabla u$ \pfb{defined in}~\eqref{eq:grad} is the computed value of the numerical results from test cases in this subsection.

In Table~\ref{tab:err_bunny}, we present the errors defined by~\eqref{eq:err_1_bunny},~\eqref{eq:err_inf_bunny}, and~\eqref{eq:err_v_bunny} for both~\ref{TB} and~\ref{RB}. The reduction of error ($ROE$) is computed similarly to $EOC$~\eqref{eq:EOC}, but with the corresponding errors. In this context, the exact solution used in $EOC$~\eqref{eq:EOC} is replaced by the function $u^0_\text{M}$, which is mesh-dependent. As a result, $ROE$ is not directly comparable to $EOC$, but it serves as an indicator of whether the observed errors tend to decrease or increase. An interesting observation is the significantly smaller errors achieved by EKBC compared to LEBC and ZNBC in both~\ref{TB} and~\ref{RB}. This suggests that EKBC maintains the accuracy of the evolution of level set function, even in complex shapes, while other methods suffer from distortions, especially near boundaries. Additionally, \nn{the deviation of $|\nabla u|$ from $1$ is markedly smaller with EKBC. Since~\ref{TB} and~\ref{RB} are rigid motions, the evolution itself preserves $|\nabla u|=1$, and EKBC supplies exactly that constraint at the inflow boundary, whereas LEBC and ZNBC supply values incompatible with it which the inflow carries into the domain, leading to inaccuracies, particularly near the boundaries;} see Figures~\ref{fig:TranBunny} and~\ref{fig:RotBunny}.

For some applications where the signed distance information is necessary in a whole domain, the typical choice of boundary conditions such as ZNBC or LEBC certainly needs an extra procedure to keep the distance property from the evolving surface, for example, the reinitialization. However, it is noteworthy that while reinitialization successfully maintains the desired distance property, it inevitably introduces unintended perturbations that manifest as undesired movement of the evolving surface.


\section{Conclusion}\label{sec:conclusion}
We propose to use the eikonal boundary condition when the evolution of the surface in the advective or normal flow equations is solved. The numerical results confirm that the experimental order convergence and the errors are comparable to the use of the Dirichlet boundary condition, which is mostly impossible to impose except in the cases of known exact analytical solutions. Moreover, the eikonal boundary condition brings the possibility of using a large $CFL$ number and makes it more practical to use the level set method in industrial applications. The numerical results of the general closed surface\nn{, for which no exact solution and hence no Dirichlet boundary condition is available,} show that the proposed \nn{boundary condition remains robust where the linearly extended and the zero Neumann conditions distort the surface.} One of the future works is to combine the proposed eikonal boundary condition with the extended velocity field~\vv{\cite{ref:AS99}} on polyhedral meshes to keep the signed distance property in addition while solving the level set equation with general velocity fields.

%

\section*{Acknowledgments}
The authors sincerely thank Dr. Branislav Basara and Dr. Reinhard Tatschl in AVL List GmbH, Austria, for supporting the University Partnership Program\footnote{See more details in AVL Advanced Simulation Technologies University Partnership Program: \href{https://www.avl.com/documents/10138/3372587/AVL_UPP_Flyer.pdf}{https://www.avl.com/documents/10138/3372587/AVL\_UPP\_Flyer.pdf}}.

\section*{Declarations}
Not applicable
\begin{itemize}
\item Conflict of interest/Competing interests: No conflicts, no competing interests
\item Availability of data and materials: The data that support the findings of this study are available on request from the corresponding author.
\item Authors' contributions: Equal
\end{itemize}

\begin{appendices}

\section{Technical Details: Computation of gradients\label{sec:app_A}}

\begin{figure}
\begin{center}
	\begin{tabular}{c}
		\includegraphics[height=4cm]{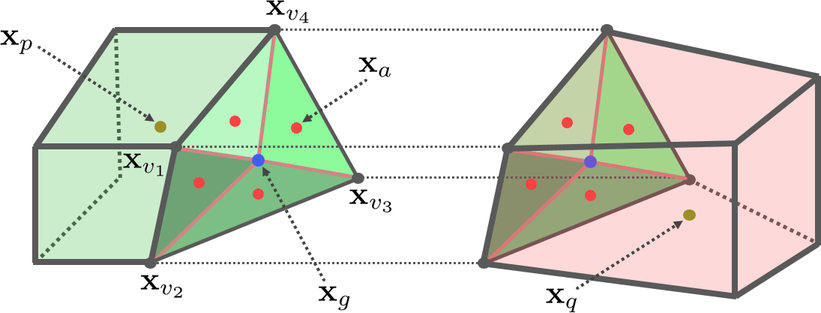}
	\end{tabular}
\end{center}
\caption{Two cells, $\Omega_p$ and $\Omega_q$ in Figure~\ref{fig:cell_dig}-(d) are redrawn separately to the left and right to describe more detail; centers of cells ($\mathbf{x}_p$ and $\mathbf{x}_q$), the center of the face ($\mathbf{x}_g$), the center of the triangle ($\mathbf{x}_a$), and the ordered vertices of the face ($V_g = \{\mathbf{x}_{v_j}: j = 1, 2, \ldots, 4\}$). All definitions are explained in the end of the Appendix.} \label{fig:tess_full}
\end{figure}

We explain technical details on how to compute the representative gradient at the center $\mathbf{x}_p$ of the polyhedron cell~$\Omega_p$,~$p \in \mathcal{I}$ and at the center $\mathbf{x}_a$ of the triangle~$e_a$,~$a \in \mathcal{F} \cup \mathcal{B}$. In~Figure~\ref{fig:tess_full}, the diagram in Figure~\ref{fig:cell_dig}-(d) is redrawn to illustrate necessary points in the rest of the Appendix.

In order to define the mentioned gradients, we begin with the concept of the gradient computed by the least-squares method. In the finite volume method, we use the function~$u$ whose value is constant on the cell~$\Omega_p$,~$p \in \mathcal{I}$. Using the least-squares method, a gradient at the center~$\mathbf{x}_p$ of the cell~$\Omega_p$ is computed:
\begin{equation}\label{eq:grad}
\nabla u_{\text{S}_p}(\mathbf{x}_p) \equiv \argmin_{\mathbf{y}} \sum_{\mathbf{x} \in \text{S}_p} \left( \omega_p(\mathbf{x}) \big( \mathbf{y} \cdot (\mathbf{x} - \mathbf{x}_p) - (u(\mathbf{x}) - u(\mathbf{x}_p))\big) \right)^2,
\end{equation}
where a weight function is defined by
\begin{align}
\omega_p(\mathbf{x})=\frac{1}{|\mathbf{x} - \mathbf{x}_p|}
\end{align}
and the set $\text{S}_p$ specifies points used to calculate the gradient; see the explicit formula to compute~$\nabla u_{\text{S}_p}(\mathbf{x}_p)$ in~\cite{ref:HMFB17}. When using the Dirichlet boundary condition, we choose 
\begin{equation}
\text{S}_p = \text{D}_p \equiv
\begin{cases} \{ \mathbf{x}_q \: | \: q \in \mathcal{N}_p \} \quad &\text{if} \:\: \mathcal{B}_p = \emptyset, \\
	\{ \mathbf{x}_q \: | \: q \in \mathcal{N}_p \} \cup \{ \mathbf{x}_b \: | \: b \in \mathcal{B}_p \} \quad &\text{if} \:\: \mathcal{B}_p \neq \emptyset.
\end{cases}
\end{equation}
When the Dirichlet boundary condition is not available, an alternative choice of values next to the cell $\Omega_p$ is
\begin{equation}
\text{S}_p = \text{L}_p \equiv \{ \mathbf{x}_q \: | \: q \in \mathcal{N}_p \}, \quad p \in \mathcal{I}.
\end{equation}
A geometrical interpretation of using $\text{L}_p$ is presented in~\cite{ref:HMFB17}. Note that $\nabla u_{\text{D}_p}(\mathbf{x}_p) \neq \nabla u_{\text{L}_p}(\mathbf{x}_p)$, for $p \in \mathcal{I}_{\text{bdr}}$ and $\nabla u_{\text{D}_p}(\mathbf{x}_p) = \nabla u_{\text{L}_p}(\mathbf{x}_p)$, for $p \in \mathcal{I}_{\text{int}}$. 

From the gradient~\eqref{eq:grad}, we prepare the necessary values to define the representative gradient at the center of the triangle~$e_a$,~$a \in \mathcal{F} \cup \mathcal{B}$ in the following steps:
\begin{enumerate}	
\renewcommand{\labelenumi}{\textbf{\theenumi}}
\renewcommand{\theenumi}{S\arabic{enumi}}
\item \label{A:step1} \textit{The values at the vertex of the cell}: We compute the values of $u$ at the vertex $\mathbf{x}_v$ on the cell by using the inverse distance weighted average:
\begin{align}\label{eq:ver_val}
	u(\mathbf{x}_v) = \frac{\displaystyle \sum_{p \in \mathcal{N}_v} \omega_p(\mathbf{x}_v) \left( u(\mathbf{x}_p) + \nabla u_{\text{S}_p}(\mathbf{x}_p) \cdot (\mathbf{x}_v - \mathbf{x}_p) \right)}{\displaystyle \sum_{p \in \mathcal{N}_v} \omega_p(\mathbf{x}_v)},
\end{align}
where $\mathcal{N}_v = \{p \in \mathcal{I} : \mathbf{x}_v \in \partial \Omega_p \}$ is the set of all indices of cells that have the vertex $\mathbf{x}_v$. 
\item \label{A:step2} \textit{The tessellation of the face on a cell}: From an internal face $e_g = \partial \Omega_p \cap \partial \Omega_q$ or a boundary face $e_g = \partial\Omega_p \cap \partial \Omega$, denoting the center of the face in Figure~\ref{fig:tess_full} as $\mathbf{x}_g$ and the ordered vertices of the face as
\begin{align*}
	\text{V}_g = \{\mathbf{x}_{v_j}: j = 1, 2, \ldots, J\}
\end{align*}
and using the cyclic notation $\mathbf{x}_{v_{J+1}} = \mathbf{x}_{v_1}$, we define a triangle by three vertices:
\begin{align*}
	e_a = \triangle_j = \triangle(\mathbf{x}_{v_j}, \mathbf{x}_{v_{j+1}}, \mathbf{x}_g).
\end{align*}
In this paper, we refer the face $e_g$ as the collection of triangles $\triangle_j$, $j=1, 2, \ldots, J$, for example, see four triangles in Figure~\ref{fig:tess_full}. A triangle $e_a$, $a\in \mathcal{F} \cup \mathcal{B}$, is always one of tessellation of a face of the cell.
\item \label{A:step3} \textit{The values at the center of internal and boundary faces}: Similarly to~\eqref{eq:grad}, we compute the value of $u$ at the center of the face using a generalized diamond-cell strategy~\cite{ref:CP95,ref:CVV99,ref:DM08,ref:MR09,ref:MO10}:
\begin{align*}
	(a_g^{*}, \mathbf{b}_g^{*}) = \argmin_{(a, \mathbf{b}) \in \mathbb{R} \times \mathbb{R}^3 }\sum_{\mathbf{x} \in \text{P}_g} \left( w_g(\mathbf{x})  \big(a + \mathbf{b} \cdot (\mathbf{x} -\mathbf{x}_g) - u(\mathbf{x}) \big) \right)^2,
\end{align*}
where $\text{P}_g$ is the set of vertices of the polygonal pyramid consisting of the bottom $\text{V}_g$ and the apex $\mathbf{x}_p$ if $e_g = \partial\Omega_p \cap \partial \Omega$ or of two polygonal pyramids consisted of the common bottom $\text{V}_g$ and the apices $\mathbf{x}_p$ and $\mathbf{x}_q$ if $e_g = \partial \Omega_p \cap \partial \Omega_q$.
\end{enumerate}
Now, we compute the representative gradient at the center of the triangle $e_a = \triangle_j = \triangle(\mathbf{x}_{v_j}, \mathbf{x}_{v_{j+1}}, \mathbf{x}_g)$ in Figure~\ref{fig:tess_full}:
\begin{align}\label{eq:rep_grad_tri}
(\alpha_a, \nabla u(\mathbf{x}_a)) = \argmin_{\substack{ (\alpha, \bm{\beta}) \in \mathbb{R} \times \mathbb{R}^3 \\ | \bm{\beta}| \leq 1 }} \sum_{\mathbf{x} \in \text{Q}_a} \left( w_a(\mathbf{x}) \big(\alpha +  \bm{\beta} \cdot (\mathbf{x} -\mathbf{x}_a) - u(\mathbf{x}) \big) \right)^2,
\end{align}
where the point $\mathbf{x}_a$ is the center of the mass in the triangle $\triangle_j$, $j=1, \ldots, J$ and $\text{Q}_a$ is the set of points:
\begin{align}
\text{Q}_a = \begin{cases}
	\{\mathbf{x}_{v_j}, \mathbf{x}_{v_{j+1}}, \mathbf{x}_g, \mathbf{x}_p\} \quad &\text{if} \: e_a \subset  e_g = \partial\Omega_p \cap \partial \Omega, \: a \in \mathcal{B} \\
	\{\mathbf{x}_{v_j}, \mathbf{x}_{v_{j+1}}, \mathbf{x}_g, \mathbf{x}_p, \mathbf{x}_q\} \quad &\text{if} \: e_a \subset  e_g = \partial\Omega_p \cap \partial \Omega_q, \: a \in \mathcal{F}
\end{cases}
\end{align}
Note that two steps,~\ref{A:step1} and~\ref{A:step2}, provide the values of $u$ used in~\eqref{eq:rep_grad_tri}. We emphasize that the representative gradient $ \nabla u(\mathbf{x}_a)$ in~\eqref{eq:rep_grad_tri} depends on the choice of the boundary conditions because the values on the vertices $u(\mathbf{x}_v)$ in~\eqref{eq:ver_val} are computed by the gradient $\nabla u_{\text{S}_p}(\mathbf{x}_p)$~\eqref{eq:grad} depending on the boundary conditions.

In the case of using the Dirichlet boundary condition, we define the representative gradient at the center of the cell $\Omega_p$:
\begin{equation}\label{eq:ABG}
\mathcal{D} u_{p} =
\displaystyle \frac{ \sum_{a \in \mathcal{F}_p \cup \mathcal{B}_p }\omega_p(\mathbf{x}_a) \nabla u(\mathbf{x}_a) }{ \sum_{a \in \mathcal{F}_p \cup \mathcal{B}_p } \omega_p(\mathbf{x}_a)},
\end{equation}
which is called the average-based gradient in~\cite{ref:HMFMB19}. When we compute the test cases using the Dirichlet boundary condition in Section~\ref{sec:NumEx}, the gradient above is applied in the algorithm~\eqref{eq:deffered_algo_DBC}. In contrast to~\eqref{eq:ABG}, in the proposed method, we differently define the representative gradient at the center of the cell $\Omega_p$:
\begin{equation}\label{eq:ABG_Soner}
\mathcal{D} u_{p} =
\displaystyle \frac{ \sum_{a \in \mathcal{F}_p \cup \mathcal{B}_p^{\nu^{+}} }\omega_p(\mathbf{x}_a) \nabla u(\mathbf{x}_a) }{ \sum_{a \in \mathcal{F}_p \cup \mathcal{B}_p^{\nu^{+}} } \omega_p(\mathbf{x}_a)}.
\end{equation}
It is used in the spatial discretization~\eqref{eq:algo_int} and~\eqref{eq:algo_bdr} or~\eqref{eq:deffered_algo_int} and~\eqref{eq:deffered_algo_bdr}. Note that there are no terms on inflow boundary triangles because it would violate the Soner boundary condition.

To complete the explanation of the notations, we define the center of the face or cell. From an internal face $e_g = \partial \Omega_p \cap \partial \Omega_q$ or a boundary face $e_g = \partial\Omega_p \cap \partial \Omega$, denoting the ordered vertices of the face as $V_g$, we define $\mathcal{C}$ as the collection of convex hull $H$ of three points; the center of mass of all vertices of the face and two consecutive points in $V_g$. Then, the center of the face is computed by
\begin{align}
\mathbf{x}_g = \frac{\sum_{H \in \mathcal{C}} |H| \bar{\mathbf{x}}_{H}}{\sum_{H \in \mathcal{C}} |H|},
\end{align}
where $\bar{\mathbf{x}}_H$ is the center of mass of the convex hull $H$. Similar to the center of the face, we also compute the center of the cell $\Omega_p$. Denoting all faces of the cell as $\partial \Omega_p = \cup_{i=1}^{I} e_{g_i}$, where $e_{g_i}$ is a face of the cell, and the ordered vertices of the face as $V_{g_i}$, we define $\mathcal{C}_i$ as the set of convex hull $H$ of four points; the center of the face, the center of mass of all vertices on the cell, and two consecutive points in $V_{g_i}$. Then, the center of the cell is computed by
\begin{align}\label{eq:app_cell_ct}
\mathbf{x}_p = \frac{\sum_{i=1}^{I} \sum_{H \in \mathcal{C}_i} |H| \bar{\mathbf{x}}_{H} }{\sum_{i=1}^{I} \sum_{H \in \mathcal{C}_i} |H|},
\end{align}
where $\bar{\mathbf{x}}_H$ is the center of mass of the convex hull $H$ in $\mathcal{C}_i$.

\section{The $EOC$ for test cases\label{sec:app_B}}
%
%
%

\begin{table}
\makebox[\linewidth][c]{%
	\begin{minipage}{1.25\textwidth}
		\centering
		\begin{tabular}{cS[table-format=1]
				S[table-format=1.2]
				S[table-format=1.2]
				S[table-format=1.2]
				S[table-format=1.2]
				S[table-format=1.2]
				S[table-format=1.2]
				S[table-format=1.2]
				S[table-format=1.2]}
			\toprule
			& & \multicolumn{4}{c}{DBC} & \multicolumn{4}{c}{EKBC}
			\\
			\cmidrule(lr){3-6} \cmidrule(lr){7-10}
			\multicolumn{1}{c}{} & \multicolumn{1}{c}{M} & \multicolumn{1}{c}{$EOC_{E^{1,\mathcal{Z}}_{\text{M}}}$} & \multicolumn{1}{c}{$EOC_{E^{\infty,\mathcal{Z}}_{\text{M}}}$} & \multicolumn{1}{c}{$EOC_{E^{\text{v}}_{\text{M}}}$} & \multicolumn{1}{c}{$EOC_{E^1_{\text{M}}}$} &
			\multicolumn{1}{c}{$EOC_{E^{1,\mathcal{Z}}_{\text{M}}}$} & \multicolumn{1}{c}{$EOC_{E^{\infty,\mathcal{Z}}_{\text{M}}}$} & \multicolumn{1}{c}{$EOC_{E^{\text{v}}_{\text{M}}}$} & \multicolumn{1}{c}{$EOC_{E^1_{\text{M}}}$} \\
			\cmidrule(lr){1-2} \cmidrule(lr){3-6} \cmidrule(lr){7-10}
			\multirow{3}{*}{\ref{TS}} & 1 & 2.16 & 2.38 & 3.54 & 2.13 & 2.27 & 2.40 & 3.51 & 2.04 \\ 
			& 2 & 1.97 & 1.88 & 2.66 & 2.03 & 1.98 & 1.87 & 2.72 & 1.70 \\ 
			& 3 & 2.16 & 2.08 & 2.53 & 2.14 & 2.15 & 2.08 & 2.61 & 1.57 \\ 
			\cmidrule(lr){1-2} \cmidrule(lr){3-6} \cmidrule(lr){7-10}
			\multirow{3}{*}{\ref{RS}} & 1 & 2.28 & 2.38 & 3.09 & 2.11 & 2.46 & 2.09 & 3.12 & 2.25 \\ 
			& 2 & 1.98 & 1.73 & 2.37 & 2.06 & 2.07 & 1.98 & 2.49 & 2.02 \\
			& 3 & 2.24 & 2.07 & 2.61 & 2.26 & 2.28 & 2.12 & 2.64 & 2.05 \\
			\cmidrule(lr){1-2} \cmidrule(lr){3-6} \cmidrule(lr){7-10}
			\multirow{3}{*}{\ref{ES}} & 1 & 2.41 & 3.03 & 2.32 & 2.28 & 2.41 & 3.03 & 2.31 & 2.30 \\
			& 2 & 2.24 & 2.24 & 2.27 & 2.03 & 2.24 & 2.24 & 2.27 & 2.03 \\
			& 3 & 2.58 & 2.40 & 2.69 & 2.24 & 2.58 & 2.40 & 2.69 & 2.24 \\
			\cmidrule(lr){1-2} \cmidrule(lr){3-6} \cmidrule(lr){7-10}
			\multirow{3}{*}{\ref{SS}} & 1 & 2.42 & 2.30 & 2.43 & 2.12 & 2.33 & 2.27 & 2.41 & 2.02 \\
			& 2 & 2.24 & 2.16 & 2.18 & 2.09 & 2.21 & 2.17 & 2.18 & 1.86 \\
			& 3 & 2.43 & 2.21 & 2.44 & 2.34 & 2.43 & 2.21 & 2.43 & 1.80 \\
			\cmidrule(lr){1-2} \cmidrule(lr){3-6} \cmidrule(lr){7-10}
			\multirow{3}{*}{\ref{TC}} & 1 & 1.14 & 0.77 & 1.88 & 1.62 & 1.15 & 0.81 & 1.75 & 1.79 \\
			& 2 & 1.25 & 0.66 & 1.95 & 1.66 & 1.24 & 0.66 & 1.96 & 1.65 \\
			& 3 & 1.38 & 0.68 & 2.20 & 1.73 & 1.38 & 0.67 & 2.21 & 1.60 \\
			\cmidrule(lr){1-2} \cmidrule(lr){3-6} \cmidrule(lr){7-10}
			\multirow{3}{*}{\ref{RC}} & 1 & 0.80 & 0.30 & 1.19 & 1.46 & 1.02 & 0.66 & 1.70 & 1.49 \\
			& 2 & 1.10 & 0.68 & 1.43 & 1.56 & 1.16 & 0.78 & 1.70 & 1.49 \\
			& 3 & 1.28 & 0.70 & 1.61 & 1.71 & 1.34 & 0.78 & 1.70 & 1.56 \\
			\cmidrule(lr){1-2} \cmidrule(lr){3-6} \cmidrule(lr){7-10}
			\multirow{3}{*}{\ref{EC}} & 1 & 1.29 & 1.00 & 1.17 & 1.65 & 1.29 & 1.00 & 1.16 & 1.68 \\
			& 2 & 1.14 & 1.04 & 1.08 & 1.39 & 1.14 & 1.04 & 1.08 & 1.39 \\
			& 3 & 1.20 & 1.10 & 1.15 & 1.44 & 1.20 & 1.10 & 1.15 & 1.44 \\
			\cmidrule(lr){1-2} \cmidrule(lr){3-6} \cmidrule(lr){7-10}
			\multirow{3}{*}{\ref{SC}} & 1 & 1.65 & 0.90 & 1.71 & 1.53 & 1.69 & 0.99 & 1.74 & 2.09 \\
			& 2 & 1.96 & 0.66 & 1.92 & 1.86 & 1.96 & 0.66 & 1.92 & 1.93 \\
			& 3 & 2.32 & 1.24 & 2.28 & 2.30 & 2.32 & 1.24 & 2.28 & 2.24 \\
			\bottomrule
		\end{tabular}
\end{minipage}}
\vspace{0.5em}
\begin{center}
	\caption{The $EOC$ for all test cases in Table~\ref{tab:test_case}, used in Section~\ref{sec:comp_EDBC}. The values are the slopes of the corresponding log-log scaled graphs in Figures~\ref{fig:conv_order_S} and~\ref{fig:conv_order_C}.\label{tab:comp_EDBC}}
\end{center}
\end{table}

For test cases in Section~\ref{sec:comp_EDBC}, we provide a supplementary table to check the $EOC$~\eqref{eq:EOC}. The values are the slopes of the corresponding graphs in Figures~\ref{fig:conv_order_S} and~\ref{fig:conv_order_C}. Table for Section~\ref{sec:comp_other}

\section{\nn{Exact solution for the expanding cube}\label{sec:app_EC}}

For the pure expansion the velocity is $\mathbf{v}=\delta\nabla u/|\nabla u|$
with $\delta>0$, so that~\eqref{eq:levelset} becomes the Hamilton--Jacobi
equation
\begin{equation}\label{eq:HJ_exp}
	u_t + \delta |\nabla u| = 0 \quad \text{in } \mathbb{R}^3\times(0,T],
	\qquad u(\cdot,0)=u^0 ,
\end{equation}
whose Hamiltonian is convex and whose unique viscosity solution is given by the
Hopf--Lax formula~\cite{ref:E98}
\begin{equation}\label{eq:hopflax}
	u(\mathbf{x},t) = \min_{|\mathbf{y}-\mathbf{x}|\le\delta t} u^0(\mathbf{y}) .
\end{equation}
In words, the value of $u$ at $\mathbf{x}$ at the time $t$ is the smallest value that the
initial condition takes on the ball of radius $\delta t$ centred at $\mathbf{x}$.

For the cube this minimum can be written down. Let $s:=\delta t$, let
$\xi_1\ge\xi_2\ge\xi_3$ be the decreasing rearrangement of $|x_i-c_i|$,
$i=1,2,3$, which loses nothing because the cube is symmetric with respect to its
centre and its three mid-planes, and let
\begin{equation}\label{eq:Dk}
	D_k := \!\!\sum_{1\le i<j\le k}\!\! (\xi_i-\xi_j)^2 , \qquad k=1,2,3 ,
\end{equation}
so that $D_1=0$, $D_2=(\xi_1-\xi_2)^2$ and $D_3-D_2=(\xi_1-\xi_3)^2+(\xi_2-\xi_3)^2$.
The point at which the minimum in~\eqref{eq:hopflax} is attained lies on a face,
on an edge or at a vertex of a cube concentric with the initial one, and which of
the three occurs is decided by $\xi_1,\xi_2,\xi_3$ and $s$ alone. Writing
$u^0_{\text{cub}}$ for $u^0_{\text{cub}}(\mathbf{x};\mathbf{c},r)$, the solution is
\begin{equation}\label{eq:exact_EC}
	u_e(\mathbf{x},t) =
	\begin{cases}
		u^0_{\text{cub}} - s,
		& \text{if } u^0_{\text{cub}} \ge s, \\[6pt]
		\xi_1 - s - r,
		& \text{else if } D_2 \ge s^2, \\[6pt]
		\dfrac{\xi_1+\xi_2-\sqrt{\,2s^2-D_2}}{2} - r,
		& \text{else if } D_3-D_2 \ge s^2, \\[10pt]
		-r,
		& \text{else if } |\mathbf{x}-\mathbf{c}| \le s, \\[6pt]
		\dfrac{\xi_1+\xi_2+\xi_3-\sqrt{\,3s^2-D_3}}{3} - r,
		& \text{else,}
	\end{cases}
\end{equation}
the cases being tested in the order in which they are written. They are, in that
order: the expansion has not yet reached $\mathbf{x}$, and the value is the
initial one lowered by $s$; the minimum is attained on a face; on an edge;
$\mathbf{x}$ lies in the plateau, where the ball already contains the centre of
the cube and the value has reached $\min u^0_{\text{cub}}=-r$; and the minimum is
attained at a vertex.



\end{appendices}


\bibliography{reference}

\end{document}